\documentclass[10pt,a4]{article}

\usepackage{graphicx}
\usepackage{latexsym}
\usepackage{amsmath}
\usepackage{amssymb}
\usepackage{amsthm}

\newcommand{\R}{{\mathbb R}}
\newcommand{\D}{{\mathbb D}}
\newcommand{\var}{\varepsilon}
\newcommand{\ve}{\varepsilon}

\renewcommand{\epsilon}{\varepsilon}

\renewcommand{\theequation}{\arabic{section}.\arabic{equation}}
\newtheorem{thm}{Theorem}[section]
\theoremstyle{definition}

\theoremstyle{Remark}

\theoremstyle{plain}
\newtheorem{lem}[thm]{Lemma}
\theoremstyle{plain}

\newtheorem{proposition}[thm]{Proposition}
\newtheorem{Remark}[thm]{Remark}

\def\bb{\begin{equation}}
\def\ee{\end{equation}}
\def\ba{\begin{equation*}}
\def\ea{\end{equation*}}
\def\bs{\begin{split}}
\def\es{\end{split}}

\newcommand{\up}{\widetilde}
\newcommand{\strip}    {\quad \mbox{in}\ \mathfrak{S},}

\newcommand{\intreal}{\int_\mathbb{R}}
\newcommand{\etaone}{\eta^{\varepsilon}_{\delta}}
\newcommand{\etathree}{\eta^{\varepsilon}_{3\delta}}
\newcommand{\intwx}[1]{\displaystyle{\intreal #1 w_x\mathrm{d}x}}
\newcommand{\intz}[1]{\displaystyle{\intreal #1 Z\mathrm{d}x}}

\begin{document}

\title
{Concentration on Surfaces for a Singularly Perturbed Neumann Problem in Three-Dimensional Domains}

\author
{
 Ying Guo
 \\
 College of mathematics and computational sciences,
 \\
  Shenzhen University,
 Shenzhen,  P. R. China, 518060.
 \\
 \\
 Jun Yang
\\
 College of mathematics and computational sciences,
 \\
  Shenzhen University,
 Shenzhen,  P. R. China, 518060.  Email: jyang@szu.edu.cn
}
\date{}

\maketitle

\begin{abstract}
 We consider the following singularly perturbed elliptic problem
$$
\varepsilon^2\triangle\tilde{u}-\tilde{u}+\tilde{u}^p=0, \  \tilde{u}>0\quad \mbox{in} \
\Omega,\ \ \
\frac{\partial\tilde{u}}{\partial \mathbf{n}}=0 \quad \mbox{on}\
\partial\Omega,
$$
where $\Omega$ is a bounded domain in $\mathbb{R}^3$ with
smooth boundary, $\varepsilon$ is a small parameter, $\mathbf{n}$ denotes the inward
normal of $ \partial\Omega$ and the exponent $p>1$.
Let $\Gamma$ be a hypersurface intersecting
$\partial\Omega$ in the right angle along its boundary $\partial\Gamma$ and satisfying
a {\em non-degenerate condition}.
We establish  the existence of a
solution $u_\varepsilon$ concentrating along a surface $\tilde{\Gamma}$ close to $\Gamma$, exponentially small in $\varepsilon$ at
any positive distance from the surface $\tilde{\Gamma}$, provided $\varepsilon$ is small and away from certain
{\em critical numbers}. The concentrating  surface $\tilde{\Gamma}$ will collapse to $\Gamma$ as $\varepsilon\rightarrow 0$.

\end{abstract}

\section{Introduction}\label{section1}
\setcounter{equation}{0}

  We consider the following problem
\begin{align}\label{originalproblemone}
\ve^2\triangle \tilde{u}-\tilde{u}+\tilde{u}^p=0, \ \  \tilde{u}>0\quad \mbox{in} \ \Omega
\quad \mbox{and}\quad
\frac{\partial\tilde{u}}{\partial \mathbf{n}}=0  \ \mbox{on}\ \partial\Omega,
\end{align}
where $\Omega$ is a bounded domain in $\mathbb{R}^N$ with smooth boundary, $\ve$ is a
small parameter, $\mathbf{n}$ denotes the inward normal of $\partial \Omega$
and the exponent $p>1$.

\subsection{Background and assumptions}\label{subsection1.1}
Problem (\ref{originalproblemone})  is  known as a stationary equation
of Keller-Segel system in chemotaxis \cite{LNT}. It can also be
viewed as a limiting stationary equation of Gierer-Meinhardt system
in biological pattern formation \cite{GM}. Problem
(\ref{originalproblemone}) has been studied extensively in recent
years. See \cite{Ni2} for backgrounds and  references.

In the pioneering papers, under the condition that $p$ is subcritical, i.e.,
$1<p<\frac{N+2}{N-2}$ when $N\geq 3$ and $ 1<p<+\infty$ when $N=2$,
Lin, Ni and Takagi\cite{LNT},  Ni and Takagi\cite{NT}-\cite{NT1},
established the existence of a least-energy solution
$U_{\ve}$ of ($\ref{originalproblemone}$)
and showed that, for $\ve$ sufficiently small, $U_{\ve}$ has only
one local maximum point $P_{\ve}\in\partial\Omega$. Moreover,
$H(P_{\ve})\rightarrow\displaystyle{\max_{P\in\partial\Omega}}H(P)$
as $\ve\rightarrow 0$, where $H(\cdot)$ is the mean curvature of
$\partial\Omega$.
Such solutions are called boundary spike-layers.

 Since then, many papers investigated further the solutions
 of ($\ref{originalproblemone}$) concentrating
 at one or multiple points of
 $\overline{\Omega}$. (These solutions are called {\bf spike-layers}.)
 A general principle is that
 the location of interior spikes is determined by the distance
 function from the boundary. We refer the reader to the
 articles \cite{BF},  \cite{dy2}, \cite{dPFW1},  \cite{GPW},   \cite{Wei2}, and references
 therein. On the other hand, boundary spikes are related to the
 mean curvature of $\partial\Omega$. This aspect is discussed in
 the papers \cite{BDS}, \cite{DY}, \cite{dPFW}, \cite{GWW}, \cite{Li}, \cite{Wei1}, and references therein.
 A good review of the subject up to 2004  can be found in  \cite{Ni2}.

The question of constructing higher-dimensional
 concentration sets has been investigated only in recent years.
 It has been conjectured in \cite{Ni2} that for any $1\leq k\leq N-1$,
 problem ($\ref{originalproblemone}$) has a solution $U_\ve$ which
 concentrates on a $k$-dimensional subset of $\overline{\Omega}$. We
 mention some results that support such a conjecture.

  In \cite{MalMon} and \cite{MalMon1}, Malchiodi and Montenegro proved
 that for $N \geq 2$,  there exists a sequence of numbers $\ve_\ell\rightarrow 0$ such
 that problem ($\ref{originalproblemone}$) has a solution $U_{\ve_{\ell}}$
 which  concentrates at the whole boundary  $\partial\Omega$ (or any component of
 $\partial\Omega$).  In \cite{Mal1, Mal2}, Malchiodi  showed the concentration
phenomena for (\ref{originalproblemone})
along a closed non-degenerate geodesic(i.e. codimension 2) of $\partial \Omega$ in
three-dimensional smooth bounded domain $\Omega$.
Mahmoudi and Malchiodi in \cite{MM} proved a full general concentration of
solutions along $k$-dimensional $(1\leq k\leq N-1)$ non-degenerate
minimal submanifolds of the boundary for $N\geq 3$ and
$1<p<\frac{N-k+2}{N-k-2}$.

In the above papers \cite{MM}-\cite{MalMon1}, the higher dimensional concentration
set is {\it on} the boundary. A natural question is if there are solutions with
high dimensional concentration set {\it inside} the domain.
For two dimensional case, the authors \cite{weiyang1, weiyang2} consider  problem (\ref{originalproblemone})
with solutions concentrating on
curves near a nondegenerate line $\Gamma'$ connecting the boundary of $\Omega$ at right angle.
The meaning of nondegeneracy of the line $\Gamma'$ can be defined similarly as
the nondegeneracy of the hypersurface $\Gamma$ in the sequel.
The reader can also refer to the survey paper by J. Wei \cite{wei}.

The main objective of the present paper is to extend the result in \cite{weiyang1} and show the existence
of concentration on interior hypersurfaces touching the boundary of $\Omega\subset\R^3$.
We make the following assumptions and notations.
The reader can also refer to the references \cite{delPinoKowalczykWei2,Sakamoto},
as well as \cite{docarmo} for some basic geometric results.
\\
{\textbf{(A1):}} Our candidate hypersurface $\Gamma \in \Omega$ is a {\it minimal} surface
that satisfies the following assumptions:
$\Gamma$ is smooth and embedded in $\R^3$,  and intersects $\partial\Omega$ in the right angle along its
boundary $\partial\Gamma=\bar{\Gamma}\cap\partial\Omega$,
which is a simple close curve in $\R^3$.
Let $\Delta^{\Gamma}$ be the Laplace-Beltrami operator on $\Gamma$ and $k_1,\,k_2$
the principal curvatures of $\Gamma$.
As a submanifold of $\R^3$, then we define the norm of the second fundamental form of $\Gamma$
by
$$
|A_{\Gamma}|^2=k_1^2+k_2^2,
$$
and the mean curvature of $\Gamma$ by $k=k_1+k_2$(with scalar $2$).
Note that $k=0$ along $\Gamma$.
We recall again that ${\bf n}$ is the inward unit normal vector on $\partial\Omega$, and hence, it
locally is also the unit normal vector of the curve $\partial\Gamma$ 
because of the perpendicularity between $\bar\Gamma$ and $\partial\Omega$.
Note that we here assume that the orientation of $\partial\Gamma$ is induced from that of $\partial\Omega$.
Since a curve on the surface $\partial\Omega$ is a geodesic if and only if its normal
vector is parallel to the normal vector ${\bf n}$ of $\partial\Omega$.
Therefore, by defining
\begin{align}
I(y)\,=\,-\big<\frac{\partial {\bf n}}{\partial\nu}, \nu\big>,
\quad
y\in\partial\Gamma\subset\partial\Omega.
\end{align}
we know that $-I$ is the curvature of the geodesic on $\partial\Omega$ passing through $ y\in\partial\Gamma$
in the normal direction $\nu(y)$ of $\Gamma$.

\noindent {\textbf{(A2):}} Let $\tau$ be the normal of the curve $\partial\Gamma$, which
is also the restriction of ${\bf n}$ on $\partial\Gamma$, $I$ is denoted in (\ref{I}).
By defining an {\it eigenvalue problem}, which will play an important role in our considerations,
\begin{align}\label{Geometriceigenvalueproblem}
\bigtriangleup^{\Gamma}f+|A_{\Gamma}|^2f=\lambda f\quad \mbox{in }\Gamma,
\qquad
\partial f/\partial\tau+I f=0\quad \mbox{on }\partial\Gamma,
\end{align}
we say that $\Gamma$ is {\it non-degenerate} if the eigenvalue
problem (\ref{Geometriceigenvalueproblem}) does not have zero eigenvalues.
The reader can refer to \cite{Sakamoto} for some explanations and examples.

\subsection{The profile function $w$}\label{subsection1.2}
 Before stating the main result,  we introduce two functions
$w$ and $Z$. Let $w$ be the unique (even) solution of
\begin{align}\label{definitionofw}
w''-w+w^p=0\ \mbox{and $w>0$ in
$\mathbb{R}$},\ w'(0)=0,\ w(\pm\infty)=0.
\end{align}
 It is well known that the
associated linearized eigenvalue problem,
\begin{align}\label{definitionofZ}
h''-h+pw^{p-1}h=\lambda h\quad  \mbox{in}\ \mathbb{R},\qquad  \int_{\mathbb{R}} h^2=1,\quad  h \in H^1 (\mathbb{R})
\end{align}
possesses a unique positive eigenvalue $\lambda_0$ with  a unique  even and
positive eigenfunction $Z$.
This follows for instance from the analysis in \cite{NT1}. In fact, we have
\begin{eqnarray} \label{lambda0}
w(x)&=&C_p\left\{\, \exp\Bigl[\,\frac{(p-1)\,x}{2}\, \Bigr]
  \, +\, \exp\Bigl[\, \frac{-(p-1)\,x}{2}\,\Bigr]\, \right\}^{\frac{-2}{p-1}},
\\
 Z&=&\Bigl[\, {\intreal{w^{p+1}}}\, \Bigr]^{-\frac{1}{2}}\, w^{p+1}, \quad \lambda_0=\frac{1}{4}(p-1)(p+3).
\end{eqnarray}
It is easy to see that for $|x|\gg 1$
\begin{align}
\label{decayofwandZ1}
w(x)    &= C_p\,e^{\,-|x|}\, -\, \frac{2C_p}{p-1}\, e^{\,-p|x|}\, +\, O(e^{\,-(2p-1)|x|}),
\\
\label{decayofwandZ2}
w'(x)&=-C_p\,e^{\,-|x|}\, +\, \frac{2\,p\,C_p}{p-1}\, e^{\,-p|x|}\, +\, O(e^{\,-(2p-1)|x|}),
\\
\label{decayofwandZ3}
Z(x)    &={\tilde C}_p\, e^{\,-(p+1)|x|}\, -\, \frac{2(p+1){\tilde C}_p}{p-1}\, e^{\,-2p|x|}\, +\, O(e^{\,-(3p-1)|x|}),\quad\quad\
\end{align}
where
\begin{eqnarray*}
C_p\,=\,\Bigl[\frac{(p+1)}{2}\Bigr]^{\frac{1}{p-1}},
\qquad
{\tilde C}_p\,=\,\Bigl[\frac{(p+1)}{2}\Bigr]^{\frac{p+1}{p-1}}\Bigl[\,{\int_{\mathbb{R}}{w^{p+1}\,{\mathrm{d}}x}}\,\Bigr]^{-\frac{1}{2}}.
\end{eqnarray*}

\subsection{Main Theorem}\label{subsection1.3}
Our main theorem can be stated as the following:
\begin{thm}
\label{main}
Assume that the minimal hypersurface $\Gamma$ satisfies the nondegeneracy
condition (\ref{Geometriceigenvalueproblem}).
There exists a sequence of small parameters $\{\ve_{l}\}_{l}$ such that
problem (\ref{originalproblemone}) has a positive solution $u_{\ve}$,
still denoting $\ve_{l}$ by $\ve$, concentrating along a surface ${\tilde\Gamma}$ near $\Gamma$.
Near $\Gamma$,  $u_\ve$ takes the form
\begin{equation}\label{uep}
u_{\epsilon}(\tilde{y}) \, \approx\,
w\, \left (\,  \frac{\mbox{dist}(\tilde{y}, {\tilde\Gamma} )}{\ve} \right ),
\end{equation}
where $w$ denotes the unique positive solution of problem (\ref{definitionofw}).
Moreover, there exists a positive number $c_0$ such that $u_{\varepsilon}$ satisfies globally,
$$u_{\ve}(\tilde{y})
\,\leq\,
\exp\bigl[\ -c_0\ \ve^{-1} \mbox{dist}(\tilde{y}, {\tilde\Gamma} )\ \bigr],
$$
and the surface ${\tilde\Gamma}$ will collapse to $\Gamma$ as $\ve\rightarrow 0$.
\qed
\end{thm}

To explain in a  few words the difficulties we have encountered, let us assume for the moment
that $\Omega= \mathbb{R} \times \Gamma$ is an infinite strip. In terms of the stretched coordinates
$(s, z)=\ve^{-1}(r,y)$ the equation would look near the surface
approximately like
$$
v_{ss}+\Delta^{\Gamma_{\var}}v-v+v^p=0,\qquad  (s, z) \in \mathfrak{S}:=\mathbb{R} \times \Gamma_{\var}.
$$
 The effects of curvature and of  the boundary conditions are
here neglected. The linearization of this problem around the profile
$w(s)$ becomes
$$
\phi_{ss}+\Delta^{\Gamma_{\var}}\phi-\phi+pw^{p-1}\phi=0,\ (s,z)\in \mathfrak{S},
$$
Let $\rho_{i},\,\omega_{i},\,  i=1,2,\cdots$, denote the eigenvalues and eigenfunctions of $-\Delta^{\Gamma}$,
then functions of the form
\bb
\begin{split}\nonumber
\phi^1_{i}=w_s(s)\,\omega_{i}(\ve z),
\qquad
\phi^2_{i}=Z(s)\,\omega_{i}(\ve z),
\end{split}
\ee are eigenfunctions associated to eigenvalues respectively
$-\rho_{i}\ve^2\ \mbox{and}\ \lambda_0-\rho_{i}\ve^2$. Many of these numbers
are {\em small} and thus {\em ``near non-invertibility" } of the linear operator
occurs. These two effects, combined in principle orthogonally
because of the $L^2$-orthogonality of $Z$ and $w_s$, are actually
coupled through the smaller order terms neglected.

In
\cite{ACF,Kowalczyk1,PR}, related singularly perturbed problems,
 the translation effect $\phi^1_{i}$ have been successfully treated
through successive improvements of the approximation and fine
spectral analysis of the actual linearized operator.
  In \cite{MalMon,MalMon1} resonance phenomena similar to the
``$\phi^2_{i}$-effect'' has been faced  in the
Neumann problem involving whole boundary concentration.
In \cite{MM}-\cite{Mal2},
\cite{delPinoKowalczykWei1},
\cite{wangweiyang},
\cite{weiyang1}-\cite{weiyang2},
the concentration on a $k-$dimensional minimal submanifold,
involving both
$\phi^1_{i}$ and $\phi^2_{i}$ effects, has been treated via arbitrary high
order approximations.

To prove Theorem \ref{main}, not only the same difficulties as that in the two dimensional case in
\cite{weiyang1}-\cite{weiyang2} are encountered but
also more obstruction appears.
In other words, in the present paper we need to handle more delicate resonance phenomena
due to the existence of $\phi^2_{i}$-effect in higher dimension,
as well as the strong interaction between the concentration set and the boundary due to
the homogeneous boundary condition in (\ref{originalproblemone}).
Some words are in order to explain the methods to handle these difficulties.

As we have stated in the above, we first neglect the boundary condition.
By the suitable rescaling of the variables $(s, z)=\ve^{-1}(r, y)$ in a type of Fermi coordinates(cf. Lemma \ref{lemma2.2}),
we then try the inner approximate solution to the problem roughly in the form,
\begin{align}\label{firstapproximation}
v( s, z)=w\big(s - f(\var z)\big)\,+\, \ve e(\var z)Z\big(s - f(\var z)\big) \,+\, \tilde{\varphi}(s, z),
\end{align}
where $Z$ is defined in (\ref{definitionofZ}), $f$ and $e$ are left as parameters,
while $\tilde{\varphi}(s, z)$ is $L^2$-orthogonal for each $z$ both to $w_s(s - f(\var z))$ and $Z(s - f(\var z))$.
Solving first in $\tilde{\varphi}$ a natural projected problem,
the standard reduction procedure will implies that
the resolution of the full problem will be reduced to a nonlinear,
nonlocal second order system of differential equations in $(f, e)$.
The nonlinear system is not solvable due to the $\phi^{2}_{i}$ effect in the case of dimension $N=3$.
This shows that the approximation in (\ref{firstapproximation}) doesn't work well as $N$ becomes large.
Hence we must improve our approximation by the recurrence method used in \cite{mmm}, \cite{Mal2}-\cite{MalMon1}, \cite{wangweiyang}.
The principle is: the better the approximation,
higher the chances of a correct inversion of the full problem to obtain a contraction mapping formulation of the nonlinear,
nonlocal second order differential equations.
To do that, we try the following form as our new approximation, (see \cite{MalMon1})
\[
v( s, z)=w\big(s - f(\var z)\big) \,+\, \ve e(\var z)Z\big(s - f(\var z)\big) \,+\, \sum_{l=1}^{3}\ve^l\varphi_l( s, z).
\]
The aim of adding the term $\sum_{l=1}^{3}\ve^l\varphi_l( s, z)$ is to cancel the error term till order $O(\ve^4)$
such that our approximation is good enough.
After very tedious but necessary computations we find that such $\varphi_l$'s must satisfy some
differential equations(cf. (\ref{defitionsolution1}), (\ref{defitionsolution2}), (\ref{defitionsolution34})),
whose solvability relies on the parameter $f$.
So we need to improve our approximation further by subtle adjusting the location of the concentration set,
namely we take the following form of inner approximate solution, (see \cite{MM}, \cite{wangweiyang})
\begin{align}\label{secondapproximation}
{\mathcal V}=w(x)\,+\,\ve e(\var z)Z(x)\,+\,\sum_{l=1}^{3}\ve^l\varphi_l(x,\var z)
\quad
\mbox{with }
\quad
x = s - \sum_{l=0}^{2} \ve^l f_l(\var z).
\end{align}
To fulfill the objective we need to analyze the corresponding linear problem very carefully
and conduct lots of computations.
This was done in Section \ref{section3}.
Note that the authors in \cite{wangweiyang} also met the same situations  for the concentration phenomena on high
dimensional hypersurface of inhomogeneous Schr$\ddot{o}$diger equation on any $N$ dimensional space.

Second, we know that the inner approximation in general
does not satisfy the boundary condition in (\ref{originalproblemone}).
Note that, for the two dimensional case of (\ref{originalproblemone}),
by using the condition that
the limit location  of concentration set, say $\Gamma'$ as before,
connect the boundary $\partial\Omega$ orthogonally,
the authors \cite{weiyang1}-\cite{weiyang2}
used a type of local coordinates (previously used by M. Kowalczyk in \cite{Kowalczyk1}) in a neighborhood of $\Gamma'$ to stretch the boundary $\partial\Omega$ and decompose the interaction
between the concentration set and boundary $\partial\Omega$
in such a way that they can improve the approximation to satisfy the boundary condition up to order of $O(\ve^2)$.
Unfortunately, the coordinate system can not be extended to the present case.
Here, in a neighborhood of the hypersurface $\Gamma$,
we will set up another type of local coordinates in Lemma \ref{lemma2.1}(used by K. Sakamoto in \cite{Sakamoto}),
called {\bf modified Fermi Coordinates} in this paper for the convenience of notion,
such that we can add boundary corrections to the inner expansion in (\ref{secondapproximation})
and finally get a good local approximation.
The reader can refer to \ref{section4}
for more details on the construction of boundary correction layers.

In the rest paper we carry out the program outlined above,
which leads to the complete proof of Theorem \ref{main}
by the infinite dimensional reduction method introduced in \cite{delPinoKowalczykWei1}.
The organization of the paper is as follows.
In Section \ref{section2}, we introduce a local modified Fermi coordinates
and then set up local formulation of problem (\ref{originalproblemone}) by suitable rescaling.
Sections \ref{section3} and \ref{section4} are devoted to the construction of
 a good approximate solution to solve the problem up to order $O(\ve^4)$.
Indeed, we first ignore the boundary condition in (\ref{originalproblemone}) and construct
an inner approximation solution by the induction method in Section \ref{section3}.
In Section 4, we add the boundary correction layers to the inner expansion and get the
final approximations involving unknown functions $(f_{2}, e)$.
In Section \ref{section5}, a gluing procedure from \cite{delPinoKowalczykWei1} reduces the nonlinear problem  to a projected problem on
the infinite strip $\mathfrak{S}$, while in Section \ref{section6} and \ref{section7}, we show that the projected
problem has a unique solution for the pair of $(f_{2},e)$ in a chosen region.
The final step is
to adjust the parameters $f_{2}$ and $e$ which is equivalent to solving a nonlocal, nonlinear coupled second
order system of differential equations for the pair $(f_{2}, e)$ with boundary
conditions. This is done in Sections \ref{section8} and \ref{section9}.

Finally, we remark that our arguments can be easily extended to show the existence of interior concentration phenomena
approaching minimal non-degenerate hypersurfaces and possessing interaction with boundary
for problem (\ref{originalproblemone})
on bounded domain of general dimension,
see Remarks \ref{remarkinnerapproximation} and \ref{remarkboundaryapproximation}.
However, for the clearness of presentation,
we here only consider the case with concentration on hypersurfaces in three dimension.

\section{Local Formulation of Problems}\label{section2}
\setcounter{equation}{0}

\medskip

The main objective of this subsection is to set up a local suitable
coordinate system near the minimal hypersurface $\Gamma$
and then express problem (\ref{originalproblemone}) in local form.

\subsection{Geometric notions}\label{subsection3.1}
Let $\gamma_0(\cdot):{\mathbb{D}}\rightarrow\Gamma\subset\Omega$ be a smooth parameterization
of the surface $\Gamma$, where $\Gamma$ is the minimal surface in Section 1
and $\mathbb{D}:=\{y\in\mathbb{R}^2:|y|<1\}$ is the unit disk.
For further computational convenience,
we will use an isothermal representation $\gamma_{0}:\mathbb{D}\rightarrow\Gamma$.
Namely,  we use $\gamma_{0}$ that satisfies
\begin{align}
\Bigg|\frac{\partial\gamma_{0}}{\partial y^{1}}\Bigg|^{2}
=\Bigg|\frac{\partial\gamma_{0}}{\partial y^{2}}\Bigg|^{2}
={\hslash}^{2}(y),
\qquad
\Bigg<\frac{\partial\gamma_{0}}{\partial y^{1}},\,\frac{\partial\gamma_{0}}{\partial y^{2}}\Bigg>=0.
\end{align}
Therefore, the tangent vectors ${\partial\gamma_{0}}/{\partial y^{1}}$ and ${\partial\gamma_{0}}/{\partial y^{2}}$
have the same length ${\hslash}(y)>0$ and are mutually orthogonal.

Recall the notations given in subsection \ref{subsection1.1}. For our future setting-up of problem (\ref{originalproblemone}),
here we provide a type of local coordinates modified from the standard Fermi coordinates,
called {\bf modified Fermi coordinates} for convenience of notions.
This coordinate system was previously used by K. Sakamoto to describe the transition layer for Allen-Cahn equation
in \cite{Sakamoto}.
\begin{lem}\label{lemma2.1}
There exist constants $\delta>0,\, r_0>0$, which depend only on $\Gamma$ and $\partial\Omega$, and
a smooth diffemorphism
\begin{align}
\gamma(\cdot,\cdot):\, (-r_0,r_0)\times\mathbb{D}_{\delta}\rightarrow\Omega^{r_0}_{\delta},
\end{align}
where $\mathbb{D}_{\delta}$ and $\Omega^{r_0}_{\delta}$ are given in (\ref{coordinated}) and (\ref{fermicoordinate}) such that
\begin{enumerate}
\item
$\gamma(0,y)=\gamma_0(y)$ for $y\in{\mathbb{D}}$,
$\quad \gamma(r,y)\in\partial\Omega$ for $y\in\partial\mathbb{D},\, |r|<r_0$;

\item
 $\frac{\partial\gamma}{\partial r}(0,y)=\nu(y)$ for $y\in\mathbb{D}$;

\item
$\gamma(r,y)$ has the following expansion, as $r\rightarrow 0$
\begin{align}\label{coordinateexpression}
\gamma(r,y)=\gamma_0(y)+r\nu(y)+\frac{r^2}{2}q_{1}(y)+\frac{r^3}{6}q_{2}(y)+O(r^4),\quad y\in\mathbb{D},
\end{align}
where $q_{1}(y)$ and $q_{2}(y)$ are vector functions orthogonal to $\nu(y)$;

\item
If we write $\gamma$ as $\gamma(r,\vartheta,\rho)$ in terms of the coordinate $(r,\vartheta,\rho)$
in (\ref{boundarypolar}),
then the derivative along the inward unit normal vector ${\bf n}$ of $\partial\Omega$ is expressed as
\begin{align}
\frac{\partial}{\partial {\bf n}}=\frac{1}{\sqrt{ g^{33}}}\Big( g^{13}\frac{\partial}{\partial r}
+ g^{23}\frac{\partial}{\partial \vartheta}
+g^{33}\frac{\partial}{\partial\rho}\Big),
\end{align}
where at $(r,\vartheta,0)$, we have the expression
\begin{align}
\begin{aligned}
 g^{13}(r,\vartheta)&=-r\Big|\frac{\partial\gamma_0}{\partial\rho}\Big|^{-2}\Big<q_{1}, \frac{\partial\gamma_0}{\partial\rho}\Big>+O(r^2),
\\
 g^{23}(r,\vartheta)&=2r\Big|\frac{\partial\gamma_0}{\partial\vartheta}\Big|^{-2}
\Big|\frac{\partial\gamma_0}{\partial\rho}\Big|^{-2}\Big<\frac{\partial\gamma_0}{\partial\vartheta}, \frac{\partial\nu}{\partial\rho}\Big>+O(r^2),
\\
g^{33}(r,\vartheta)&=\Big|\frac{\partial\gamma_0}{\partial\rho}\Big|^{-2}
+2r\Big|\frac{\partial\gamma_0}{\partial\rho}\Big|^{-4}\Big<\frac{\partial\gamma_0}{\partial\rho}, \frac{\partial\nu}{\partial\rho}\Big>+O(r^2).
\end{aligned}
\end{align}
\end{enumerate}
\end{lem}

\proof
The proof can be found in \cite{Sakamoto}.
For the convenience of readers, we also provide the details here.
We extend $\gamma_0$ smoothly  to
\begin{align}\label{coordinated}
\mathbb{D}_{\delta}:=\{\,y\in\mathbb{R}^2:|y|<1+\delta\,\},
\end{align}
for some fixed constant $\delta>0$.
The extension is still denoted by $\gamma_0$ and its image by $\Gamma_{\delta}$.
Let $\nu(y)\in\mathbb{R}^3$ be the unit normal of $\Gamma_{\delta}$ at $\gamma_0(y)\in\Gamma_{\delta}$.
We now define a neighborhood $\Omega^{r_0}_{\delta}$ of $\Gamma_{\delta}$ by the Fermi coordinates
\begin{align}\label{fermicoordinate}
\Omega^{r_0}_{\delta}:=\big\{\tilde y\in\mathbb{R}^3:\,\tilde y=\gamma_0(y)+r\nu(y),\,
|r|<r_0,\,
y\in\mathbb{D}_{\delta}\big\}.
\end{align}
Now we choose $r_0>0$ in (\ref{fermicoordinate}) small enough so that
$$
(\Omega^{r_0}_{\delta}\cap\Omega)
\cap\Big\{\,\gamma_{0}(y)+r\nu(y):\, |r|<r_{0},\, |y|=1+\delta\,\Big\}=\emptyset\,.
$$
When we deal with the portion of $\partial\Omega$ in $\Omega^{r_0}_{\delta}$,
we use the coordinate $(\vartheta,\rho)\in\partial\mathbb{D}\times[0,\delta)\subset\mathbb{D}$
to parameterize a small region nearby $\partial\Gamma$ in $\Gamma$,
where the points $(\vartheta,0)$ is sent to the boundary of surface $\partial\Gamma$ by $\gamma_0$
and there holds
\begin{align}
\Big<\frac{\partial\gamma_0}{\partial\vartheta},\,\frac{\partial\gamma_0}{\partial\rho}\Big>=0\ \mbox{ for }\rho=0.
\end{align}

We define the Fermi coordinates in (\ref{fermicoordinate}) by
$\bar{\gamma}:(-r_{0}, r_{0})\times\mathbb{D}_{\delta}\longmapsto\mathbb{R}^{3}$, i.e.
$$
\bar{\gamma}(r,y)=\gamma_{0}(y)+r\nu(y).
$$
Note that this coordinate system is not suitable for the boundary of $\Omega$.
Whence the main objective of this proof is to make a suitable modification of ${\bar\gamma}$.

Denote by $S$ the preimage of $\Omega^{r_0}_{\delta}\cap\partial\Omega$:
$$
S=\bar{\gamma}^{-1}(\Omega^{r_0}_{\delta}\cap\partial\Omega).
$$
Since $\partial\Omega\perp\Gamma$ by \textbf{(A1)}, we have
\begin{align}\label{perps}
S\perp_{\partial\mathbb{D}}(\{0\}\times\mathbb{D}).
\end{align}
We also denote by $C$ the preimage of $\Omega^{r_0}_{\delta}\cap\Omega$:
$$
C=\bar{\gamma}^{-1}(\Omega^{r_0}_{\delta}\cap\Omega),
$$
and by $C(r)$ the $r-$slice of $C$:
$$
C(r)=\{\,y\in\mathbb{D}_{\delta}:\,(r,y)\in C\,\}
\quad
\mbox{for }|r|<r_{0}.
$$
For later use, we set
$$
C_\delta=\bar{\gamma}^{-1}(\Omega^{r_0}_{\delta}),
\qquad
C_\delta(r)=\{\,y\in\mathbb{D}_{\delta}:\,(r,y)\in C_\delta\,\}
\quad
\mbox{for }|r|<r_{0},
$$

Since $\partial\Omega$ and $\gamma_{0}$ are smooth, $C(r)$ is a smooth domain, diffeomorphic to $C(0)=\mathbb{D}$.
Therefore, there exists a family of smooth diffeomorphisms
\begin{align*}
\widetilde{Y}(r,\cdot):\,\mathbb{D}\longmapsto C(r)
\end{align*}
parametrized smoothly by $r\in(-r_{0}, r_{0})$.
Thanks to (\ref{perps}), we can choose $\widetilde{Y}$ so that
\begin{align*}
\widetilde{Y}(0,y)=y, \quad
\frac{\partial \widetilde{Y}}{\partial r}(0,y)=0,\quad y\in\mathbb{D}.
\end{align*}
Furthermore, we make a smooth extension
\begin{align*}
Y(r,\cdot):\,\mathbb{D}_{\delta}\longmapsto C_{\delta}(r)
\end{align*}
such that
\begin{align}\label{defineY}
\begin{aligned}
Y(0,y)\,=\,y,
\qquad
\frac{\partial Y}{\partial r}(0,y)\,=\,0,\quad y\in\mathbb{D}_{\delta},
\\
Y(r,y)\,=\,{\widetilde Y}(r, y),\quad (r, y)\in(-r_0,r_0)\times\mathbb{D}.
\end{aligned}
\end{align}

Let us now define the desired diffeomorphism $\gamma$ by, called {\bf modified Fermi coordinates}
\begin{align}\label{definer}
\gamma(r,y)=\bar{\gamma}(r,Y(r,y))=\gamma_{0}(Y(r,y))+r\nu(Y(r,y))
\end{align}
for $(r,y)\in(-r_{0}, r_{0})\times\mathbb{D}_{\delta}.$
It is now straightforward to verify that $\gamma$
in (\ref{definer}) satisfies Lemma \ref{lemma2.1}(1). By elementary computations and (\ref{defineY}), we find that
\begin{align}
\begin{aligned}\label{vpq}
\frac{\partial\gamma}{\partial r}(0,y)&=\nu(y)=\nu(Y(0,y)),
\\
q_{1}(y)=\frac{\partial^{2}\gamma}{\partial r^{2}}(0,y)
&=\sum^{2}_{j=1}\frac{\partial\gamma_{0}}{\partial Y^{j}}\frac{\partial^{2}Y^{j}}{\partial r^{2}}(0,y)\perp\nu(y),
\\
q_{2}(y)=\frac{\partial^{3}\gamma}{\partial r^{3}}(0,y)
&=\sum^{2}_{j=1}\Big(\frac{\partial\gamma_{0}}{\partial Y^{j}}\frac{\partial^{3}Y^{j}}{\partial r^{3}}(0,y)
+3\frac{\partial\nu}{\partial Y^{j}}\frac{\partial^{2}Y^{j}}{\partial r^{2}}(0,y)\Big)\perp\nu(y),
\end{aligned}
\end{align}
which will show the validity of the statements (2) and (3).

To prove (4) of Lemma \ref{lemma2.1},
we use the coordinates $(\vartheta, \rho)$ introduced in the above.
Recall that $(\vartheta, \rho=0)$ parametrizes $\partial\mathbb{D}$ and $\rho$ is chosen so that
$$
\Big\langle\frac{\partial\gamma_{0}}{\partial\vartheta},
\frac{\partial\gamma_{0}}{\partial\rho}\Big\rangle=0\,\,\quad \mbox{at}\,\,\rho=0.
$$
For $y\in\mathbb{D}$ near $\partial\mathbb{D}$, we express $\gamma(r,y)$ by
\begin{align}\label{boundarypolar}
\gamma(r,y)=\gamma(r,\vartheta,\rho).
\end{align}
We also denote by $\mathbf{n}(r,\vartheta)$ the unit inward normal vector of $\partial\Omega$ at $\gamma(r,\vartheta,0).$
Note that at $\rho=0$(i.e. $\mbox{on}\,\,\partial\Omega\cap\Omega^{r_0}_{\delta})$, vectors
$$
\frac{\partial\gamma}{\partial r}\,,\,\frac{\partial\gamma}{\partial\vartheta}\,,
\,\frac{\partial\gamma}{\partial\rho}\in \mathbb{R}^{3},
$$
constitute a basis for $\mathbb{R}^{3}$. Hence $\mathbf{n}(r,\vartheta)$ is expressed as
\begin{align}\label{definen}
\mathbf{n}=a\frac{\partial\gamma}{\partial r}+b\frac{\partial\gamma}{\partial\vartheta}
+c\frac{\partial\gamma}{\partial\rho}\,\quad \mbox{at}\,\,\rho=0,
\end{align}
where $c>0$.
Since $\Big\{\frac{\partial\gamma}{\partial r}\,,\,\frac{\partial\gamma}{\partial\vartheta}\Big\}$
spans
the tangent space of $\partial\Omega$ at $\gamma(r,\vartheta,0)$, we have
\begin{align}\label{definperp}
\Big\langle\frac{\partial\gamma}{\partial r}, \mathbf{n}\Big\rangle=0,
\qquad
\Big\langle\frac{\partial\gamma}{\partial\vartheta}, \mathbf{n}\Big\rangle=0,
\qquad
\langle \mathbf{n},\mathbf{n}\rangle=1\,\quad \mbox{at}\,\,\rho=0.
\end{align}
From (\ref{definen}) and (\ref{definperp}), we easily obtain
$$
a=\frac{\tilde{g}^{13}}{\sqrt{\tilde{g}^{33}}},
\qquad
b=\frac{\tilde{g}^{23}}{\sqrt{\tilde{g}^{33}}},
\qquad
c=\frac{\tilde{g}^{33}}{\sqrt{\tilde{g}^{33}}}=\sqrt{\tilde{g}^{33}},
$$
where the formulas depend on the inverse of metric matrix at $\rho=0$,
i.e. $(\tilde{g}^{ij})=(\tilde{g}_{ij})^{-1}$.
More precisely, the metric coefficients $g_{ij}$'s are given by,
\begin{align*}
\tilde{g}_{11}&=\Big\langle\frac{\partial\gamma}{\partial r}, \frac{\partial\gamma}{\partial r}\Big\rangle,
\qquad
\tilde{g}_{12}=\tilde{g}_{21}=\Big\langle\frac{\partial\gamma}{\partial r}, \frac{\partial\gamma}{\partial\vartheta}\Big\rangle,
\\
\tilde{g}_{22}&=\Big\langle\frac{\partial\gamma}{\partial\vartheta}, \frac{\partial\gamma}{\partial \vartheta}\Big\rangle,
\qquad
\tilde{g}_{13}=\tilde{g}_{31}=\Big\langle\frac{\partial\gamma}{\partial r}, \frac{\partial\gamma}{\partial\rho}\Big\rangle,
\\
\tilde{g}_{33}&=\Big\langle\frac{\partial\gamma}{\partial\rho}, \frac{\partial\gamma}{\partial\rho}\Big\rangle,
\qquad
\tilde{g}_{23}=\tilde{g}_{32}=\Big\langle\frac{\partial\gamma}{\partial\vartheta}, \frac{\partial\gamma}{\partial\rho}\Big\rangle.
\end{align*}
Therefore, $\frac{\partial}{\partial \mathbf{n}}$ is given by
$$
\frac{\partial}{\partial \mathbf{n}}=\frac{1}{\sqrt{\tilde{g}^{33}}}\Big(\tilde{g}^{13}\frac{\partial}{\partial r}
+\tilde{g}^{23}\frac{\partial}{\partial\vartheta}+\tilde{g}^{33}\frac{\partial}{\partial\rho}\Big).
$$

Let us now take Taylor expansions of $\tilde{g}^{jk}(r, \vartheta, 0)$ in $r$ at $r=0.$ From the expansion of $\gamma(r, \vartheta, 0)$ in
Lemma \ref{lemma2.1}(3), we have, at $\rho=0$
\begin{align}\label{express}
\begin{aligned}
\frac{\partial\gamma}{\partial r}&=\nu(\vartheta,0)+rq_{1}(\vartheta,0)+\frac{r^{2}}{2}q_{2}(\vartheta,0)+O(r^{3}),
\\
\frac{\partial\gamma}{\partial \vartheta}&=\frac{\partial\gamma_{0}}{\partial\vartheta}(\vartheta,0)
+r\frac{\partial\nu}{\partial\vartheta}(\vartheta,0)
+\frac{r^{2}}{2}\frac{\partial q_{1}}{\partial\vartheta}(\vartheta,0)+O(r^{3}),
\\
\frac{\partial\gamma}{\partial \rho}&=\frac{\partial\gamma_{0}}{\partial\rho}(\vartheta,0)
+r\frac{\partial\nu}{\partial\rho}(\vartheta,0)
+\frac{r^{2}}{2}\frac{\partial q_{1}}{\partial\rho}(\vartheta,0)+O(r^{3}).
\end{aligned}
\end{align}
By using the orthogonalities $\frac{\partial\gamma_{0}}{\partial y}\perp\nu\,,\,
\frac{\partial\nu}{\partial y}\perp\nu\,,\,q_{1}\perp\nu$ and $\frac{\partial\gamma_{0}}{\partial\vartheta}
\perp\frac{\partial\gamma_{0}}{\partial\rho}$, together with (\ref{express}),
we find that, at $\rho=0$
\begin{align}
\begin{aligned}\label{defg}
g=
\left(
  \begin{array}{ccc}
    \tilde{g}_{11} & \tilde{g}_{12} & \tilde{g}_{13}
    \\
    \tilde{g}_{21} & \tilde{g}_{22} & \tilde{g}_{23}
    \\
    \tilde{g}_{31} & \tilde{g}_{32} & \tilde{g}_{33}
    \\
  \end{array}
\right)
\,=\,&\,\left(
  \begin{array}{ccc}
    1 &                   0                                    & 0
    \\
    0 & \big|\frac{\partial\gamma_{0}}{\partial\vartheta}\big|^{2} & 0
    \\
    0 &                   0               & \big|\frac{\partial\gamma_{0}}{\partial\rho}\big|^{2}
    \\
  \end{array}
\right)
\\
&\,+\,r
\left(
  \begin{array}{ccc}
  0 & \langle q_{1},  \frac{\partial\gamma_{0}}{\partial\vartheta}\rangle
    & \langle q_{1},  \frac{\partial\gamma_{0}}{\partial\rho}\rangle
  \\
  \langle q_{1},  \frac{\partial\gamma_{0}}{\partial\vartheta}\rangle
  & -2\tilde{L} & -2\tilde{M}
  \\
  \langle q_{1},  \frac{\partial\gamma_{0}}{\partial\rho}\rangle
  & -2\tilde{M} & -2\tilde{N}
  \\
  \end{array}
\right)
+O(r^{2}),
\end{aligned}
\end{align}
where
\begin{align*}
-\tilde{L}&=\Big\langle \frac{\partial\gamma_{0}}{\partial\vartheta}, \frac{\partial\nu}{\partial\vartheta}\Big\rangle, \,\,\quad
-\tilde{N}=\Big\langle \frac{\partial\gamma_{0}}{\partial\rho}, \frac{\partial\nu}{\partial\rho}\Big\rangle,
\\
-2\tilde{M}&=\Big\langle \frac{\partial\gamma_{0}}{\partial\vartheta}, \frac{\partial\nu}{\partial\rho}\Big\rangle
+\Big\langle \frac{\partial\gamma_{0}}{\partial\vartheta}, \frac{\partial\nu}{\partial\rho}\Big\rangle.
\end{align*}
Therefore, we obtain
\begin{align*}
\left(
  \begin{array}{ccc}
    \tilde{g}^{11} & \tilde{g}^{12} & \tilde{g}^{13}
    \\
    \tilde{g}^{21} & \tilde{g}^{22} & \tilde{g}^{23}
    \\
    \tilde{g}^{31} & \tilde{g}^{32} & \tilde{g}^{33}\\
  \end{array}
\right)
\,=\,&\,
\left(
  \begin{array}{ccccc}
    1
    &\vdots &
    0
    &\vdots &
    0
    \\
    0
    &\vdots &
    \big|\frac{\partial\gamma_{0}}{\partial\vartheta}\big|^{-2}
    &\vdots &
    0
    \\
    0
    &\vdots &
    0
    &\vdots &
    \big|\frac{\partial\gamma_{0}}{\partial\rho}\big|^{-2}\\
  \end{array}
\right)
\\
&+r
\left(
  \begin{array}{ccccc}
  0
  &\vdots&
  -\big|\frac{\partial\gamma_{0}}{\partial\vartheta}\big|^{-2}\langle q_{1},  \frac{\partial\gamma_{0}}{\partial\vartheta}\rangle
  &\vdots &
  -\big|\frac{\partial\gamma_{0}}{\partial\rho}\big|^{-2}\langle q_{1},  \frac{\partial\gamma_{0}}{\partial\rho}\rangle
  \\
  -\big|\frac{\partial\gamma_{0}}{\partial\vartheta}\big|^{-2}\langle q_{1},  \frac{\partial\gamma_{0}}{\partial\vartheta}\rangle
  &\vdots&
  2\big|\frac{\partial\gamma_{0}}{\partial\vartheta}\big|^{-4}\tilde{L}
  &\vdots&
  2\big|\frac{\partial\gamma_{0}}{\partial\vartheta}\big|^{-2}\big|\frac{\partial\gamma_{0}}{\partial\rho}\big|^{-2}\tilde{M}
  \\
  -\big|\frac{\partial\gamma_{0}}{\partial\rho}\big|^{-2}\langle q_{1},  \frac{\partial\gamma_{0}}{\partial\rho}\rangle
  &\vdots&
  2\big|\frac{\partial\gamma_{0}}{\partial\vartheta}\big|^{-2}\big|\frac{\partial\gamma_{0}}{\partial\rho}\big|^{-2}\tilde{M}
  &\vdots&
  2\big|\frac{\partial\gamma_{0}}{\partial\rho}\big|^{-4}\tilde{N}\\
  \end{array}
\right)
\\
\,&\,+O(r^{2}).
\end{align*}
This completes the proof of Lemma \ref{lemma2.1}.
\qed

We can use the modified Fermi coordinate system in (\ref{definer}) to express the Laplacian operator
$\Delta$ in (\ref{originalproblemone}) locally in terms of $(r,y)\in(-r_0,r_0)\times\D$.
Since $\tilde{y}=\gamma(r,y)$, the standard metric in
$\Omega^{r_0}_{\delta}\subset\R^{3}$ is pulled back to $g_{j k}(r,y)$:
\begin{align}
g_{j k}(r,y)=\Big\langle\frac{\partial\gamma}{\partial y^{j}},\,\frac{\partial\gamma}{\partial y^{k}}\Big\rangle, \quad  j,k=0,1,2,
\end{align}
where $y^{0}$ stands for $r$.
Hence, the {\it Laplace-Beltrami operator} and {\it gradient operator }
are defined in local coordinates by
\begin{align}\label{laplace}
\Delta\equiv\frac{1}{\sqrt{\mbox{det} g}}\,\partial_i\Bigl(\sqrt{\mbox{det} g}\, g^{ij}\partial_j\Bigr),
\quad
\nabla h\equiv g^{ij}\partial_i h\partial_j,
\end{align}
where $h$ is any smooth function and the coefficients $g^{ij}$ are
the entries of the inverse matrix of $g=(g_{ij})$.
There also hold similar expressions for
the {\it Laplace-Beltrami operator} and {\it gradient operator } on $\Gamma$.

For later use, we also recall  the Weyl's asymptotic formula, referring for example to \cite{chavel}, or to \cite{LiYau} and \cite{minakPleijel} for
further details. Let $\rho_i,\, \omega_i, i=1,2,\cdots,$ denote the eigenvalues and eigenfunctions
of $-\Delta^{\Gamma}$(ordered to
be non-decreasing in $i$ and counted with the multiplicity), then we have that
\begin{align}\label{weylonmanifold}
\rho_i\rightarrow \frac{C\,i^{2/(N-1)}}{\mbox{Vol}(\Gamma)}\quad\mbox{as }i\rightarrow\infty,
\end{align}
where $\mbox{Vol}(\Gamma)$ is the volume of $(\Gamma, {\tilde g})$ and ${\tilde g}$
is induced from $g$,
$C$ is a constant depending only on the dimension $N-1$.

\subsection{Local formulation of the scaled problem}\label{subsection2.2}
If we set $u(\bar{z})=\tilde{u}(\varepsilon \bar{z})$, then problem (\ref{originalproblemone})
is equivalent to the {\bf scaled problem}
\begin{align}\label{problemafterscaling0}
\Delta_{\bar{z}} u-u+u^{p}=0   \quad \mbox{in} \,\,\Omega_{\varepsilon},
\qquad
\frac{\partial u}{\partial \mathbf{n}_{\varepsilon}}=0  \quad   \mbox{on}\,\,\partial\Omega_{\varepsilon}.
\end{align}
After rescaling, we denote $\Omega_{\varepsilon}=\Omega/\ve$ and also $\Gamma_{\varepsilon}=\Gamma/\varepsilon$.

\medskip
To get local form of problem (\ref{problemafterscaling0})
and construct the approximation to a solution of (\ref{problemafterscaling0}),
which concentrates near $\Gamma_{\varepsilon}$,
after rescaling, we also introduce the scaled modified Fermi coordinates in the neighborhood
of $\Gamma_{\varepsilon}$ by
\begin{align}\label{coordinatessz}
(s,z)=(s,z_1,z_2)\in\Big(-\frac{r_0}{\ve},\,\frac{r_0}{\ve}\Big)\times\mathbb{D}_{\ve},
\end{align}
where $\mathbb{D}_{\var}:=\{z\in\mathbb{R}^2:|z|<1/\var\}$ for some fixed constant $\var>0$.
More precisely, we denote
\begin{align}
\gamma_{\var}(s,z)=\gamma(\ve s, \ve z)/\ve.
\end{align}
However, for the portion of $\partial\Omega_{\var}$ in ${\bar\Omega}_{\var}$,
we use the coordinate $(\theta,\eta)\in\partial\mathbb{D}_{\var}\times[0,\delta/\var)\subset\mathbb{D}_{\var}$
to parameterize a small region of $\partial\Gamma_{\var}$ in $\Gamma_{\var}$,
where the points $(\theta,0)$ is sent to the boundary of surface $\partial\Gamma_{\var}$ by $\gamma_0$ and there holds
\begin{align}
\Big<\frac{\partial\gamma_0}{\partial\theta},\frac{\partial\gamma_0}{\partial\eta}\Big>=0\ \mbox{ for }\eta=0.
\end{align}

As a direct consequence of Lemma  \ref{lemma2.1}, there holds
\begin{lem}\label{lemma2.2}
There exist a constant $r_0>0$, which depend only on $\Gamma$ and $\partial\Omega$,
such that
\begin{enumerate}
\item
$\gamma_{\var}(0,z)=\gamma_0(\var z)/\var$ for $z\in\mathbb{D}_{\var}$,
$\quad \gamma_{\var}(s,z)\in\partial\Omega_{\var}$ for $ z\in\partial\mathbb{D}_{\var},\,
|s|<\frac{r_0}{\var}$;

\item
 $\frac{\partial\gamma_{\var}}{\partial s}(0,z)=\nu(\var z)$ for $z\in\mathbb{D}_{\var}$;

\item
$\gamma_{\var}(s,z)$ has the following expansion, as $s\rightarrow 0$
\begin{align}\label{coordinateexpressionnew}
\gamma_{\var}(s,z)=\gamma_0(\var z)/\var+s\nu(\var z)+\frac{\var s^2}{2}q_{1}(\var z)
+\frac{\var^{2}s^3}{6}q_{2}(\var z)+O(\var^{3}s^4),\quad  z\in\mathbb{D}_{\var},
\end{align}
where $q_{1}(\var z)$ and $q_{2}(\var z)$ are vector functions defined in (\ref{vpq}),
which are orthogonal to $\nu(\var z)$.

\item
If we write $\gamma_{\var}$ as $\gamma_{\var}( s,\theta,\eta)$ in terms of the coordinate $(s,\theta,\eta)$,
then the derivative along the inward unit normal vector ${\bf n}_{\var}$ of $\partial\Omega_{\var}$ is expressed as
\begin{align}
\frac{\partial}{\partial {\bf n}_{\var}}=\frac{1}{\sqrt{g^{33}}}\Big(g^{13}\frac{\partial}{\partial s}
+g^{23}\frac{\partial}{\partial \theta}
+g^{33}\frac{\partial}{\partial\eta}\Big),
\end{align}
where at $(s,\theta,0)$, we have the expression
\begin{align}
\begin{aligned}
g^{13}(s,\theta)&=-\var s\Big|\frac{\partial\gamma_0}{\partial\rho}\Big|^{-2}\Big<q_{1},
\frac{\partial\gamma_0}{\partial\rho}\Big>+O(\var^{2}s^2),
\\
g^{23}(s,\theta)&=2\var s\Big|\frac{\partial\gamma_0}{\partial\vartheta}\Big|^{-2}
\Big|\frac{\partial\gamma_0}{\partial\rho}\Big|^{-2}\Big<\frac{\partial\gamma_0}{\partial\vartheta},
 \frac{\partial\nu}{\partial\rho}\Big>+O(\var^{2}s^2),
\\
g^{33}(s,\theta)&=\Big|\frac{\partial\gamma_0}{\partial\rho}\Big|^{-2}
+2\var s\Big|\frac{\partial\gamma_0}{\partial\rho}\Big|^{-4}\Big<\frac{\partial\gamma_0}{\partial\rho},
\frac{\partial\nu}{\partial\rho}\Big>+O(\var^{2}s^2).
\end{aligned}
\end{align}
\end{enumerate}
\end{lem}
\proof The reader can refer to Lemma \ref{lemma2.1} for the details of the proof.
\qed

Obviously, in the scaled modified Fermi coordinates, it is of importance to express Laplace-Beltrami operator
in (\ref{problemafterscaling0}) in local form. In fact, we obtain the form by careful calculations
\begin{align}\label{delta}
\Delta=\frac{\partial^{2}}{\partial s^{2}}+\Delta^{\Gamma_{\var}}+B_{0}(s,\var z)+B_{1}(s,\var z)+B_{2}(s,\var z),
\end{align}
where
\begin{align*}
B_{0}(s,\var z)&=\var k\frac{\partial}{\partial s}
-\var^{2}s\big(k^{2}_{1}(\var z)+k^{2}_{2}(\var z)\big)\frac{\partial}{\partial s}
+\var^{3}b(\var z)s^{2}\frac{\partial}{\partial s}\,,
\\
B_{1}(s,\var z)&=-\nabla^{\Gamma_{\var}}_{q_{1\var}}
-2 s\nabla^{\Gamma_{\var}}_{q_{1\var}}\Big(\frac{\partial}{\partial s}\Big)\,,
\\
B_{2}(s,\var z)&=\var^{2}\sum_{i=1}^{2}a_{i}(\var s,\var z)\frac{\partial}{\partial z_{i}}
+\var^{2}\sum_{i=1}^{2}b_{i}(\var s,\var z)\frac{\partial^{2}}{\partial z_{i}\partial s}+B_{3}(s,\var z)\,.
\end{align*}
In the above, we have denoted
\begin{align*}
\nabla^{\Gamma_{\var}}_{q_{1\var}}&=\var{\hslash}^{-2}(\var z)\sum^{2}_{j=1}
\Big\langle q_{1}(0,\var z),\,\frac{\partial\gamma_{0}}
{\partial y^{j}}\Big\rangle\frac{\partial}{\partial z^{j}}.
\end{align*}
The functions $b(\var z),\, a_{i}(\var s,\var z),b_{i}(\var s,\var z),i=1,2$ satisfy:
\begin{align*}
\Big|b(\var z),\, a_{i}(\var s,\var z),b_{i}(\var s,\var z)\Big|\leq C\big(1+|\var s|^{4}\big).
\end{align*}
Note that $O(\varepsilon)$-term  in $B_{0}(s,\var z)$ is actually absent because
$\Gamma$ is minimal (cf. ({\bf A1})).
Here we have denoted $\bigtriangleup^{\Gamma_{\var}}$ the Laplace-Beltrami operator on $\Gamma_{\var}$.
The vector $q_{1}$ is defined in (\ref{coordinateexpressionnew}) and the differential operator $B_{3}(s,\var z)$
is of size $O(\varepsilon^{4})$.
Whence, by using of Lemma \ref{lemma2.2},
we can get local form of problem (\ref{problemafterscaling0}).
This will be done in more details in Sections \ref{section3} and \ref{section4}.

\section{Inner Approximate Solutions}\label{section3}

In this section, we neglect the boundary condition in (\ref{problemafterscaling0})
and then find a local inner approximate solution to solve the first equation in (\ref{problemafterscaling0})
up to order of $O(\ve^4)$.
In fact, by the scaled modified Fermi coordinates in Lemma \ref{lemma2.2},
we write down the local form of the first equation in problem (\ref{problemafterscaling0})
and then extend it to a differential equation with unknown functions
defined on the infinite strip ${\mathfrak S}$ in $\R^3$ with notations
\begin{align}
\begin{aligned}\label{mathfraks0}
\mathfrak{S}\,=\, \{\ (x,z): x\in \mathbb{R},\, z\in\Gamma_{\var} \},
\qquad
\partial \mathfrak{S}=\{\ (x,z): x\in \mathbb{R},\, z\in\partial\Gamma_{\var} \}.
\end{aligned}
\end{align}

\subsection{Inner formulation of the scaled problem}
As we have mentioned in subsection \ref{subsection2.2},
in terms of the local coordinate system $(s,z)$(cf. (\ref{coordinatessz})),
in the neighborhood of $\Gamma_\ve$,
the differential equation in (\ref{problemafterscaling0}) is locally expressed as
\begin{align}\label{problemafterscaling1}
u_{ss}\,+\,\Delta^{\Gamma_{\var}}u\,-\, u\,+\,u^p+\, B(u)=0,
\end{align}
where the linear differential operator $B(u)$ is defined by
\begin{align*}
B(u)&=B_{0}(u)+B_{1}(u)+B_{2}(u),
\end{align*}
by the notations
\begin{align*}
B_{0}(u)&=-\ve^2\, s(k^{2}_{1}(\var z)+k^{2}_{2}(\var z))\frac{\partial u}{\partial s}
+\varepsilon^{3}b(\var z)s^{2}\frac{\partial u}{\partial s},
\\
B_{1}(u)&=-\bigtriangledown^{\Gamma_{\var}}_{q_{1\var}}(u)-2\,s\bigtriangledown^{\Gamma_{\var}}_{q_{1\var}}(\frac{\partial u}{\partial s}),\\
B_{2}(u)&=\var^{2}\sum_{i=1}^{2}a_{i}(\var s,\var z)\frac{\partial u}{\partial z_{i}}
+\var^{2}\sum_{i=1}^{2}b_{i}(\var s,\var z)\frac{\partial^{2} u}{\partial z_{i}\partial s}+B_{3}(u).
\end{align*}

We assume that, in the $(s,z)$ coordinates, the location of concentration of the solution
is characterized by the surface
\begin{align}\label{layerf}
\widetilde{\Gamma}_{\varepsilon}:\,s=\sum^{2}_{i=0}\varepsilon^{i}f_{i}(\var z).
\end{align}

\begin{Remark}
The smooth functions $f_{0},\,f_{1}$ are to be determined in Sections \ref{section3} and \ref{section4}.
More precisely, with the help of the validity of the nondegeneracy of $\Gamma$
in subsection \ref{subsection1.1}, we will choose $f_0$ by solving the equation (\ref{conditionf0}) with boundary condition (\ref{boundaryf0}), and also the equation (\ref{conditionf1}) with boundary condition (\ref{boundaryf1}) for $f_1$.
In fact, the nondegeneracy of $\Gamma$ implies that $f_0$ is identically zero.
While the unknown parameter $f_{2}$ is to be chosen by a type of reduction procedure,
which is equivalent to solving a system of differential equations in Section \ref{section9}
(cf.  (\ref{nonlocalsystemoffandeone})-(\ref{nonlocalsystemoffandefour})).
\qed
\end{Remark}

In the sequel, we always assume that $f_{2}$ satisfies the uniform constraint
\begin{align}\label{conditionforf}
\|f_{2}\|_{a}=\|f_{2}\|_{\mathrm{L}^{\infty}(\Gamma)}+\|\nabla^{\Gamma}f_{2}\|_{\mathrm{L}^{\infty}(\Gamma)}
+\|\Delta^{\Gamma}f_{2}\|_{\mathrm{L}^{q}(\Gamma)}\leq\varepsilon^{\frac{1}{2}}.
\end{align}

We consider a further changing of variables and define a new function $v(x,z)$ as follows
\begin{align}\label{coordinatesxz}
u(s,z)=v(x,z)\quad\mbox{with}\quad x=s-\sum^{2}_{i=0}\varepsilon^{i}f_{i}(\var z),\quad z=z.
\end{align}
We now want to express the problem in the new coordinates. Whence we need the following formulas
\begin{align*}
u_{s}&=v_{x},
\qquad
u_{ss}=v_{xx},
\\
\triangle^{\Gamma_{\var}}(u)&=\var^{2}\Big(\sum\limits_{i=0}^{2}\varepsilon^{i}\nabla^{\Gamma}f_{i}\Big)^{2}v_{xx}-
\var^{2}\Big(\sum\limits_{i=0}^{2}\varepsilon^{i}\triangle^{\Gamma}f_{i}\Big)v_{x}-
\var\Big(\sum\limits_{i=0}^{2}\varepsilon^{i}\nabla^{\Gamma}f_{i}\Big)\nabla^{\Gamma_{\var}}v_{x}+\triangle^{\Gamma_{\var}}v,
\\
\nabla^{\Gamma_{\var}}_{q_{1\var}}(u)&=-\var\Big(\sum\limits_{i=0}^{2}\varepsilon^{i+1}\nabla^{\Gamma}_{q_{1}}f_{i}\Big)v_{x}
+\var\nabla^{\Gamma_{\var}}_{q_{1}}v,\,\,\,\,\,\,
\qquad
\nabla^{\Gamma_{\var}}_{q_{1\var}}(u_{s})=-\var\Big(\sum\limits_{i=0}^{2}\varepsilon^{i+1}\nabla^{\Gamma}_{q_{1}}f_{i}\Big)v_{xx}
+\var\nabla^{\Gamma_{\var}}_{q_{1}}(v_{x}),
\\
\frac{\partial u}{\partial z_{i}}&=-\Big(\sum\limits_{j=0}^{2}\varepsilon^{j+1}\frac{\partial f_{j}}
{\partial y_{i}}\Big)v_{x}+\frac{\partial v}{\partial z_{i}},\,\,\,\,\,\,
\qquad
\frac{\partial^{2} u}{\partial z_{i}\partial s}=-\Big(\sum\limits_{j=0}^{2}\varepsilon^{j+1}\frac{\partial f_{j}}{\partial y_{i}}\Big)v_{xx}+\frac{\partial^{2} v}{\partial z_{i}\partial x}.
\end{align*}

Locally, this gives that $u$ solves (\ref{problemafterscaling1}) if and only if the function $v$ defined
in (\ref{coordinatesxz}) solves the following problem
\begin{equation}\label{problemafterscaling2}
S(v)\equiv v_{xx}+\Delta^{\Gamma_{\var}}v-v+v^{p}+B(v)=0.
\end{equation}
Note that we will consider problem (\ref{problemafterscaling2})
on the whole ${\mathfrak S}$(cf. (\ref{mathfraks0})).
In the above we have denoted the linear operator
\begin{equation}
B(v)=B_{4}(v)+B_{5}(v)+B_{6}(v).
\end{equation}
The linear operator $B_{4}(v),B_{5}(v)$\,and\,$ B_{6}(v)$ can be expressed explicitly by
\begin{align*}
B_{4}(v)=&\var^{2}\Big(\sum\limits_{i=0}^{2}\varepsilon^{i}\nabla^{\Gamma}f_{i}\Big)^{2}v_{xx}-
\var^{2}\Big(\sum\limits_{i=0}^{2}\varepsilon^{i}\triangle^{\Gamma}f_{i}\Big)v_{x}-
\var\Big(\sum\limits_{i=0}^{2}\varepsilon^{i}\nabla^{\Gamma}f_{i}\Big)\nabla^{\Gamma_{\var}}v_{x}\\
& -\varepsilon^{2}\Big(x+\sum\limits_{i=0}^{2}\varepsilon^{i}f_{i}\Big)(k^{2}_{1}+k_{2}^{2})v_{x}
+\varepsilon^{3}b\Big(x+\sum^{2}_{i=0}\varepsilon^{i}f_{i}\Big)^{2}v_{x},
\\
B_{5}(v)=&2\varepsilon\Big(x+\sum\limits_{i=0}^{2}\varepsilon^{i+1}f_{i}\Big)
\Big[\Big(\sum\limits_{i=0}^{2}\varepsilon^{i}\nabla^{\Gamma}_{q_{1}}f_{i}\Big)v_{xx}-\nabla^{\Gamma_{\var}}_{q_{1}}(v_{x})\Big]
+\var\Big[\Big(\sum\limits_{i=0}^{2}\varepsilon^{i+1}\nabla^{\Gamma}_{q_{1}}f_{i}\Big)v_{x}-\nabla^{\Gamma_{\var}}_{q_{1}}v\Big],
\\
B_{6}(v)=&\var^{2}\sum\limits_{i=1}^{2}a_{i}\Big[-\Big(\sum\limits_{j=0}^{2}\varepsilon^{j+1}\frac{\partial f_{j}}{\partial y_{i}}\Big)v_{x}
+\frac{\partial v}{\partial z_{i}}\Big]+\var^{2}\sum\limits_{i=1}^{2}b_{i}
\Big[-\Big(\sum\limits_{j=0}^{2}\varepsilon^{j+1}\frac{\partial f_{j}}{\partial y_{i}}\Big)v_{xx}
+\frac{\partial^{2} v}{\partial z_{i}\partial x}\Big]+B_{3}(v).
\end{align*}

\subsection{Inner approximate solutions}
In the subsection, we want to use coordinates $(x,z)$ defined in (\ref{coordinatesxz})
 to construct a suitable approximation to a solution expressed in the form,
\begin{align}\label{basicapproximation}
\mathcal{V}=w(x)+\varepsilon e(\var z)Z(x)+\sum\limits_{i=1}^{3}\varepsilon^{i}\varphi_{i}(x,\var z),
\end{align}
where $w$ and $Z$ are two functions given by (\ref{definitionofw}) and (\ref{definitionofZ}).
In the above expression, we have denoted $\varphi_{i},\,i=1,2,3$,\,
smooth bounded functions to be determined in the sequel. As we have mentioned, the unknown parameters
 $f_{2}$(cf. (\ref{layerf})) and $e$ will be chosen in the last section
 by solving a system of differential
 equations(cf. (\ref{nonlocalsystemoffandeone})-(\ref{nonlocalsystemoffandefour})).
 In all what follows, we shall assume the validity of the following uniform constraints on the parameter $e$
\begin{align}\label{conditionfore}
\|e\|_{b}=\|e\|_{\mathrm{L}^{\infty}(\Gamma)}+\varepsilon\|\nabla^{\Gamma}e\|_{\mathrm{L}^{q}(\Gamma)}
+\varepsilon^{2}\|\Delta^{\Gamma}e\|_{\mathrm{L}^{q}(\Gamma)}\leq\varepsilon^{\frac{1}{2}}.
\end{align}
For simplicity of notations, define
\begin{align}\label{regionoffande}
\digamma=\Bigl\{\ (f_{2},e)|
       \mbox{ the functions } f_{2} \mbox{ and } e \mbox{ satisfy }
       (\ref{conditionforf})\mbox{ and } (\ref{conditionfore})\ \mbox{respectively} \ \Bigr\}.
\end{align}
Now the key point is to choose suitable correction
 terms $\varphi_{1},\,\varphi_2,\,\varphi_{3}$, and then prove
 that the approximate solution $\mathcal{V}$ solve problem (\ref{problemafterscaling2})
 up to order $O(\varepsilon^{4})$.

Formally, we have
\begin{align}\label{expansion}
(w+\Theta)^{p}=w^{p}\Bigl[1+p\frac{\Theta}{w}+\frac{p(p-1)}{2}\Big(\frac{\Theta}{w}\Big)^{2}+...
+C_{k,p}\Big(\frac{\Theta}{w}\Big)^{k}+O\Big(\Big|\frac{\Theta}{w}\Big|^{k+1}\Big)\Bigr].
\end{align}
Setting $\Theta=\varepsilon(eZ+\varphi_{1})+\sum\limits_{i=2}^{3}\varepsilon^{i}\varphi_{i}$ and separating the powers of $\varepsilon$, we get
\begin{align*}
\Big(w+\varepsilon eZ+\sum\limits_{i=1}^{3}\varepsilon^{i}\varphi_{i}\Big)^{p}=
w^{p}\sum\limits_{l=0}^{4}\varepsilon^{l}\sum\limits_{j_{1},...,j_{3},\sum ij_{i}={l}}C_{l,j_{1},...j_{3}}
\frac{(eZ+\varphi_{1})^{j_{1}}\varphi_{2}^{j_{2}}\varphi_{3}^{j_{3}}}{w^{j_{1}+...+j_{3}}}
+w^{p}O\Big(\Big|\frac{\Theta}{w}\Big|^{5}\Big).
\end{align*}
Whence, using elementary calculation, we collect the powers of $\varepsilon$ up to order 4 in the last formula, and then get the estimate
\begin{align*}
\mathcal{H}_{0}&\,=\,\Biggl|\,\Big(w+\varepsilon eZ+\sum\limits_{i=1}^{3}\varepsilon^{i}\varphi_{i}\Big)^{p}
\,-\,w^{p}\sum\limits_{l=0}^{4}\varepsilon^{l}\sum\limits_{j_{1},...j_{3},\sum ij_{i}={l}}C_{l,j_{1},...j_{3}}
\frac{(eZ+\varphi_{1})^{j_{1}}\varphi_{2}^{j_{2}}\varphi_{3}^{j_{3}}}{w^{j_{1}+...+j_{3}}}\,\Biggr|
\\
&\leq C_{4,p}w^{p}\Bigl[\varepsilon^{5}\Big(1+\frac{|\Theta|}{w}\Big)^{4}+\Big|\frac{\Theta}{w}\Big|^{5}\Bigr].
\end{align*}
More precisely, using (\ref{expansion}), we make a decomposition
\begin{align}\label{1.32}
\begin{aligned}
&w^{p}\sum\limits_{l=0}^{4}\varepsilon^{l}\sum\limits_{j_{1},...j_{3},\sum ij_{i}={l}}C_{l,j_{1},...j_{3}}
\frac{(eZ+\varphi_{1})^{j_{1}}\varphi_{2}^{j_{2}}\varphi_{3}^{j_{3}}}{w^{j_{1}+...+j_{3}}}
\\
&\qquad=w^{p}+pw^{p-1}\varepsilon eZ+\sum\limits_{i=1}^{3}\varepsilon^{i}pw^{p-1}\varphi_{i}
+\frac{1}{2}\varepsilon^{2}p(p-1)w^{p-2}(eZ+\varphi_{1})^{2}
\\
&\qquad\quad+p(p-1)w^{p-2}(eZ+\varphi_{1})\sum\limits_{i=3}^{4}\varepsilon^{i}\varphi_{i-1}+\mathcal{H}_{3}
\\
&\qquad\equiv\mathcal{H}_{1}+\mathcal{H}_{2}+\mathcal{H}_{3},
\end{aligned}
\end{align}
where we have denoted
\begin{align}
\mathcal{H}_{1}&=w^{p}+pw^{p-1}\varepsilon eZ+\sum\limits_{i=1}^{3}\varepsilon^{i}pw^{p-1}\varphi_{i}\,,
\end{align}
\begin{align}\label{expasion2}
\begin{aligned}
\mathcal{H}_{2}&=\frac{1}{2}\varepsilon^{2}p(p-1)w^{p-2}(eZ+\varphi_{1})^{2}
+p(p-1)w^{p-2}(eZ+\varphi_{1})\sum\limits_{i=3}^{4}\varepsilon^{i}\varphi_{i-1}.
\end{aligned}
\end{align}
In the above, we have also denoted that
\begin{align}\label{expasion3}
\mathcal{H}_{3}=\sum\limits_{l=3}^{4}\varepsilon^{l}\mathfrak{D}_{l},
\end{align}
where for every $l=3,4,$ the component $\mathfrak{D}_{l}$ is independent of the terms $\varphi_{l-1},...,\varphi_{3}$.

Putting $\mathcal{V}$ into (\ref{problemafterscaling2}) and expanding formally, we derive that
\begin{align}\label{problemafterscaling3}
\begin{aligned}
B(w)&+\var\big(\var^{2}\Delta^{\Gamma}e+\lambda_{0}e\big)Z+\sum\limits_{i=1}^{3}\varepsilon^{i}\big[\varphi_{i,xx}-\varphi_{i}+pw^{p-1}\varphi_{i}\big]
\\
&+B(\varepsilon eZ)+\sum\limits_{i=1}^{3}\var^{i}\Delta^{\Gamma_{\var}}\varphi_{i}+\sum\limits_{i=1}^{3}\varepsilon^{i}B(\varphi_{i})+\mathcal{H}_{2}+\mathcal{H}_{3}=0.
\end{aligned}
\end{align}

\subsection{Inner errors}\label{subsection3.2}
In this subsection, we accept that $\varphi_i$'s as known functions for the moment.
Then we compute all error terms in (\ref{problemafterscaling3})
and make decomposition of all components into suitable forms according to the order of $\ve$.

First of all, we calculate the term
\begin{align}
B(w)=B_{4}(w)+B_{5}(w)+B_{6}(w).
\end{align}
Direct calculation gives that
\begin{align*}
B_{4}(w)+B_{5}(w)=\sum\limits_{i=1}^{4}\varepsilon^{i}\check{S}_{i}
+\sum\limits_{i=1}^{4}\varepsilon^{i}\hat{S}_{i}+\sum\limits_{i=5}^{6}\varepsilon^{i}S_{i},
\end{align*}
with expressions defined by
\begin{align*}
\check{S}_{1}&=\hat{S}_{1}=0,\\
\check{S}_{2}&=-\Delta^{\Gamma}f_{0}w_{x}-(k_{1}^{2}+k_{2}^{2})f_{0}w_{x}, \qquad
\\
\hat{S}_{2}&=(\nabla^{\Gamma}f_{0})^{2}w_{xx}-(k_{1}^{2}+k_{2}^{2})xw_{x}+2f_{0}\nabla^{\Gamma}_{q_{1}}f_{0}w_{xx},
\\
\check{S}_{3}&=-\Delta^{\Gamma}f_{1}w_{x}-(k_{1}^{2}+k_{2}^{2})f_{1}w_{x}+b(f^{2}_{0}+x^{2})w_{x}, \qquad
\\
\hat{S}_{3}&=2f_{0}\nabla^{\Gamma}_{q_{1}}f_{1}w_{xx}+2f_{1}\nabla^{\Gamma}_{q_{1}}f_{0}w_{xx}+2bf_{0}xw_{x},
\\
\check{S}_{4}&=-\Delta^{\Gamma}f_{2}w_{x}-(k_{1}^{2}+k_{2}^{2})f_{2}w_{x}, \qquad
\\
\hat{S}_{4}&=(\nabla^{\Gamma}f_{1})^{2}w_{xx}+2f_{0}\nabla^{\Gamma}_{q_{1}}f_{2}w_{xx},
+2f_{1}\nabla^{\Gamma}_{q_{1}}f_{1}w_{xx}+2f_{2}\nabla^{\Gamma}_{q_{1}}f_{0}w_{xx}+2bf_{1}xw_{x}.
\end{align*}
Note that for any $i=1,\cdots,4$, $\check{S}_{i}$ is an odd function in the variable $x$, while $\hat{S}_{i}$ is an even function in the variable $x$. On the other hand, for any $i=5,6$, the high order term $S_{i}$ is combination of powers of the parameters of $f_{0},f_{1},f_{2}$ and their derivatives with smooth bounded coefficients.

At the meantime, the linear operator $B_{6}(w)$ comes from $B_{6}(v)$ and can be expressed explicitly by
\begin{align*}
B_6(w)&=\sum_{i=3}^4\ve^i{\mathbf{\check B}}_i(f_0,\cdots,f_{i-3})
\,+\,\sum_{i=3}^4\ve^i{\mathbf{\hat B}}_i(f_0,\cdots,f_{i-3})
\,+\,\ve^{5}{\mathbf B}_{5}(f_0,f_1,f_{2}),
\end{align*}
where the terms ${\mathbf{ B}}_{5}(f_0,f_1,f_{2})$ and
${\mathbf{{\check B}}}_i(f_0,\cdots,f_{l-3}),\, {\mathbf{{\hat B}}}_i(f_0,\cdots,f_{l-3}),i=3,4$
are combination of powers of the parameters $f_0,f_1,f_{2}$ and their derivatives with smooth bounded coefficients. Moreover, ${\mathbf{\check B}}_i(f_0,\cdots,f_{i-3}), i=3,4$ are odd functions in the
variable $x$, while ${\mathbf{\hat B}}_i(f_0,\cdots$, $f_{i-3}),\,\,i=3,4$ are even functions in the
variable $x$.
As a conclusion, we get
\begin{align}
\begin{aligned}
B(w)&=\sum\limits_{i=1}^{4}\varepsilon^{i}\check{S}_{i}
+\sum\limits_{i=1}^{4}\varepsilon^{i}\hat{S}_{i}+\sum\limits_{i=5}^{6}\varepsilon^{i}S_{i}
+\sum_{i=3}^4\ve^i{\mathbf{\check B}}_i
+\sum_{i=3}^4\ve^i{\mathbf{\hat B}}_i
+\ve^{5}{\mathbf{ B}}_{5}.\label{definitionofS(w)}
\end{aligned}
\end{align}

Second, we compute the error
\begin{align}
\var(\var^{2}\Delta^{\Gamma} eZ+\lambda_{0}eZ)+B(\varepsilon eZ),
\end{align}
with $B(\varepsilon eZ)=B_{4}(\varepsilon eZ)+B_{5}(\varepsilon eZ)+B_{6}(\varepsilon eZ)$. There also holds
\begin{align*}
B_{4}(\varepsilon eZ)+B_{5}(\varepsilon eZ)+B_{6}(\varepsilon eZ)=\sum\limits_{i=1}^{3}\varepsilon^{i}\check{T}_{i}
+\sum\limits_{i=1}^{3}\varepsilon^{i}\hat{T}_{i}+\sum\limits_{i=4}^{6}\varepsilon^{i}T_{i}.
\end{align*}
In the above, we have denoted the following forms
\begin{align*}
\check{T}_{1}&=\hat{T}_{1}=0,\qquad \check{T}_{2}=\hat{T}_{2}=0,
\\
\check{T}_{3}&=-\Delta^{\Gamma}f_{0}eZ_{x}-f_{0}(k_{1}^{2}+k_{2}^{2})eZ_{x}
-\nabla^{\Gamma}f_{0}\nabla^{\Gamma}eZ_{x}-2f_{0}\nabla^{\Gamma}_{q_{1}}eZ_{x}, \qquad
\\
\hat{T}_{3}&=(\nabla^{\Gamma}f_{0})^{2}eZ_{xx}
-(k_{1}^{2}+k_{2}^{2})exZ_{x}-\nabla^{\Gamma}_{q_{1}}eZ+2f_{0}\nabla^{\Gamma}_{q_{1}}f_{0}eZ_{xx}
-2\nabla^{\Gamma}_{q_{1}}exZ_{x},
\\
T_{i}&=\breve{b}_{i}(f_{0},f_{1},f_{2},e),\quad i=4,5,6.
\end{align*}
In the above, for any $i=4,5,6$, the term $\breve{b}_{i}(f_{0},f_{1},f_{2},e)$ is combination of powers of the parameter
$f_{0},f_{1},f_{2},e$ and their derivatives with smooth bounded coefficients. Moreover, for any $i=1,2,3$,
$\check{T}_{i}$ is an odd function in the variable $x$, while $\hat{T}_{i}$ is  an even function in the variable $x$.

In summary, we have that
\begin{align}\label{definitionofs(v)}
\begin{aligned}
S(\mathcal{V})=&\sum\limits_{i=1}^{4}\varepsilon^{i}\check{S}_{i}+\sum\limits_{i=1}^{4}\varepsilon^{i}\hat{S}_{i}
+\sum\limits_{i=5}^{6}\varepsilon^{i}S_{i}
+\sum_{i=3}^4\ve^i{\mathbf{\check B}}_i+\sum_{i=3}^4\ve^i{\mathbf{\hat B}}_i
\\
&+\ve^{5}{\mathbf{ B}}_{5}
+\sum\limits_{i=1}^{3}\varepsilon^{i}\check{T}_{i}
+\sum\limits_{i=1}^{3}\varepsilon^{i}\hat{T}_{i}
+\sum\limits_{i=4}^{6}\varepsilon^{i}T_{i}
\\
&+\varepsilon^{3}\Delta^{\Gamma}eZ+\varepsilon \lambda_{0}eZ+\sum\limits_{i=1}^{3}\varepsilon^{i}[\varphi_{i,xx}-\varphi_{i}+pw^{p-1}\varphi_{i}]
\\
&+\sum\limits_{i=1}^{3}\var^{i}\Delta^{\Gamma_{\var}}\varphi_{i}+\sum\limits_{i=1}^{3}\varepsilon^{i}B(\varphi_{i})+\mathcal{H}_{2}+\mathcal{H}_{3}.
\end{aligned}
\end{align}

Now, we shall write the error terms involving correction terms $\varphi_{1},...,\varphi_{3}$ in a suitable form. In the next subsection, for any given $i=1,...,3$, we will choose $\varphi_{i}$ as the form
\begin{align*}
a_{i1}(\var z)b_{i1}(x)+a_{i2}(\var z)b_{i2}(x),
\end{align*}
for some generic smooth functions $a_{i1},a_{i2},b_{i1}\,$(odd)\,and\,$b_{i2}\,$(even). Moreover, the terms $a_{i1}$
and $a_{i2}$ do not depend on the unknown parameters $f_{0},f_{1}$\,and\,$f_{2}$. The reader can refer to (\ref{solution2}) and (\ref{solution3}).

Whence,  we make a decomposition as
\begin{align}
\sum_{i=1}^{3}\var^i\Delta^{\Gamma_{\var}}\varphi_{i}
\,=\,\sum_{j=3}^4\var^j{\mathcal{\check N}}_j
\,+\,\sum_{j=3}^4\var^j{\mathcal{\hat N}}_j
\,+\,\var^{5}{\mathcal N}_{5}.
\end{align}
where
\begin{align*}
{\mathcal{\check N}}_j\,=\,&{\mathcal{\check N}}_j(\varphi_1,\cdots,\varphi_{j-2},f_0,\cdots,f_{j-3}),
\\
{\mathcal{\hat N}}_j\,=\,&{\mathcal{\hat N}}_j(\varphi_1,\cdots,\varphi_{j-2},f_0,\cdots,f_{j-3}),
\\
{\mathcal{N}}_{5}\,=\,&{\mathcal N}_{5}(\varphi_1,\cdots,\varphi_{3},f_0,\cdots,f_{k-2}).
\end{align*}
Moreover, ${\mathcal {\check N}}_j$ is an odd function in the variable $x$, while ${\mathcal {\hat N}}_j$ is an
even functions in the variable $x$.

From the definition of the operator $B$ in (\ref{problemafterscaling2}), we also write
\begin{align*}
\sum\limits_{i=1}^{3}\varepsilon^{i}B(\varphi_{i})=\sum\limits_{j=3}^{4}\varepsilon^{j}\check{G}_{j}
+\sum\limits_{j=3}^{4}\varepsilon^{j}\hat{G}_{j}+\sum\limits_{j=5}^{10}\varepsilon^{j}G_{j},
\end{align*}
where, for any $j=3,4,$ the components $\check{G}_{j}$ and $\hat{G}_{j}$ do not depend on the correction
terms $\varphi_{j-1},...,\varphi_{3}$ and the unknown parameters $f_{j-2},...,f_{2}$. In other words, we have
\begin{align*}
\check{G}_{j}&=\check{G}_{j}(\varphi_{1},...,\varphi_{j-2},f_{0},...,f_{j-3}),
\\
\hat{G}_{j}&=\hat{G}_{j}(\varphi_{1},...,\varphi_{j-2},f_{0},...,f_{j-3}),
\\
G_{j}&=G_{j}(\varphi_{1},...,\varphi_{3},f_{0},...,f_{2}).
\end{align*}
Moreover,\,$\check{G}_{j}$ is an odd function in the variable $x$, while $\hat{G}_{j}$ is  an even function in the variable $x$.

For later use, using (\ref{expasion2}) and (\ref{expasion3}), we decompose $\mathcal{H}_{3}$ into even parts and odd parts, and then write
$\mathcal{H}_{2}+\mathcal{H}_{3}$ as
\begin{align}
\begin{aligned}
\mathcal{H}_{2}+\mathcal{H}_{3}&=\frac{1}{2}\varepsilon^{2}p(p-1)w^{p-2}(eZ+\varphi_{1})^{2}
\\
&\quad+p(p-1)w^{p-2}(eZ+\varphi_{1})\sum\limits_{i=3}^{4}\varepsilon^{i}\varphi_{i-1}
+\sum\limits_{j=3}^{4}\varepsilon^{j}(\check{\mathfrak{D}}_{j}+\hat{\mathfrak{D}}_{j}),
\end{aligned}
\end{align}
For $j=3,4,$ the components $\check{\mathfrak{D}}_{j}$ and $\hat{\mathfrak{D}}_{j}$ are independent of the terms $\varphi_{j-1},...,\varphi_{3}$,\\
i.e.
\begin{align}
\check{\mathfrak{D}}_{j}=\check{\mathfrak{D}}_{j}(\varphi_{1},...,\varphi_{j-2}),
\qquad
\hat{\mathfrak{D}}_{j}=\hat{\mathfrak{D}}_{j}(\varphi_{1},...,\varphi_{j-2}),
\end{align}
Moreover,\,\,$\check{\mathfrak{D}}_{j}$ is an odd function in the variable $x$, while $\hat{\mathfrak{D}}_{j}$ is
an even function in the variable $x$.

\subsection{Determinations of the inner correction terms $\varphi_i$'s}\label{subsection3.3}
In this subsection,  by a recurrence procedure, we will choose suitable parameters $f_{0},\,f_{1}$
so that we can really find the correction terms $\varphi_{2},\, \varphi_{3}$
and then improve the approximation.
In fact, it will be shown that $\varphi_1$ is identically zero.

This can be done in the following way.
It is worth mentioning that the term
$$
\varepsilon^{3}\Delta^{\Gamma}eZ+\varepsilon\lambda_{0}eZ
$$
lies in the approximate kernel of  the linearized problem of problem (\ref{problemafterscaling2}) at ${\mathcal V}$.
We ignore this term for the moment
and then cancel other components of the error in (\ref{definitionofs(v)})
with order of $\ve$ lower than $4$ by choosing suitable
correction terms $\varphi_{1},...,\varphi_{3}$.
Whence, for given $y\in\Gamma$,
we then consider the problems
\begin{align}\label{defitionsolution1}
\begin{aligned}
\varphi_{1,xx}-\varphi_{1}+pw^{p-1}\varphi_{1}=0  \quad \mbox{in} \,\,\mathbb{R},
\\
\varphi_{1}(\pm\infty)=0,\quad \int_{\mathbb{R}}\varphi_{1}w_{x}\,\mathrm{d}x=0;
\end{aligned}
\end{align}
\begin{align}\label{defitionsolution2}
\begin{aligned}
\varphi_{2,xx}-\varphi_{2}+pw^{p-1}\varphi_{2}&=-\check{S}_{2}-\hat{S}_{2}-\frac{1}{2}p(p-1)w^{p-2}(eZ+\varphi_{1})^{2}
 \quad \mbox{in} \,\,\mathbb{R},
\\
\varphi_{2}(\pm\infty)&=0,\quad \int_{\mathbb{R}}\varphi_{2}w_{x}\,\mathrm{d}x=0;
\end{aligned}
\end{align}
and
\begin{align}\label{defitionsolution34}
\begin{aligned}
\varphi_{3,xx}-\varphi_{3}+pw^{p-1}\varphi_{3}&=-\check{S}_{3}-\hat{S}_{3}-p(p-1)w^{p-2}
(eZ+\varphi_{1})\varphi_{2}-\mathcal{C}_{3}   \quad \mbox{in} \,\,\mathbb{R},
\\
\varphi_{3}(\pm\infty)&=0,\quad \int_{\mathbb{R}}\varphi_{3}w_{x}\,\mathrm{d}x=0.
\end{aligned}
\end{align}
In the above, we have denoted $\mathcal{C}_{3}=\check{\mathcal{C}}_{3}+\hat{\mathcal{C}}_{3}$,
with the odd part and even part given by
\begin{align}
\check{\mathcal{C}}_{3}=\mathbf{\check{B}}_{3}+\check{T}_{3}+\mathcal{\check{N}}_{3}+\check{G}_{3}+\check{\mathfrak{D}}_{3},\,\,\, \,\,\,
\hat{\mathcal{C}}_{3}=\mathbf{\hat{B}}_{3}+\hat{T}_{3}+\mathcal{\hat{N}}_{3}+\hat{G}_{3}+\hat{\mathfrak{D}}_{3}.
\end{align}

Using (\ref{defitionsolution1}),  it is easy to show $\varphi_{1}=0$ and we finish the first step.
To proceed the second step and cancel the error terms of order $O(\varepsilon^{2})$ for the
improvement of the approximation,
we should choose the correction term $\varphi_{2}$ by solving problem (\ref{defitionsolution2}).
For this purpose, first, we collect all terms of order $O(\varepsilon^{2})$ in $S(\mathcal{V})$, which has the form $\varepsilon^{2}\mathcal{A}_{2}$ with
\begin{align*}
\mathcal{A}_{2}=\check{S}_{2}+\hat{S}_{2}+\frac{1}{2}p(p-1)w^{p-2}(eZ)^{2}.
\end{align*}
We denote the odd part and even part respectively by $\varepsilon^{2}\check{\mathcal{A}}_{2}$
and $\varepsilon^{2}\hat{\mathcal{A}}_{2}$ with
\begin{align*}
\check{\mathcal{A}}_{2}&=-\Delta^{\Gamma}f_{0}w_{x}-(k_{1}^{2}+k_{2}^{2})f_{0}w_{x},
\\
\hat{\mathcal{A}}_{2}&=(\nabla^{\Gamma}f_{0})^{2}w_{xx}-(k_{1}^{2}+k_{2}^{2})xw_{x}+2f_{0}\nabla^{\Gamma}_{p}f_{0}w_{xx}
+\frac{1}{2}p(p-1)w^{p-2}(eZ)^{2}.
\end{align*}
Then, we consider the problem
\begin{align}\label{solutionvarphi2}
\begin{aligned}
-\varphi_{2,xx}+\varphi_{2}-pw^{p-1}\varphi_{2}&=\check{\mathcal{A}}_{2}+\hat{\mathcal{A}}_{2}   \quad \mbox{in} \,\,\mathbb{R},
\\
\varphi_{2}(\pm\infty)=0,\quad &\int_{\mathbb{R}}\varphi_{2}w_{x}\,\mathrm{d}x=0,
\end{aligned}
\end{align}
as it is well known, which is uniquely solvable provided that
\begin{align}\label{orthogonalitycondition}
\int_{\mathbb{R}}(\check{\mathcal{A}}_{2}+\hat{\mathcal{A}}_{2})w_{x}\,\mathrm{d}x=0.
\end{align}
In fact, using the fact that $w$ is an even function in the variable $x$, we have
\begin{align*}
\int_{\mathbb{R}}\hat{\mathcal{A}}_{2}w_{x}\,\mathrm{d}x=0.
\end{align*}
On the other hand, we get
\begin{align}\label{orthogonalityw2}
\begin{aligned}
\int_{\mathbb{R}}\check{\mathcal{A}}_{2}w_{x}\,\mathrm{d}x&=-\Delta^{\Gamma}f_{0}\int_{\mathbb{R}}w_{x}^{2}\,\mathrm{d}x
-(k_{1}^{2}+k_{2}^{2})f_{0}\int_{\mathbb{R}}w_{x}^{2}\,\mathrm{d}x
\\
&=-\big[\Delta^{\Gamma}f_{0}+(k_{1}^{2}+k_{2}^{2})f_{0}\big]\int_{\mathbb{R}}w_{x}^{2}\,\mathrm{d}x.
\end{aligned}
\end{align}
While (\ref{orthogonalityw2}) implies that
\begin{align*}
\int_{\mathbb{R}}\check{\mathcal{A}}_{2}w_{x}\,\mathrm{d}x=0,
\end{align*}
is equivalent to the following differential equation
\begin{align}\label{conditionf0}
\Delta^{\Gamma}f_{0}+(k_{1}^{2}+k_{2}^{2})f_{0}=0\quad\mbox{on }\Gamma.
\end{align}
In fact, we find that $f_0$ is identically zero by combining the equation (\ref{conditionf0})
with the boundary condition (\ref{boundaryf0}) with the help of the non-degeneracy
of $\Gamma$ in (\ref{Geometriceigenvalueproblem}).
Whence, the solution to (\ref{defitionsolution2}) can be expressed as
\begin{align}\label{solution2}
\varphi_{2}(x,z)=\psi_{21}(x,\var z)+\psi_{22}(x,\var z),
\end{align}
where $\psi_{21}(x,\var z)$ is in fact identically zero and $\psi_{22}(x,\var z)$ is an even function in the variable $x$.
The components in $\psi_{21}(x,\var z)$ and $\psi_{22}(x,\var z)$
are independent of the parameters $f_{1}$\,and\,$f_{2}$.

In the same way, in order to cancel the error terms of order $O(\varepsilon^{3})$ and improve the approximation by solving problem (\ref{defitionsolution34}), we collect all terms of order $O(\varepsilon^{3})$ in $S(\mathcal{V})$, which has the form $\varepsilon^{3}\mathcal{A}_{3}$ with
\begin{align*}
\mathcal{A}_{3}=\check{S}_{3}+\hat{S}_{3}+\frac{1}{2}p(p-1)w^{p-2}(eZ)^{2}
+p(p-1)w^{p-2}(eZ+\varphi_{1})\varphi_{2}+\mathcal{C}_{3}.
\end{align*}
We denote the odd part and even part respectively by $\varepsilon^{3}\check{\mathcal{A}}_{3}$
and $\varepsilon^{3}\hat{\mathcal{A}}_{3}$.
As the arguments in solving (\ref{solutionvarphi2}),
we need an orthogonality condition like (\ref{orthogonalitycondition}).
Hence, we compute the projection of $\check{\mathcal{A}}_{3}$ and $\hat{\mathcal{A}}_{3}$
onto the kernel of the operator
$\partial^2/{\partial x^2}-1+pw^{p-1}$, which is spanned by $w_{x}$. In fact, we obtain
\begin{align*}
\check{\mathcal{A}}_{3}&=-\Delta^{\Gamma}f_{1}w_{x}-(k_{1}^{2}+k_{2}^{2})f_{1}w_{x}+\check{\mathcal{C}}_{3}
+g(f_{0}).
\end{align*}
Here \,$g(f_{0})$\, is a  smooth bounded function independent of the unknown parameters $f_{1}$\, and $f_{2}$.\\
Since the term $\hat{\mathcal{A}}_{3}$ is even in the variable $x$ and $w_{x}$ is odd in the variable $x$, there holds
\begin{align*}
\int_{\mathbb{R}}\hat{\mathcal{A}}_{3}w_{x}\,\mathrm{d}x=0.
\end{align*}
On the other hand, using the same arguments as in (\ref{orthogonalityw2}), we get
\begin{align}\label{orthogonalityw34}
\begin{aligned}
\int_{\mathbb{R}}\check{\mathcal{A}}_{3}w_{x}\,\mathrm{d}x&=-\Delta^{\Gamma}f_{1}\int_{\mathbb{R}}w_{x}^{2}\,\mathrm{d}x
-(k_{1}^{2}+k_{2}^{2})f_{1}\int_{\mathbb{R}}w_{x}^{2}\,\mathrm{d}x+\int_{\mathbb{R}}\check{\mathcal{C}}_{3}w_{x}\,\mathrm{d}x
+\int_{\mathbb{R}}gw_{x}\,\mathrm{d}x
\\
&=-\big[\Delta^{\Gamma}f_{1}+(k_{1}^{2}+k_{2}^{2})f_{1}\big]\int_{\mathbb{R}}w_{x}^{2}\,\mathrm{d}x
+\mathrm{d}_{1}(f_{0},e,\varphi_{1},\varphi_{2}).
\end{aligned}
\end{align}
Here $\mathrm{d}_{1}$ is a smooth bounded function independent of the unknown parameters $f_{1},f_{2}$ and the
correction terms $\varphi_{3}$ and is Lipschitz continuous with respect to its parameters.
By setting
\begin{align}
b_{1}=\int_{\mathbb{R}}w_{x}^{2}\,\mathrm{d}x
\end{align}
we derive from (\ref{orthogonalityw34}) that
\begin{align*}
\int_{\mathbb{R}}\check{\mathcal{A}}_{3}w_{x}\,\mathrm{d}x=0,
\end{align*}
is equivalent to the following differential equation
\begin{align}\label{conditionf1}
\Delta^{\Gamma}f_{1}+(k_{1}^{2}+k_{2}^{2})f_{1}=\frac{\mathrm{d}_{1}}{b_{1}}  \quad\mbox{on }\Gamma.
\end{align}
Now, we also choose $f_1$ by combining the equation (\ref{conditionf1})
with the boundary condition (\ref{boundaryf1})
with the help of the non-degeneracy condition (\ref{Geometriceigenvalueproblem}).
The solution to (\ref{defitionsolution34}) can be expressed as
\begin{align}\label{solution3}
\varphi_{3}(x,z)=\psi_{31}(x,\var z)+\psi_{32}(x,\var z),
\end{align}
where $\psi_{31}(x,\var z)$ is an odd function in the variable $x$ and $\psi_{32}(x,\var z)$ is an even function in the variable $x$.
The component in $\psi_{31}(x,\var z)$ and $\psi_{32}(x,\var z)$
are independent of the parameter $f_{2}$.

\begin{Remark}\label{remarkinnerapproximation}
The recurrence procedure we described above is the same as the arguments in \cite{wangweiyang}.
Whence, we can improve the local approximate solution to solve the first equation in (\ref{problemafterscaling0})
up to $O(\ve^m)$ for any positive integer $m$.
\end{Remark}

\section{Boundary Correction Layers and Further Improvement of Approximations}\label{section4}
\setcounter{equation}{0}

For any given $(f_{2}, e)\in\digamma$,
we have the inner expansion $\mathcal{V}$ in (\ref{basicapproximation}).
This approximation in general does not satisfy the boundary condition in (\ref{problemafterscaling0}).
In order to improve the approximation,
we need to write problem (\ref{problemafterscaling0})
in local coordinates in the neighborhood of $\partial\Gamma\in\R^3$,
and then add boundary correction layers to the inner expansion (\ref{basicapproximation}).
In fact, by recalling the scaled modified Fermi coordinates
with variables $\theta,\, \eta,\, s$ in Lemma \ref{lemma2.2}
and also the translated variable $x$ in (\ref{coordinatesxz}),
we will extend the local form of problem (\ref{problemafterscaling0})
to the infinite strip $\mathbf{\Lambda}$ in $\mathbb{R}^3$  with notations
\begin{align}
\begin{aligned}\label{lambda00}
\mathbf{\Lambda}\,=\,&\,\{\ (x,\theta,\eta):\, x\in \mathbb{R},\,
(\theta,\eta)\in\partial\mathbb{D}\times\mathbb{R}_{+}\},
\\
\partial \mathbf{\Lambda}\,=\,&\,\{\ (x,\theta,\eta):\, x\in \mathbb{R},\, \theta\in\partial\mathbb{D},\,\eta=0\ \}.
\end{aligned}
\end{align}

\subsection{Local formulation of the scaled problem near boundary}
Here is the local form of problem (\ref{problemafterscaling0}) in the neighborhood of $\partial\Gamma_\ve$:
\begin{lem}\label{lemma2.9}
In terms of the coordinate system $(s,\theta,\eta)$, the equation in (\ref{problemafterscaling0}) is expressed as
\begin{align}\label{boundaryequation}
\begin{aligned}
u_{ss}\,&+\,\frac{1}{l_1^2(\var\theta)}u_{\theta\theta}\,+\,\frac{1}{l_2^2(\var\theta)}u_{\eta\eta}\,-\,u\,+\,u^{p}
+\widetilde{B}(u)\,=\,0,
\end{aligned}
\end{align}
where $\widetilde{B}(u)=\widetilde{B}_{0}(u)+\widetilde{B}_{1}(u)+\widetilde{B}_{2}(u)$ with components
\begin{align*}
\widetilde{B}_{0}(u)\,=\,&\,\varepsilon\bigg[-2\eta\frac{A(\var\theta)}{l^{4}_{2}(\var\theta)}u_{\eta\eta}
+\bigg(\frac{C(\var\theta)}{l^{2}_{1}(\var\theta)l^{2}_{2}(\var\theta)}-\frac{A(\var\theta)}{l^{4}_{2}(\var\theta)}\bigg)u_{\eta}\bigg]
-\varepsilon^{2}s(k^{2}_{1}+k^{2}_{2})u_{s}
\\
&+\varepsilon\bigg(\frac{E(\var\theta)}{l^{2}_{1}(\var\theta)l^{2}_{2}(\var\theta)}
-\frac{R(\var\theta)}{l^{4}_{1}(\var\theta)}\bigg)u_{\theta}+\varepsilon^{3}b(\var\theta,\varepsilon \eta)s^{2}u_{s},
\\
\\
\widetilde{B}_{1}(u)\,=\,&\,-\varepsilon\frac{I(\var\theta)}{l_{1}(\var\theta)l_{2}(\var\theta)}u_{\eta}
-\varepsilon\frac{F(\var\theta)}{l_{1}(\var\theta)l_{2}(\var\theta)}u_{\theta}
-2\varepsilon s \frac{I(\var\theta)}{l_{1}(\var\theta)l_{2}(\var\theta)}u_{s \eta}
-2\varepsilon s \frac{F(\var\theta)}{l_{1}(\var\theta)l_{2}(\var\theta)}u_{s \theta},
\\
\\
\widetilde{B}_{2}(u)\,=\,&\,\var^{2}a_{1}u_{\theta}+\var^{2}a_{2}u_{\eta}
+\var^{2}b_{1}u_{s \theta}+\var^{2}b_{2}u_{s \eta}+\widetilde{B}_{3}(u),
\end{align*}
we have defined for $\var\theta=\vartheta$
\begin{align}
\begin{aligned}\label{I}
l_1(\vartheta)=\Big|\frac{\partial\gamma_0(\var\theta,0)}{\partial\vartheta}\Big|>0,
\qquad
&l_2(\vartheta)=\Big|\frac{\partial\gamma_0(\var\theta,0)}{\partial\rho}\Big|>0,
\\
A(\vartheta)=\Big<\frac{\partial^2\gamma_0(\var\theta,0)}{\partial\rho^2},
\,  \frac{\partial\gamma_0(\var\theta,0)}{\partial\rho}\Big>,
\qquad
&I(\vartheta)=\Big<q_{1}(\var\theta,0),\,  \frac{\partial\gamma_0(\var\theta,0)}{\partial\rho}\Big>,
\\
C(\vartheta)=\Big<\frac{\partial^2\gamma_0(\var\theta,0)}{\partial\vartheta\partial\rho},
\,  \frac{\partial\gamma_0(\var\theta,0)}{\partial\vartheta}\Big>,
\qquad
&R(\vartheta)=\Big<\frac{\partial^2\gamma_0(\var\theta,0)}{\partial\vartheta^{2}},
\, \frac{ \partial\gamma_0(\var\theta,0)}{\partial\vartheta}\Big>,
\\
E(\vartheta)=\Big<\frac{\partial^2\gamma_0(\var\theta,0)}{\partial\vartheta\partial\rho},
\, \frac{\partial\gamma_0(\var\theta,0)}{\partial\rho}\Big>,
\qquad
&F(\vartheta)=\Big<q_{1}(\var\theta,0),  \frac{\partial\gamma_0(\var\theta,0)}{\partial\vartheta}\Big>,
\end{aligned}
\end{align}
and  the differential operator $\widetilde{B}_{3}(s,\var\theta,\var\eta)$ is of size $O(\varepsilon^{4})$.
Here, $a_{1},a_{2},b_{1},b_{2}$ are the coefficients of the operator $B_{2}$ in (\ref{delta}).

The boundary condition in (\ref{problemafterscaling0}) is recast as
\begin{align}
\begin{aligned}\label{boundarycondition}
\varepsilon I(\var\theta)su_{s}&-u_{\eta}
-\frac{2\varepsilon s}{l_2(\,\var\theta)^2\,}\Big<\frac{\partial\nu}{\partial\rho},\,
\frac{\partial\gamma_0(\var\theta,0)}{\partial\rho}\Big>u_{\eta}
+\varepsilon\frac{M(\var\theta)}{l^{2}_{1}(\var\theta)}su_{\theta}
\\
&+\varepsilon^{2}\frac{I(\var\theta)}{l^{2}_{2}(\var\theta)}\Big<\frac{\partial\nu}{\partial\rho},\,
\frac{\partial\gamma_0(\var\theta,0)}{\partial\rho}\Big>s^{2}u_{s}
+\varepsilon^{2}b_{4}(\var\theta)s^{2}u_{\eta}
\\
&+\varepsilon^{3}b_{5}(\var\theta)s^{3}u_{s}
+\varepsilon^{2}b_{6}(\var s,\var\theta)u_{\theta}+\widetilde{D}(u)=0,
\end{aligned}
\end{align}
where we have denoted
$$
M(\var\theta)=\Big<\frac{\partial\gamma_0(\var\theta,0)}{\partial\vartheta},\,\frac{\partial\nu}{\partial\rho}\Big>
+
\Big<\frac{\partial\gamma_0(\var\theta,0)}{\partial\rho},\,\frac{\partial\nu}{\partial\vartheta}\Big>.
$$
The differential operator $\widetilde{D}(\var s,\var\theta)$ is of size $O(\varepsilon^{4})$  and functions
$b_{4}(\var\theta),\,b_{5}(\var\theta),\,b_{6}(\var s,\var\theta)$ satisfy:
\begin{align*}
|b_{4}(\var\theta),\,b_{5}(\var\theta),\,b_{6}(\var s,\var\theta)|\,\leq\, C(1+|\var s|^{4}).
\end{align*}
\end{lem}
\proof
The reader can refer to Section 4 in \cite{Sakamoto} and the references therein for the details of the proof.
\qed

\medskip
As we have done in (\ref{coordinatesxz}),
we introduce a further changing of variables and define a new function $v(x,\theta,\eta)$ as follows
\begin{align}\label{coordinatexeta}
u(s,\theta,\eta)=v(x,\theta,\eta)\quad\mbox{with}\quad
x=s-\sum\limits_{i=0}^{2}\varepsilon^{i}f_{i}(\var\theta,\varepsilon\eta),
\quad \theta=\theta,\quad \eta=\eta.
\end{align}
We now want to express the problem in the new coordinates. Whence we need the following formulas
\begin{align*}
\begin{aligned}
u_{s}&=v_{x},
\qquad
u_{ss}=v_{xx},
\\
u_{s \theta}&=\Big(-\sum\limits_{i=0}^{2}\varepsilon^{i+1}\nabla_{\vartheta}f_{i}\Big)v_{xx}+v_{x\theta},
\qquad\quad
u_{s \eta}=\Big(-\sum\limits_{i=0}^{2}\varepsilon^{i+1}\nabla_{\rho}f_{i}\Big)v_{xx}+v_{x\eta},
\\
u_{\eta}&=\Big(-\sum\limits_{i=0}^{2}\varepsilon^{i+1}\nabla_{\rho}f_{i}\Big)v_{x}+v_{\eta},
\qquad\quad\quad
u_{\theta}=\Big(-\sum\limits_{i=0}^{2}\varepsilon^{i+1}\nabla_{\vartheta}f_{i}\Big)v_{x}+v_{\theta},
\\
u_{\theta\theta}&=\Big(\sum\limits_{i=0}^{2}\varepsilon^{i+1}\nabla_{\vartheta}f_{i}\Big)^{2}v_{xx}
+\Big(-\sum\limits_{i=0}^{2}\varepsilon^{i+2}\Delta_{\vartheta}f_{i}\Big)v_{x}
+\Big(-\sum\limits_{i=0}^{2}\varepsilon^{i+1}\nabla_{\vartheta}f_{i}\Big)v_{x\theta}+v_{\theta\theta},
\\
u_{\eta\eta}&=\Big(\sum\limits_{i=0}^{2}\varepsilon^{i+1}\nabla_{\rho}f_{i}\Big)^{2}v_{xx}
+\Big(-\sum\limits_{i=0}^{2}\varepsilon^{i+2}\Delta_{\rho}f_{i}\Big)v_{x}
+\Big(-\sum\limits_{i=0}^{2}\varepsilon^{i+1}\nabla_{\rho}f_{i}\Big)v_{x\eta}+v_{\eta\eta}.
\end{aligned}
\end{align*}
Using Lemma \ref{lemma2.9},
we get the local form of the problem (\ref{problemafterscaling0}) in the new coordinates
and then extend it to the infinite strip ${\mathbf\Lambda}$
\begin{align}\label{boundaryforv}
S(v)\equiv v_{xx}+\,\frac{1}{l_1^2(\var\theta)}v_{\theta\theta}+\,\frac{1}{l_2^2(\var\theta)}v_{\eta\eta}\,-\,v\,+\,v^{p}+\widetilde{B}(v)
\,=\,0\qquad \mbox{in} \,\,\mathbf{\Lambda},
\end{align}
where we have denoted
\begin{align}\label{definitionofBv}
\widetilde{B}(v)=&\widetilde{B}_{0}(v)+\widetilde{B}_{1}(v)+\widetilde{B}_{2}(v),
\end{align}
\begin{align*}
\widetilde{B}_{0}(v)=&\frac{1}{l^{2}_{2}(\var\theta)}\Bigg[\Big(\sum\limits_{i=0}^{2}\varepsilon^{i+1}\nabla_{\rho}f_{i}\Big)^{2}v_{xx}
+\Big(-\sum\limits_{i=0}^{2}\varepsilon^{i+2}\Delta_{\rho}f_{i}\Big)v_{x}
+\Big(-\sum\limits_{i=0}^{2}\varepsilon^{i+1}\nabla_{\rho}f_{i}\Big)v_{x\eta}\Bigg]
\\
&-2\varepsilon\eta\frac{A(\var\theta)}{l^{4}_{2}(\var\theta)}\Bigg[\Big(\sum\limits_{i=0}^{2}\varepsilon^{i+1}\nabla_{\rho}f_{i}\Big)^{2}v_{xx}
+\Big(-\sum\limits_{i=0}^{2}\varepsilon^{i+2}\Delta_{\rho}f_{i}\Big)v_{x}
+\Big(-\sum\limits_{i=0}^{2}\varepsilon^{i+1}\nabla_{\rho}f_{i}\Big)v_{x\eta}+v_{\eta\eta}\Bigg]
\\
&+\varepsilon\bigg(\frac{C(\var\theta)}{l^{2}_{1}(\var\theta)l^{2}_{2}(\var\theta)}-\frac{A(\var\theta)}{l^{4}_{2}(\var\theta)}\bigg)
\Bigg[\,\Big(-\sum\limits_{i=0}^{2}\varepsilon^{i+1}\nabla_{\rho}f_{i}\Big)v_{x}+v_{\eta}\,\Bigg]
-\varepsilon^{2}\Bigg(x+\sum\limits_{i=0}^{2}\varepsilon^{i}f_{i}\Bigg)(k^{2}_{1}+k^{2}_{2})v_{x}
\\
&+\varepsilon^{2}\frac{1}{l^{2}_{1}(\var\theta)}\Bigg[\Big(\sum\limits_{i=0}^{2}\varepsilon^{i}\nabla_{\vartheta}f_{i}\Big)^{2}v_{xx}
+\Big(-\sum\limits_{i=0}^{2}\varepsilon^{i}\Delta_{\vartheta}f_{i}\Big)v_{x}
+\Big(-\sum\limits_{i=0}^{2}\varepsilon^{i}\nabla_{\vartheta}f_{i}\Big)v_{x\theta}\Bigg]
\\
&+\varepsilon\bigg(\frac{E(\var\theta)}{l^{2}_{1}(\var\theta)l^{2}_{2}(\var\theta)}-\frac{R(\var\theta)}{l^{4}_{2}(\var\theta)}\bigg)
\Bigg[\Big(-\sum\limits_{i=0}^{2}\varepsilon^{i+1}\nabla_{\vartheta}f_{i}\Big)v_{x}+v_{\theta}\Bigg]
+\varepsilon^{3}b(\var\theta,\var\eta)\Bigg(x+\sum\limits_{i=0}^{2}\varepsilon^{i}f_{i}\Bigg)v_{x}\,,
\end{align*}
\begin{align*}
\widetilde{B}_{1}(v)=&-\varepsilon\frac{I(\var\theta)}{l_{1}(\var\theta)l_{2}(\var\theta)}
\Big[\Big(-\sum\limits_{i=0}^{2}\varepsilon^{i+1}\nabla_{\rho}f_{i}\Big)v_{x}+v_{\eta}\Big]
-\varepsilon\frac{F(\var\theta)}{l_{1}(\var\theta)l_{2}(\var\theta)}
\Big[\Big(-\sum\limits_{i=0}^{2}\varepsilon^{i+1}\nabla_{\vartheta}f_{i}\Big)v_{x}+v_{\theta}\Big]
\\
&-2\varepsilon\Big(x+\sum\limits_{i=0}^{2}\varepsilon^{i}f_{i}\Big) \frac{I(\var\theta)}{l_{1}(\var\theta)l_{2}(\var\theta)}
\Big[\Big(-\sum\limits_{i=0}^{2}\varepsilon^{i+1}\nabla_{\rho}f_{i}\Big)v_{xx}+v_{x\eta}\Big]
\\
&-2\varepsilon\Big(x+\sum\limits_{i=0}^{2}\varepsilon^{i}f_{i}\Big)\frac{F(\var\theta)}{l_{1}(\var\theta)l_{2}(\var\theta)}
\Big[\Big(-\sum\limits_{i=0}^{2}\varepsilon^{i+1}\nabla_{\vartheta}f_{i}\Big)v_{xx}+v_{x \theta}\Big],
\end{align*}
\begin{align*}
\widetilde{B}_{2}(v)=\,&\var^{2}a_{1}\Big[\Big(-\sum\limits_{i=0}^{2}\varepsilon^{i+1}\nabla_{\vartheta}f_{i}\Big)v_{x}+v_{\theta}\Big]
+\var^{2}a_{2}\Big[\Big(-\sum\limits_{i=0}^{2}\varepsilon^{i+1}\nabla_{\rho}f_{i}\Big)v_{x}+v_{\eta}\Big]
\\
&+\var^{2}b_{1}\Big[\Big(-\sum\limits_{i=0}^{2}\varepsilon^{i+1}\nabla_{\vartheta}f_{i}\Big)v_{xx}+v_{x \theta}\Big]
+\var^{2}b_{2}\Big[\Big(-\sum\limits_{i=0}^{2}\varepsilon^{i+1}\nabla_{\rho}f_{i}\Big)v_{xx}+v_{x\eta}\Big]+\widetilde{B}_{3}(v).
\end{align*}
The boundary condition is
\begin{align}
\begin{aligned}\label{boundaryconditionv}
\varepsilon\Big(x+&\sum\limits_{i=0}^{2}\varepsilon^{i}f_{i}\Big)I(\var\theta)v_{x}
+\Big(\sum\limits_{i=0}^{2}\varepsilon^{i+1}\nabla_{\rho}f_{i}\Big)v_{x}-\,v_{\eta}
\\
&-\frac{2\varepsilon}{l_2^2(\var\theta)}\Big< \frac{\partial\nu}{\partial\rho},
\frac{\partial\gamma_0(\var\theta,0)}{\partial\rho}\Big>\Big(x+\sum\limits_{i=0}^{2}\varepsilon^{i}f_{i}\Big)
\Big[\Big(-\sum\limits_{i=0}^{2}\varepsilon^{i+1}\nabla_{\rho}f_{i}\Big)v_{x}+v_{\eta}\Big]
\\
&+\varepsilon^{2}\frac{I(\var\theta)}{l^{2}_{2}(\var\theta)}\Big<\frac{\partial\nu}{\partial\rho},
\frac{\partial\gamma_0(\var\theta,0)}{\partial\rho}\Big>\Big(x+\sum\limits_{i=0}^{2}\varepsilon^{i}f_{i}\Big)^{2}v_{x}
+\varepsilon^{3}b_{5}(\var\theta)\Big(x+\sum\limits_{i=0}^{2}\varepsilon^{i}f_{i}\Big)^{3}v_{x}
\\
&+\varepsilon^{2}b_{4}(\var\theta)\Big(x+\sum\limits_{i=0}^{2}\varepsilon^{i}f_{i}\Big)^{2}
\Big[\Big(-\sum\limits_{i=0}^{2}\varepsilon^{i+1}\nabla_{\rho}f_{i}\Big)v_{x}+v_{\eta}\Big]
\\
&+\varepsilon\frac{M(\var\theta)}{l^{2}_{1}(\var\theta)}\Big(x+\sum\limits_{i=0}^{2}\varepsilon^{i}f_{i}\Big)
\Big[\Big(-\sum\limits_{i=0}^{2}\varepsilon^{i+1}\nabla_{\vartheta}f_{i}\Big)v_{x}+v_{\theta}\Big]
\\
&+\varepsilon^{2}b_{6}(x,\var\theta)\Big[\Big(-\sum\limits_{i=0}^{2}
\varepsilon^{i+1}\nabla_{\vartheta}f_{i}\Big)v_{x}+v_{\theta}\Big]+\widetilde{D}(v)=\,0\, \qquad  \mbox{on} \,\, \partial\mathbf{\Lambda}.
\end{aligned}
\end{align}

\subsection{Boundary correction layers}
In this subsection, we will find some boundary correction layer terms,
say $\phi_i$'s, by also the recurrence method.
The method is basically the same as that in \cite{weiyang1} and \cite{weiyang2}.

By recalling (\ref{basicapproximation}), we take
\begin{align}\label{u1}
u_1\equiv\mathcal{V}=w(x)+\varepsilon eZ+\sum\limits_{i=2}^{3}\varepsilon^{i}\varphi_{i},
\end{align}
as the first approximate solution of the problem (\ref{boundaryforv}) and (\ref{boundaryconditionv}) on $\mathbf{\Lambda}$.
Then we compute
\begin{align}\label{definitionofS(u_{1})}
S(u_{1})=\varepsilon^{3}\Delta^{\Gamma}e Z+\varepsilon\lambda_{0}eZ+\widetilde{B}_{3}(u_1)
\quad  \mbox{in }\, \mathbf{\Lambda}.
\end{align}
 On the boundary, the errors become
\begin{align*}
\Phi\,\equiv\,\ve\Phi_1\,+\,\ve^2\Phi_2\,+\,\ve^3\Phi_3\,+\,\widetilde{D}(u_{1})=0
\quad \mbox{on }\,\, \partial\mathbf{\Lambda},
\end{align*}
where we have denoted
\begin{align*}
\Phi_1\,=\,&\,\Big(I(\var\theta)x+I(\var\theta)f_{0}+\nabla_{\rho}f_{0}\Big)w_x,
\\
\\
\Phi_2\,=\,&\,\nabla_{\rho}f_{1}w_{x}
+I(\var\theta)f_{1}w_{x}-\nabla_{\rho}e\,Z+\nabla_{\rho}f_{0}eZ_{x}
\\
&+I(\var\theta)\Big(x+f_{0}\Big)eZ_{x}+\frac{2}{l_2^2(\,\var\theta)\,}\Big<\frac{\partial\nu}{\partial\rho},\,
\frac{\partial\gamma_{0}}{\partial\rho}\Big>\Big(x+f_{0}\Big)\nabla_{\rho}f_{0}w_{x}
\\
&+\frac{I(\var\theta)}{l^{2}_{2}(\var\theta)}\Big<\frac{\partial\nu}{\partial\rho},\,
\frac{\partial\gamma_0(\var\theta,0)}{\partial\rho}\Big>\Big(x+f_{0}\Big)^{2}w_{x}
-\frac{M(\var\theta)}{l^{2}_{1}(\var\theta)}\Big(x+f_{0}\Big)\nabla_{\vartheta}f_{0}w_{x},
\end{align*}
\begin{align*}
\Phi_3\,=\,&\,\nabla_{\rho}f_{2}w_{x}+I(\var\theta)f_{2}w_{x}
+\nabla_{\rho}f_{1}eZ_{x}+I(\var\theta)f_{1}eZ_{x}-\varphi_{2,\eta}
\,+\,\varepsilon^{3}b_{5}(\var\theta)\Big(x+f_{0}\Big)^{3}w_{x}
\\
&+\frac{2}{l_2^2(\,\var\theta)\,}\Big<\frac{\partial\nu}{\partial\rho},
\frac{\partial\gamma_{0}}{\partial\rho}\Big>\Big[\Big(x+f_{0}\Big)\Big(\nabla_{\rho}f_{1}w_{x}
+\nabla_{\rho}f_{0}eZ_{x}-\nabla_{\rho}eZ\Big)+f_{1}\nabla_{\rho}f_{0}w_{x}\Big]
\\
&+\frac{I(\var\theta)}{l^{2}_{2}(\var\theta)}\Big<\frac{\partial\nu}{\partial\rho},\,
\frac{\partial\gamma_0(\var\theta,0)}{\partial\rho}\Big>\Big[\,\Big(x+f_{0}\Big)^{2}eZ_{x}+2f_{1}xw_{x}\,\Big]
\,-\,b_{4}(\var\theta)\Big(x+f_{0}\Big)^{2}\nabla_{\rho}f_{0}w_{x}
\\
&\,-\frac{M(\var\theta)}{l^{2}_{1}(\var\theta)}\Big[\,\Big(x+f_{0}\Big)\Big(\nabla_{\vartheta}f_{1}w_{x}
+\nabla_{\vartheta}f_{0}eZ_{x}-\nabla_{\vartheta}eZ\Big)+f_{1}\nabla_{\vartheta}f_{0}w_{x}\,\Big]
\\
&-b_{6}(x,\var\theta)\nabla_{\vartheta}f_{0}w_{x}.
\end{align*}

To improve the approximation,
we need to cancel the terms $\Phi_i$'s by adding boundary correction terms to
the inner approximate expansion in (\ref{u1}).
This can be done step by step.
On the boundary it is natural to take
\begin{align}\label{boundaryf0}
\nabla_{\rho}f_{0}+I(\var\theta)f_{0}=0,
\end{align}
which will lead to the cancelation of $(I(\var\theta)f_{0}+\nabla_{\rho}f_{0}\big)w_x$
in the components of the error of first order of $\var$.
Recall that by combining of (\ref{conditionf0}) and (\ref{boundaryf0})
we have chosen $f_{0}\equiv 0$ in Section \ref{section3}.
On the other hand, to cancel the first order term like $\varepsilon I(\var\theta)xw_x$
on the boundary, we shall introduce a boundary layer term, say $\phi_1$.
Indeed, we first introduce the term
$$
\tilde{\phi} (x,\theta,\eta)=\var\, A_{0}(\var\theta)\,b\big(\sqrt{\lambda_{0}}\,l_{2}(\var \theta)\eta\big)Z(x),
$$
with $\tilde{\phi} (x,\theta,\eta)$ satisfying
\begin{align*}
\tilde{\phi}_{xx}+\frac{1}{l^{2}_{1}(\var\theta)}\tilde{\phi}_{\theta\theta}
+\frac{1}{l^{2}_{2}(\var\theta)}\tilde{\phi}_{\eta\eta}
-\tilde{\phi} +pw^{p-1}\tilde{\phi}=O(\var^{3}) \qquad &\mbox{in} \,\, \mathbf{\Lambda},
\\
\frac{\partial \tilde{\phi}}{\partial \eta}=-c_{0}(\var\theta)Z,
\qquad \mbox{on}\,\,\partial\mathbf{\Lambda},
\end{align*}
where $c_0(\var\theta)$ is a function in the parameter $\theta$ of the form
\begin{align}\label{definitionofCzeroandCone}
c_0(\var\theta)=I(\var\theta)\int_\mathbb{R} xw_xZ\,\mathrm{d}x.
\end{align}
In fact, we can choose
\begin{align}\label{definitionofphitilde1}
\tilde{\phi} (x,\theta,\eta)
\,=\,\varepsilon\frac{-c_{0}(\var \theta)}{\sqrt{\lambda_{0}}\,l_{2}(\var \theta)}
\sin\big(\sqrt{\lambda_{0}}\,l_{2}(\var \theta)\eta\big)Z(x)
\,\equiv\,\var\phi_{11}(x,\theta,\eta).
\end{align}
Then, by Corollary \ref{cor1}, there exists a   unique solution (denoted by $\phi_{12}$) of the following problem
\begin{align*}
\phi_{12,xx}+\frac{1}{l^{2}_{1}(\var\theta)}\phi_{12,\theta\theta}+\frac{1}{l^{2}_{2}(\var\theta)}\phi_{12,\eta\eta}
-\phi_{12} +pw^{p-1}\phi_{12}=0 \qquad &\mbox{in} \,\, \mathbf{\Lambda},
\\
\frac{\partial \phi_{12}}{\partial \eta}=I(\var\theta)x w_{x}-c_{0}(\var\theta)Z
\qquad& \mbox{on}\,\,\partial\mathbf{\Lambda}.
\end{align*}
Moreover, $\phi_{12}$ is even in the variable $x$.
By choosing a smooth cut-off function $\chi_0$ in the form
\begin{align}\label{cut-off}
\chi_0(t)=1\quad\mbox{if}\quad |t|<1,
\qquad
\chi_0(t)=0\quad\mbox{if}\quad |t|>2,
\end{align}
we can set the {\bf first boundary layer term} by
\begin{align}\label{chi0}
\phi_1= \varepsilon\,\chi_{0}\,(\varepsilon \eta)\,\Psi_{1}(x,\theta,\eta)
\quad\mbox{with}\quad
\Psi_{1}(x,\theta,\eta)=\phi_{11} (x,\theta,\eta)+\phi_{12} (x,\theta,\eta).
\end{align}
Thus we finish the first step.

To proceed the second step,
let $u_2=u_{1}+\phi_1$ be the second approximate solution.
We again compute the new error
\begin{align}\label{definitionofS(u_{1}+phione)}
\begin{aligned}
S(u_{2})=S(u_{1})&+2\,\var^2\,\frac{1}{l^{2}_{2}(\var\theta)}\chi_0'\,\Psi_{1,\eta}
+\varepsilon^3\frac{1}{l^{2}_{2}(\var\theta)}\chi_0''\,\Psi_{1}+\frac{p(p-1)}{2}w^{p-2}(\varepsilon eZ+\phi_{1})^{2}
\\
&+\varepsilon^{2}\bigg(\frac{C(\var\theta)}{l_{1}^{2}(\var\theta)l_{2}^{2}(\var\theta)}
-\frac{A(\var\theta)}{l_{2}^{4}(\var\theta)}-\frac{I(\var\theta)}{l_{1}(\var\theta)l_{2}(\var\theta)}\bigg)\chi_{0}\Psi_{1,\eta}
\\
&-\varepsilon^{2}\bigg(\frac{1}{l_{2}^{2}(\var\theta)}
+2x\frac{I(\var\theta)}{l_{1}(\var\theta)l_{2}(\var\theta)}\bigg)\chi_{0}\Psi_{1,x\eta}
\\
&+\mathcal{C}_{0}(\phi_{11})+L_{0}(\phi_{1})
+N_0(\phi_{1})\qquad \mbox{in}\,\,\mathbf{\Lambda},
\end{aligned}
\end{align}
where $\mathcal{C}_{0}(\phi_{11}),\,L_{0}(\phi_{1})$ and $N_0(\phi_{1})$ are of size $O(\varepsilon^{3})$  and $S(u_{1})$ is defined in (\ref{definitionofS(u_{1})}),
\begin{align*}
\mathcal{C}_{0}(\phi_{11})=S(\phi_{11})&+\var^{2}\bigg(\frac{E(\var\theta)}{l^{2}_{1}(\var\theta)l^{2}_{2}(\var\theta)}
-\frac{R(\var\theta)}{l^{4}_{2}(\var\theta)}-\frac{F(\var\theta)}{l_{1}(\var\theta)l_{2}(\var\theta)}\bigg)\phi_{11,\theta}
\\
&-2\var^{2}x\frac{F(\var\theta)}{l_{1}(\var\theta)l_{2}(\var\theta)}\phi_{11,x\theta},
\\
L_{0}(\phi_{1})=\widetilde{B}(\phi_1)&-\varepsilon^{2}\bigg(\frac{C(\var\theta)}{l_{1}^{2}(\var\theta)l_{2}^{2}(\var\theta)}
-\frac{A(\var\theta)}{l_{2}^{4}(\var\theta)}-\frac{I(\var\theta)}{l_{1}(\var\theta)l_{2}(\var\theta)}\bigg)\chi_{0}\Psi_{1,\eta}
\\
&+\varepsilon^{2}\bigg(\frac{1}{l_{2}^{2}(\var\theta)}+2x\frac{I(\var\theta)}{l_{1}(\var\theta)l_{2}(\var\theta)}\bigg)\chi_{0}\Psi_{1,x\eta}\,,
\\
N_0(\phi_{1})=(u_{1}+&\phi_1)^p-u_{1}^p-pw^{p-1}\phi_{1}-\frac{p(p-1)}{2}w^{p-2}(\varepsilon eZ+\phi_{1})^{2}.
\end{align*}
Note that $f_{0}\equiv 0$. On the boundary, the  new error becomes
\begin{align*}
&\varepsilon^{2}\bigg(\nabla_{\rho}f_{1}\,w_{x}+I(\var\theta)f_{1}\,w_{x}+I(\var\theta)e\,xZ_{x}-\nabla_{\rho}e\,Z\bigg)
+\varepsilon^{2}I(\var\theta)x\,\phi_{12,x}
\\
&+\varepsilon^{2}\frac{I(\var\theta)}{l^{2}_{2}(\var\theta)}\Big<\frac{\partial\nu}{\partial\rho},
\frac{\partial\gamma_0(\var\theta,0)}{\partial\rho}\Big>x^{2}w_{x}
+\varepsilon^{2}\frac{M(\var\theta)}{l^{2}_{1}(\var\theta)}\phi_{12,\theta}
\\
&+\varepsilon^{3}\Big[\nabla_{\rho}f_{2}w_{x}+I(\var\theta)f_{2}w_{x}
+\nabla_{\rho}f_{1}(eZ_{x}+\phi_{12,x})+I(\var\theta)f_{1}(eZ_{x}+\phi_{12,x})-\varphi_{2,\eta}\Big]
\\
&+\varepsilon^{3}\frac{2x}{l_2^2(\,\var\theta)\,}\Big<\frac{\partial\nu}{\partial\rho},
\frac{\partial\gamma_{0}}{\partial\rho}\Big>\Big(\nabla_{\rho}f_{1}w_{x}-\nabla_{\rho}eZ\Big)
\\
&+\varepsilon^{3}\frac{I(\var\theta)}{l^{2}_{2}(\var\theta)}\Big<\frac{\partial\nu}{\partial\rho},
\frac{\partial\gamma_0(\var\theta,0)}{\partial\rho}\Big>\bigg[x^{2}(eZ_{x}+\phi_{12,x})+2f_{1}xw_{x}\bigg]
\\
&+\varepsilon^{3}b_{4}(\var\theta)x^{2}\Psi_{1,\eta}
+\varepsilon^{3}b_{5}(\var\theta)x^{3}w_{x}-\varepsilon^{3}b_{6}(x,\var\theta)\nabla_{\vartheta}f_{0}w_{x}
\\
&-\varepsilon^{3}\frac{M(\var\theta)}{l^{2}_{1}(\var\theta)}x\big(\nabla_{\vartheta}f_{1}w_{x}
-\nabla_{\vartheta}eZ\big)+\widetilde{D}(u_{2})=0  \quad \mbox{on}\,\, \partial\mathbf{\Lambda}.
\end{align*}
In order to cancel the error terms of order $O(\varepsilon^{2})$ on the boundary,
we do the same thing as that in the first step.
We first set
\begin{align}\label{boundaryf1}
\nabla_{\rho}f_{1}+I(\var\theta)f_{1}=0,
\end{align}
and
\begin{align}\label{boundarye}
-\int_{\R}\Big[\,I(\var\theta)e\,x Z_{x}-\nabla_{\rho}e\,Z\,\Big]\,Z\,\mathrm{d}x
\,=\,\nabla_{\rho}e+\frac{1}{2}I(\var\theta)e\,=\,0.
\end{align}
Then we choose
\begin{align}\label{definitionofphitilde2}
\hat{\phi} (x,\theta,\eta)
\,=\,\varepsilon^2\frac{-c_{1}(\var \theta)-c_{2}(\var \theta)}{\sqrt{\lambda_{0}}\,l_{2}(\var \theta)}
\sin\big(\sqrt{\lambda_{0}}\,l_{2}(\var \theta)\eta\big)Z(x)
\equiv\var^2\phi_{21}(x,\theta,\eta),
\end{align}
where $c_1(\var\theta)$ and $c_2(\var\theta)$ are functions in the parameter $\theta$ of the form
\begin{align}\label{definitionofConeandCone}
c_1(\var\theta)=I(\var\theta)\int_\mathbb{R} x\phi_{12,x}Z\,\mathrm{d}x,
\qquad
c_2(\var\theta)=-\frac{M(\var\theta)}{l^{2}_{1}(\var\theta)}\int_\mathbb{R} \phi_{12,\theta}Z\,\mathrm{d}x.
\end{align}
Now, by Lemma \ref{lemma0.7}, there exists a   unique solution (denoted by $\phi_{22}$) of the following problem
\begin{align*}
\phi_{22,xx}+\frac{1}{l^{2}_{1}(\var\theta)}\phi_{22,\theta\theta}+\frac{1}{l^{2}_{2}(\var\theta)}\phi_{22,\eta\eta}
-{\phi}_{22} +pw^{p-1}\phi_{22}=h_{1}\,\, \qquad  \mbox{in} \,\, \mathbf{\Lambda},
\qquad\qquad
\\
\\
\frac{\partial \phi_{22}}{\partial \eta}\,=\,\,
I(\var\theta)e\,x Z_{x}-\nabla_{\rho}e\,Z
+I(\var\theta)x\phi_{12,x}-c_{1}(\var\theta)Z
-\frac{M(\var\theta)}{l^{2}_{1}(\var\theta)}\phi_{12,\theta}-c_{2}(\var\theta)Z  \,\,\qquad \, \mbox{on}\,\,\partial\mathbf{\Lambda},
\end{align*}
where $h_1$ is a function in the form
\begin{align*}
h_{1}=&-2\,\frac{1}{l^{2}_{2}(\var\theta)}\chi_0'\,\Psi_{1,\eta}-\frac{p(p-1)}{2}w^{p-2}(\varepsilon eZ+\phi_{1})^{2}
\\
&+\bigg(\frac{C(\var\theta)}{l_{1}^{2}(\var\theta)l_{2}^{2}(\var\theta)}
-\frac{A(\var\theta)}{l_{2}^{4}(\var\theta)}-\frac{I(\var\theta)}{l_{1}(\var\theta)l_{2}(\var\theta)}\bigg)\chi_{0}\Psi_{1,\eta}
-\bigg(\frac{1}{l_{2}^{2}(\var\theta)}
+2x\frac{I(\var\theta)}{l_{1}(\var\theta)l_{2}(\var\theta)}\bigg)\chi_{0}\Psi_{1,x\eta}.
\end{align*}
Moreover $\phi_{22}$ is even in the variable $x$.
By recalling the cut-off function in (\ref{cut-off}),
we set the {\bf second boundary layer term} by
\begin{align}
\phi_{2}= \varepsilon^{2}\chi_{0}(\varepsilon \eta)\,\Psi_{2}(x,\theta,\eta)
\quad\mbox{with}\quad
\Psi_{2}(x,\theta,\eta)=\phi_{21} (x,\theta,\eta)+\phi_{22} (x,\theta,\eta).
\end{align}

For the completeness, we here give the third step although it is the same as we have done in the above.
Let $u_{3}=u_{2}+\phi_{2}$ be the second approximate solution. We compute the new error
\begin{align*}
S(u_{3})=&S(u_{2})+\varepsilon^{3}\frac{1}{l^{2}_{2}(\var\theta)}\chi_0''\Psi_{1}
+2\,\varepsilon^3\,\frac{1}{l^{2}_{2}(\var\theta)}\chi_0'\,\Psi_{2,\eta}
+\varepsilon^4\frac{1}{l^{2}_{2}(\var\theta)}\chi_0''\,\Psi_{2}
\\
&+p(p-1)w^{p-2}\varepsilon eZ\phi_{2}+\varepsilon^{3}\bigg(\frac{C(\var\theta)}{l_{1}^{2}(\var\theta)l_{2}^{2}(\var\theta)}
-\frac{A(\var\theta)}{l_{2}^{4}(\var\theta)}-\frac{I(\var\theta)}{l_{1}(\var\theta)l_{2}(\var\theta)}\bigg)\chi_{0}\Psi_{2,\eta}
\\
&-\varepsilon^{3}\bigg(\frac{1}{l_{2}^{2}(\var\theta)}
+2x\frac{I(\var\theta)}{l_{1}(\var\theta)l_{2}(\var\theta)}\bigg)\chi_{0}\Psi_{2,x\eta}
-4\varepsilon^{3}\eta\frac{A(\var\theta)}{l_{2}^{4}(\var\theta)}\chi_{0}'\Psi_{1,\eta}
\\
&-\var^{3}x\big(k_{1}^{2}+k_{2}^{2}\big)\chi_{0}\Psi_{1,x}+\mathcal{C}_{0}+\mathcal{C}_{1}+L_{1}(\phi_{2})+N_{1}(\phi_{2})
\qquad \mbox{in}\,\,\mathbf{\Lambda},
\end{align*}
where $\mathcal{C}_{1}(\phi_{21}),\,L_{0}(\phi_{2})$ and $N_0(\phi_{2})$ are of size $O(\varepsilon^{4})$  and $S(u_{2})$ is defined in (\ref{definitionofS(u_{1}+phione)}),
\begin{align*}
\mathcal{C}_{1}(\phi_{21})&=\var S(\phi_{21})+\var^{3}\bigg(\frac{E(\var\theta)}{l^{2}_{1}(\var\theta)l^{2}_{2}(\var\theta)}
-\frac{R(\var\theta)}{l^{4}_{2}(\var\theta)}-\frac{F(\var\theta)}{l_{1}(\var\theta)l_{2}(\var\theta)}\bigg)\phi_{21,\theta}
\\
&-2\var^{3}x\frac{F(\var\theta)}{l_{1}(\var\theta)l_{2}(\var\theta)}\phi_{21,x\theta},
\\
L_{1}(\phi_{2})=&\widetilde{B}(\phi_2)+L_{0}(\phi_{1})-\varepsilon^{3}\bigg(\frac{C(\var\theta)}{l_{1}^{2}(\var\theta)l_{2}^{2}(\var\theta)}
-\frac{A(\var\theta)}{l_{2}^{4}(\var\theta)}-\frac{I(\var\theta)}{l_{1}(\var\theta)l_{2}(\var\theta)}\bigg)\chi_{0}\Psi_{2,\eta}
\\
&+\varepsilon^{3}\bigg(\frac{1}{l_{2}^{2}(\var\theta)}+2x\frac{I(\var\theta)}{l_{1}(\var\theta)l_{2}(\var\theta)}\bigg)\chi_{0}\Psi_{2,x\eta}
\\
&+4\varepsilon^{3}\eta\frac{A(\var\theta)}{l_{2}^{4}(\var\theta)}\chi_{0}'\Psi_{1,\eta}-\var^{4}x(k_{1}^{2}+k_{2}^{2})\chi_{0}\Psi_{2,x},
\\
N_{1}(\phi_{2})=&(u_{2}+\phi_2)^p-u_{2}^p-pw^{p-1}\phi_{2}-p(p-1)w^{p-2}\varepsilon eZ\phi_{2}.
\end{align*}
On the boundary, the new error becomes
\begin{align*}
\varepsilon^{3}\big(\nabla_{\rho}f_{2}&\,w_{x}+I(\var\theta)f_{2}\,w_{x}+I(\var\theta)x\,\varphi_{2,x}-\varphi_{2,\eta}\big)
\\
&+\varepsilon^{3}I(\var\theta)x\,\phi_{22,x}+\varepsilon^{3}b_{4}(\var\theta)x^{2}\Psi_{1,\eta}
+\varepsilon^{3}b_{5}(\var\theta)x^{3}w_{x}
\\
&-\varepsilon^{3}\frac{M(\var\theta)}{l^{2}_{1}(\var\theta)}x\big(\nabla_{\vartheta}f_{1}w_{x}
-\nabla_{\vartheta}eZ\big)
+\widetilde{D}(u_{3})=0  \quad \mbox{on}\,\, \partial\mathbf{\Lambda}.
\end{align*}
In order to cancel the error terms of order $O(\varepsilon^{3})$ on the boundary,
we first set
\begin{align}\label{boundaryf2}
\nabla_{\rho}f_{2}+I(\var\theta)f_{2}+\mathfrak{R}(\theta)=0,
\end{align}
by introducing the term
\begin{align}\label{mathfrakR}
\mathfrak{R}(\theta)\,=\,
\frac{M(\var\theta)}{l^{2}_{1}(\var\theta)}\nabla_{\vartheta}e\int_\mathbb{R} xZw_{x}
\,\mathrm{d}x\Bigg/\int_\mathbb{R}w_{x}^{2}\,\mathrm{d}x.
\end{align}
Then we also introduce the term of the form
\begin{align}\label{definitionofphitilde3}
\check{\phi} (x,\theta,\eta)
\,=\,\var\frac{-c_{3}(\var \theta)}{\sqrt{\lambda_{0}}\,l_{2}(\var \theta)}
\sin\Big(\sqrt{\lambda_{0}}\,l_{2}(\var \theta)\,\eta\Big)Z(x)
\equiv\var\phi_{31}(x,\theta,\eta),
\end{align}
where $c_3(\var\theta)$ is a function of the form
\begin{align}\label{definitionofCthreeandCone}
c_3(\var\theta)=\int_\mathbb{R}c_{4}(x,\var\theta) Z\,\mathrm{d}x.
\end{align}
In the above, we have denoted
\begin{align}\label{defintR}
\begin{aligned}
c_{4}(x,\var\theta)=&I(\var\theta)x\,\varphi_{2,x}+I(\var\theta)x\,\phi_{22,x}-\varphi_{2,\eta}+b_{4}(\var\theta)x^{2}
\\
&+b_{5}(\var\theta)x^{3}-\frac{M(\var\theta)}{l^{2}_{1}(\var\theta)}\nabla_{\vartheta}f_{1}xw_{x},
\end{aligned}
\end{align}
Again, by Lemma \ref{lemma0.7}, there exists a   unique solution (denoted by $\phi_{32}$) of the following problem
\begin{align*} \label{phi32}
\phi_{32,xx}+\frac{1}{l^{2}_{1}(\var\theta)}\phi_{32,\theta\theta}+&\frac{1}{l^{2}_{2}(\var\theta)}\phi_{32,\eta\eta}-{\phi}_{32} +pw^{p-1}\phi_{32}=h_{2}\,\, \qquad  \mbox{in} \,\, \mathbf{\Lambda},
\\
\frac{\partial \phi_{32}}{\partial \eta}&=c_{4}(x,\var\theta)
-c_{3}(\var\theta)Z  \,\,\qquad \qquad  \mbox{on}\,\,\partial\mathbf{\Lambda},
\end{align*}
where $h_2$ is a function of the form
\begin{align*}
h_{2}=&-(\frac{C(\var\theta)}{l_{1}^{2}(\var\theta)l_{2}^{2}(\var\theta)}
-\frac{A(\var\theta)}{l_{2}^{4}(\var\theta)}-\frac{I(\var\theta)}{l_{1}(\var\theta)l_{2}(\var\theta)})\chi_{0}\Psi_{2,\eta}
-p(p-1)w^{p-2}\varepsilon eZ\phi_{2}
\\
&+\Big(\frac{1}{l_{2}^{2}(\var\theta)}+2x\frac{I(\var\theta)}{l_{1}(\var\theta)l_{2}(\var\theta)}\Big)\chi_{0}\Psi_{2,x\eta}
+4\eta\frac{A(\var\theta)}{l_{2}^{4}(\var\theta)}\chi_{0}'\Psi_{1,\eta}
-\frac{1}{l^{2}_{2}(\var\theta)}\chi_0''\Psi_{1}
\\
&+\mathcal{C}_{0}/\var^{3}+x(k_{1}^{2}+k_{2}^{2})\chi_{0}\Psi_{1,x}.
\end{align*}
Moreover $ \phi_{22}$ is even in the variable $x$.
We can set the {\bf third boundary layer term} by
\begin{align}
\phi_3= \varepsilon^{3}\chi_{0}(\varepsilon \eta)\,\Psi_{3}(x,\theta,\eta)
\quad\mbox{with}\quad
\Psi_{3}(x,\theta,\eta)=\phi_{31} (x,\theta,\eta)+\phi_{32} (x,\theta,\eta).
\end{align}

\begin{Remark}\label{remarkboundaryapproximation}
In fact, we can proceed the above arguments step by step to find boundary correction
terms and get rid of the error terms up to order $O(\ve^m)$ for any positive integer $m$.
\end{Remark}

\subsection{Summary}
We conclude that for any given parameter pair $(f_{2},e)\in\digamma$,
our {\bf final approximate solution} to the problem (\ref{problemafterscaling0})
near the surface $\Gamma_\ve$  is expressed in the local form by
\begin{align}\label{chi1}
u_4=w(x)+\varepsilon eZ+\sum\limits_{i=2}^{3}\varepsilon^{i}\varphi_{i}+\phi_1+\phi_2+\phi_3.
\end{align}

For a suitable perturbation term $\phi$,
if we locally set $u_{4}+\phi$ as the solution to problem (\ref{problemafterscaling0}),
the problem can be recast as follows
\begin{align}\label{originaldefinitionofL1}
S(u_4+\phi)=S(u_4)+L_2(\phi)+B(\phi)+N_2(\phi)=0,
\end{align}
where the linear and nonlinear operators are in the form
\begin{align}
\begin{aligned}\label{definitionofltwo1}
L_2(\phi)&=\phi_{xx}+\Delta^{\Gamma_{\var}}\phi-\phi+pu_4^{p-1}\phi,
\\
N_2(\phi)&=(u_4+\phi)^p-u_4^p-pu_4^{p-1}\phi,
\end{aligned}
\end{align}
with boundary condition
\begin{align}\label{boundaryphitwo}
\widetilde{D}_{0}(\phi)-\phi_\eta+\widetilde{D}(u_4+\phi) =g.
\end{align}
The error of the approximation is
\begin{align}\label{E1}
\begin{aligned}
E_1\,=\,&S(u_4)
\\
\,=\,&\ve^3\Delta^{\Gamma} eZ+\ve \lambda_0 eZ
+\ve^{4}\big(\,\hat{\mathbf{B}}_4+\check{\mathbf{B}}_4
+\hat{G}_{4}+\check{G}_{4}+\hat{\mathfrak{D}}_{4}
+\check{\mathfrak{D}}_{4}\,\big)
\\
&+\ve^{4}p(p-1)w^{p-2}eZ\phi_{3}
+\sum\limits_{i=4}^{6}\varepsilon^{i}S_{_{i}}+\sum\limits_{i=4}^{6}\varepsilon^{i}T_{_{i}}
+\sum\limits_{i=5}^{10}\varepsilon^{i}G_{i}+\varepsilon^{5}\mathbf{B}_{5}.
\end{aligned}
\end{align}
In the above, we also have denoted
\begin{align}\label{definitionofDzero}
\begin{aligned}
\widetilde{D}_{0}(\phi)=&\varepsilon I(\var\theta)x\phi_{x}+\varepsilon^{2}\frac{M(\var\theta)}{l^{2}_{1}(\var\theta)}x\phi_{\theta}
+\varepsilon^{3}b_{5}(\var\theta)x^{3}\phi_{x}
\\
&+\varepsilon^{3}\frac{M(\var\theta)}{l^{2}_{1}(\var\theta)}f_{1}\phi_{\theta}
+\varepsilon^{3}b_{4}(\var\theta)x^{2}\phi_{\eta},
\end{aligned}
\end{align}
\begin{align}
\begin{aligned}\label{definitionofg}
g(x,\theta)=&\varepsilon^{4}I(\var\theta)x\varphi_{3,x}+\varepsilon^{4}I(\var\theta)x\phi_{32,x}
+\sum\limits_{i=1}^{2}\varepsilon^{i+3}\Big(f_{i}I(\var\theta)+\nabla_{\rho}f_{i}\Big)\varphi_{i+2,x}
\\
&+\varepsilon^{2}b_{4}(\var\theta)\Big(x+\sum\limits_{i=1}^{2}\varepsilon^{i}f_{i}\Big)^{2}\Big(\varepsilon^{2}\nabla_{\rho}eZ
+\sum\limits_{i=2}^{4}\varepsilon^{i}\varphi_{i,\eta}\Big)-\varepsilon^{4}\varphi_{3,\eta}
\\
&-\varepsilon^{2}b_{4}(\var\theta)\Big(x+\sum\limits_{i=1}^{2}\varepsilon^{i}f_{i}\Big)^{2}
\Big(\sum\limits_{i=1}^{2}\varepsilon^{i+1}\nabla_{\rho}f_{i}\Big)\Big(\varepsilon eZ_{x}
+\sum\limits_{i=2}^{4}\varepsilon^{i}\varphi_{i,x}\Big)
\\
&-\varepsilon^{2}\frac{M(\var\theta)}{l^{2}_{1}(\var\theta)}\Big(x+\sum\limits_{i=1}^{2}\varepsilon^{i}f_{i}\Big)
\Big(\sum\limits_{i=1}^{2}\varepsilon^{i}\nabla_{\vartheta}f_{i}\Big)
\Big(\sum\limits_{i=2}^{4}\varepsilon^{i}\varphi_{i,x}+\sum\limits_{i=2}^{3}\varepsilon^{i}\phi_{i2,x}\Big)
\\
&+\varepsilon^{2}\frac{M(\var\theta)}{l^{2}_{1}(\var\theta)}\Big(x+\sum\limits_{i=1}^{2}\varepsilon^{i}f_{i}\Big)
\Big(\sum\limits_{i=2}^{4}\varepsilon^{i}\varphi_{i,\theta}+\sum\limits_{i=2}^{3}\varepsilon^{i}\phi_{i2,\theta}\Big)
\\
&+\varepsilon^{3}b_{6}(x,\var\theta)\Big(\varepsilon \nabla_{\vartheta}eZ
+\sum\limits_{i=2}^{4}\varepsilon^{i}\varphi_{i,\theta}+\sum\limits_{i=1}^{3}\varepsilon^{i}\phi_{i2,\theta}\Big)
\\
&+\varepsilon^{3}b_{6}(x,\var\theta)\Big(-\sum\limits_{i=1}^{2}
\varepsilon^{i}\nabla_{\vartheta}f_{i}\Big)\Big(\varepsilon eZ_{x}+\sum\limits_{i=2}^{4}\varepsilon^{i}\varphi_{i,x}
+\sum\limits_{i=1}^{3}\varepsilon^{i}\phi_{i2,x}\Big).
\end{aligned}
\end{align}

For later use, we decompose the error as two components
\begin{align}\label{E11}
E_1=E_{11}+E_{12},
\end{align}
where we have denoted
$$
E_{11}=\ve^3 \Delta^{\Gamma}eZ+\ve \lambda_0 e Z\quad   \mbox{ and }\ \ E_{12}=E_1-E_{11}.
$$
%
%
%
%

\subsection{Size of the errors in weighted Sobolev norms}
To estimate the size of error, we have to introduce some suitable weighted Sobolev norms.
Here we use the same norms as those in \cite{dkwy}:
for a function $h(x,z)$ defined on a set ${\mathbb E}\in \mathbb{R}^{3}$,
and for $0<\varrho<\frac{1}{100}$ and $4<q\leq+\infty$, we set
\begin{align}\label{norm}
\begin{aligned}
\|h\|_{q,\varrho;\,{\mathbb E}}=\sup_{(x,z)\in {\mathbb E}}e^{\varrho|x|}\|h\|_{L^{q}(B((x,z),1))},\quad
\\
\|h\|_{2,q,\varrho;\,{\mathbb E}}=\sum^{2}_{j=0}\sup_{(x,z)\in {\mathbb E}}e^{\varrho|x|}\|D^{j}h\|_{L^{q}(B((x,z),1))}.
\end{aligned}
\end{align}
Here $B((x,z),1)$ denotes the ball of radius 1 centered at $(x,z)$.

For the application of contraction mapping theorem in the procedure of finding the perturbation term $\phi$,
we need analyze the properties of $g$, $E_{11}$ and $E_{12}$.
The reader can refer to Section \ref{section7}.
From the uniform bound of $e$ in (\ref{conditionfore}), it is easy to see  that
\begin{align}\label{boundnessforEoneone}
||E_{11}||_{q,\varrho;\mathfrak{S}}\leq C \varepsilon^{\frac{1}{2}+1-\frac{3}{q}}.
\end{align}
where $\mathfrak{S}$ is defined in (\ref{mathfraks0}).

 All terms in $E_{12}$ carry $\ve^4$ in front, we then claim that
\begin{align}\label{boundnessforEonetwo}
||E_{12}||_{q,\varrho;\mathfrak{S}}\leq C \varepsilon^{\frac{1}{2}+4-\frac{3}{q}}.
\end{align}
A rather delicate term in $E_{12}$ is
the one carrying $\bigtriangleup^{\Gamma}f_{2}$ since we only assume a uniform bound on
$||\bigtriangleup^{\Gamma}f_{2}||_{L^{q}(\Gamma)}$. For example, we have a term
$K_1=\ve^4\bigtriangleup^{\Gamma}f_{2}w_x$ in $S(w)$ which has bound like
\begin{align*}
 ||K_1||_{q,\varrho;\mathfrak{S}}\leq C \varepsilon^{\frac{1}{2}+4-\frac{3}{q}}.
\end{align*}
Similarly, we have the following estimates
\begin{align}\label{boundnessforEonetwo1}
 || g||_{q,\varrho;\partial\mathfrak{S}} \leq C \varepsilon^{\frac{1}{2}+4-\frac{3}{q}}.
\end{align}
Other terms can be estimated in the similar way. Moreover, for the
Lipschitz dependence of the term of error $E_{12}$ on the parameters
$f_{2}$ and $e$ for the norms defined in (\ref{conditionforf}) and
(\ref{conditionfore}), we have the validity of the estimate
\begin{align}\label{Lipshitzoferroonfande}
&||E_{12}(f_{2},e)-E_{12}({\tilde f}_{2},{\tilde e})||_{q,\varrho;\mathfrak{S}}
\,\leq\, C\varepsilon^{\frac{1}{2}+4-\frac{3}{q}}\bigl(\,||f_{2}-{\tilde f}_{2}||_a+||e-{\tilde e}||_b\,\bigr).
\end{align}
Similarly we obtain
\begin{align}
||g(f_{2},e)-g({\tilde f}_{2},{\tilde e})||_{q,\varrho;\partial\mathfrak{S}}
\,\leq\, C \varepsilon^{\frac{1}{2}+4-\frac{3}{q}}\bigl(\,||f_{2}-{\tilde f}_{2}||_a+||e-{\tilde e}||_b\,\bigr).
\end{align}

\section{The Gluing Procedure}\label{section5}
\setcounter{equation}{0}

Recall that, in Sections \ref{section3} and \ref{section4},
we consider problem (\ref{problemafterscaling0}) in a small neighborhood of $\Gamma_\ve$
and find a local approximate solution.
In this section, to get a real solution to (\ref{problemafterscaling0}) by the perturbation method,
we use a gluing technique (as in  \cite{delPinoKowalczykWei1})
to reduce the problem in $\Omega_\ve$
to a projected problem on the infinite strip $ \mathfrak{S}$(cf. (\ref{mathfraks0})) in $\mathbb{R}^3$.

 Let $\sigma<r_0/100$ be a fixed number, where $r_0$ is a constant defined in
(\ref{fermicoordinate}). We consider a smooth
cut-off function $\eta_\sigma(t)$ where $t\in \mathbb{R}_{+}$ such that $
\eta_\sigma(t)= 1$ for $ 0 \leq t \leq \sigma $ and $ \eta(t)=0$ for $ t>2 \sigma$.
Set  $\eta^{\ve}_{\sigma}(s)=\eta_\sigma(\ve |s|)$, where $s$ is the normal coordinate to $\Gamma_\ve$.
Let $u_4 (s,z)$ denote the approximate solution constructed near the
surface $\Gamma_\ve$ in the coordinates $(s,z)$. We define our first global approximation to be simply
\begin{align}\label{defW}
W=\eta^{\ve}_{3\sigma}(s)u_4.
\end{align}
Obviously, $W$ is a function defined on $\Omega_\ve$, which is extended globally as $0$ outside
the $6\sigma/\ve$-neighborhood of $\Gamma_\varepsilon.$

  For  $u=W+\hat{\phi}$ where $\hat{\phi}$ globally defined in $\Omega_\ve$,
denote $$S(u)=\Delta_{\tilde{y}}u-u+u^p \quad \mbox{in} \ \Omega_\ve.$$
Then $u$ satisfies (\ref{originalproblemone}) if
and only if
\begin{align}\label{problemthreepointthree}
\widetilde{\mathcal{L}}(\hat{\phi})=-\widetilde{E}-\widetilde{N}(\hat{\phi})\quad \mbox{in} \ \Omega_\ve,
\end{align}
with boundary condition
\begin{align}\label{problemthreepointfour}
\frac{\partial\hat{\phi}}{\partial {\mathbf n}_\ve}+\frac{\partial W}{\partial {\mathbf n}_\ve}=0 \quad \mbox{on} \ \partial\Omega_\ve,
\end{align}
where
\begin{align*}
\widetilde{E}=S(W),
\quad
\widetilde{\mathcal{L}}(\hat{\phi})=\Delta_{\tilde{y}} \hat{\phi}-\hat{\phi}+pW^{p-1}\hat{\phi},
\\
\widetilde{N}(\hat{\phi})=(W+\hat{\phi})^p-W^p-pW^{p-1}\hat{\phi}.\quad
\end{align*}

  We will look for $\hat{\phi}$ in the following form
$$
\hat{\phi}=\etathree\phi+\psi,
$$
where, in the coordinates $(x,z)$ of the form (\ref{coordinatesxz}), we
assume that $\phi$ is defined in the whole strip $\mathfrak{S}$.
Obviously, (\ref{problemthreepointthree})-(\ref{problemthreepointfour}) is
equivalent to the following problem
\begin{align}
\etathree\Bigl(\Delta_{\tilde{y}}\phi-\phi+pW^{p-1}\phi\Bigr)=\etaone\Bigl{[}-\up{N}(\etathree\phi+\psi)
   -\up{E}- pW^{p-1}\psi\Bigr{]},
\end{align}
\begin{align}
\begin{aligned}\label{problemforpsi}
\Delta_{\tilde{y}}\psi-\psi+(1-\etaone)pW^{p-1}\psi&=-\ve^2 (\Delta_{\tilde{y}}\etathree)\phi
   -2\ve(\nabla_{\tilde{y}}\etathree)(\nabla_{\tilde{y}} \phi)\\
           &\quad-(1-\etaone)\up{N}(\etathree\phi+\psi)-(1-\etaone)\up{E}.
\end{aligned}
\end{align}
  On the boundary, we get
\begin{align}
\etathree \frac{\partial \phi }{\partial {\mathbf n}_\ve}+\etaone\frac{\partial W }{\partial {\mathbf n}_\ve} =0,
\end{align}
\begin{align}\label{boundaryforpsi}
\frac{\partial \psi }{\partial {\mathbf n}_\ve}+(1-\etaone)\frac{\partial W }{\partial {\mathbf n}_\ve}
      +\ve \frac{\partial \etathree }{\partial {\mathbf n}_\ve}\phi=0.
\end{align}
The key  observation is that, after solving (\ref{problemforpsi}) and (\ref{boundaryforpsi}),
the problem can be transformed to
the following nonlinear problem involving the parameter $\psi$
\begin{align}
\up{\mathcal{L}}(\phi)=\etaone\Bigl{[}\ -\up{N}(\etathree\phi+\psi) -\up{E}- p W^{p-1}\psi\ \Bigr{]} \strip
\end{align}
\begin{align} \frac{\partial \phi }{\partial {\mathbf n}_\ve}+\etaone\frac{\partial W }{\partial {\mathbf n}_\ve}=0\quad
    \mbox{on}\ \partial\mathfrak{S}.
\end{align}
Notice that the operator $\up{\mathcal{L}}$ in $\Omega_{\ve}$ may be taken as any
compatible extension outside the $6\sigma/\ve$-neighborhood of
$\Gamma/\ve $ in the strip $\mathfrak{S}$ and the operator $\partial/\partial \mathbf{n}_{\var}$
may be taken as any compatible extension outside the $6\sigma/\ve$-neighborhood of
$\Gamma/\ve$ on the boundary $\partial\mathfrak{S}$.

  First, we solve, given a small $\phi$, problem
(\ref{problemforpsi}) and (\ref{boundaryforpsi}) for $\psi$. Assume
now that $\phi$ satisfies the following decay property
\begin{align}\label{decayofphi}
\bigl{|}\nabla\phi(y)\bigr{|}+\bigl{|}\phi(y)\bigr{|}
    \leq e^{-\gamma/\ve}\ \ \mbox{if}\ \  |s|>\sigma/\ve,
\end{align}
for certain constant $\gamma>0$.
The solvability can be done in the following way: let us observe that
$W$ is exponentially small for $|s|>\sigma/\ve$, where $s$ is the
normal coordinate to $\Gamma/\ve$.  Then the problem
\begin{align}\nonumber
\triangle\psi-\bigl[1-(1-\etaone)pW^{p-1}\bigr]\psi=h \quad \mbox{in}\ \Omega_\ve,
\end{align}
\begin{align}\nonumber
\frac{\partial \psi }{\partial {\mathbf n}_\ve}=-(1-\etaone)\frac{\partial W }{\partial {\mathbf n}_\ve}-\ve \frac{\partial
\etathree }{\partial {\mathbf n}_\ve}\phi\quad \mbox{on}\ \partial\Omega_\ve,
\end{align}
has a unique bounded solution $\psi$ whenever $||h||_{\infty}\leq +\infty$.
Moreover, $$||\psi||_{\infty}\leq C||h||_{\infty}.$$
  Since $\up{N}$ is power-like with power greater than one, a direct
application of contraction mapping principle yields that
(\ref{problemforpsi}) and (\ref{boundaryforpsi}) has a unique
(small) solution $\psi=\psi(\phi)$ with
\begin{align}\label{boundforpsi}
||\psi(\phi)||_{L^{\infty}}
\leq C\ve \bigl[\ ||\phi||_{L^{\infty}(|s|>\sigma/\ve)}
+||\nabla\phi||_{L^{\infty}(|s|>\sigma/\ve)}+e^{-\sigma/\ve}\ \bigr],
\end{align}
where ${|s|>\sigma/\ve}$ denotes the complement in $\Omega_{\ve}$ of
$\sigma/\ve$-neighborhood of $\Gamma/\ve$.   Moreover, the
nonlinear operator $\psi$ satisfies a Lipschitz condition of the form
\begin{align}\label{lipshitzforpsi}
||\psi(\phi_1)-\psi(\phi_2)||_{L^{\infty}}
\leq C\ve \bigl[\ ||\phi_1-\phi_2||_{L^{\infty}(|s|>\sigma/\ve)}
+||\nabla\phi_1-\nabla\phi_2||_{L^{\infty}(|s|>\sigma/\ve)}\ \bigr].\ \
\end{align}

  Therefore, from the above discussion, the full problem has been reduced to solving the
following (nonlocal) problem in the infinite strip $\mathfrak{S}$
\begin{align}\label{problemofphione}
\mathcal{L}_2(\phi)=\etaone\Bigl{[}-\up{N}(\etathree\phi+\psi(\phi)) -\up{E}- p W^{p-1}\psi(\phi)\Bigr{]} \strip
\end{align}
\begin{align}\label{problemofphitwo}
\mathfrak{B}(\phi)+\etaone\frac{\partial W }{\partial {\mathbf n}_\ve}=0\quad \quad
    \mbox{on}\ \partial\mathfrak{S},
\end{align}
for $\phi\in W^{2,q}(\mathfrak{S})$ satisfying condition (\ref{decayofphi}). Here $\mathcal{L}_2$
denotes a linear operator that coincides with $\up{\mathcal{L}}$ on the region $|s|<8\sigma/\ve$,
$\mathfrak{B}$ denotes the inward normal derivatives of $\mathfrak{S}$ that coincides with inward normal
$\partial/\partial {\mathbf n}_\ve $ of $\Omega_\ve$ on the region $|s|<8\sigma/\ve$.

The definitions of these operators can be shown as follows.
The operator $\up{\mathcal{L}}$ for $|s|<8\sigma/\ve$ is given in coordinates $(x,z)$ by formula
(\ref{coordinatesxz}). We extend it for functions $\phi$ defined in the strip $\mathfrak{S}$
in terms of $(x,z)$ as the following
\begin{align}\label{definitionofLtwo}
\mathcal{L}_2(\phi)=L_2(\phi)+\chi(\ve |x|)B(\phi)\qquad\,\mbox{in}\,\mathfrak{S}
\end{align}
where $\chi(r)$ is a smooth cut-off function which equals 1 for
$0\leq r<10\sigma$ and vanishes identically for $r>20\sigma$, $L_2$ and $B$
are the operators defined in (\ref{definitionofltwo1}) and (\ref{problemafterscaling2}).
Similarly, the boundary conditions can be written as
\begin{align}
\chi(\ve |x|) \widetilde{D}_{0}(\phi)-\frac{\partial\phi}{\partial\tau_{\var}}
+\chi(\ve |x|) \widetilde{D}(W+\phi) =\chi(\ve |x|)\ g\quad\mbox{on}\ \partial\mathfrak{S}
\end{align}
where the operators $\widetilde{D}_{0}$ and $\widetilde{D}$ are
defined in (\ref{boundaryphitwo}) and (\ref{boundaryconditionv}), $\tau_{\var}$ is the unit inward normal of $\partial\mathfrak{S}$.

Rather than solving problem (\ref{problemofphione})-(\ref{problemofphitwo})
satisfying the boundary condition, we deal with the following projected problem:
for each pair of parameters $f_{2}$ and $e$ in $\digamma$, finding functions
$\phi\in W^{2,q}(\mathfrak{S}),\ c, d\in L^{q}(\Gamma) $ and $\Lambda_{1},\Lambda_{2}$ such that
\begin{align}\label{problemforphitwoone}
\mathcal{L}_2(\phi)=-\chi E_1-\chi \widetilde{N}_3(\phi)
+c(\ve z)\ \chi w_x+d(\ve z)\ \chi Z\quad \,\mbox{in}\,\mathfrak{S},
\end{align}
\begin{align}\label{problemforphitwotwo}
\chi \widetilde{D}_0(\phi) -\frac{\partial\phi}{\partial\tau_{\var}} +\chi \widetilde{D}(W+\phi)=\chi g \quad\mbox{on}\ \partial\mathfrak{S},
\end{align}
\begin{align}\label{problemforphitwothree}
\intreal \phi(x,z)w_x\mathrm{d}x=\Lambda_{1},\,\,
\intreal \phi(x,z)Z(x)\mathrm{d}x=\Lambda_{2}\,\quad\,\mbox{in}\,\,\Gamma_{\var},
\end{align}
where $\widetilde{N}_3(\phi)= \up{N}(\etathree \phi+\psi(\phi)) +p W^{p-1}\psi(\phi)$.
In Sections \ref{section6} and \ref{section7}, we will prove that this problem has a unique solution $\phi$ whose norm is
 controlled by the $\|\cdot\|_{q,\varrho}$-norm, not of the error component $E_{11}$,
 but rather of the components $E_{12}$ and $g$. Moreover, $\phi$ will satisfies (\ref{decayofphi}).
The reader can refer to the conclusion in Proposition \ref{proposition7.1}.

  After this has been done, our task is to adjust the parameters $f_{2}\ \mbox{and}\ e$ such that
  the functions $c \ \mbox{and}\ d$ are identically zero.
  It is equivalent to solving a nonlocal, nonlinear coupled second order system of differential
equations for the pair $(f_{2},e)$ with suitable boundary conditions.
Indeed,  we will derive the system of differential equations for the unknown functions $f_2,\, e$
in Section \ref{section8},
and then show the solvability of this system on the infinite dimensional
space $\digamma$(cf. (\ref{regionoffande})) in Section \ref{section9}.

\section{The Invertibility of $\mathcal{L}_2$}\label{section6}
\setcounter{equation}{0}

  Let $\mathcal{L}_2$ be the operator defined by
 (\ref{definitionofLtwo}) and $g$ be the functions in (\ref{definitionofg}). Note that the function
$\chi(\ve|x|)$ is even in the definition of $\mathcal{L}_2$. In this section,
 We study the following linear problem: for given $h\in L^{q}(\mathfrak{S})$ and $g\in L^{q}(\partial\mathfrak{S}) $,
 finding functions $\phi\in W^{2,q}(\mathfrak{S}),\ c, d \in L^{q}(\partial\mathfrak{S})$ and $\Lambda_{1},\Lambda_{2}$ such that
\begin{align}\label{4point1}
\begin{aligned}
\mathcal{L}_2(\phi)&=\chi h+c(\ve z)\, \chi w_x+d(\ve z)\, \chi Z \quad\,\,\mbox{in}\,\mathfrak{S},
\\
\chi \widetilde{D}_0(\phi)&-\frac{\partial\phi}{\partial\tau_{\var}}+\chi \widetilde{D}(W+\phi)
=\chi g\quad\,\,\mbox{on}\,\partial\mathfrak{S},
\\
\intreal\phi w_x(x)&\mathrm{d}x=\Lambda_{1},\ \ \  \
\intreal\phi Z(x)\mathrm{d}x=\Lambda_{2}\quad\,\mbox{in}\,\,\Gamma_{\var}.
\end{aligned}
\end{align}

\begin{proposition}\label{operatorTtwo}
If $\ \sigma$ in the definition of $\mathcal{L}_2$ is chosen small enough and $h\in L^{q}(\mathfrak{S})$,
then there exists a constant $C>0$, independent of $\ve$,
such that for all small $\ve$, the problem (\ref{4point1}) has a unique
solution $\phi=T_1(h, g)$ with suitable $\Lambda_1$ and $\Lambda_2$ which satisfy
\begin{align}
\begin{aligned}
||\phi||_{q,\varrho}&\leq C(||h||_{L^{q}(\mathfrak{S})} + || g ||_{L^{q}(\partial\mathfrak{S})}),
\\
||\Lambda_{i}||_{L^{q}(\Gamma_{\var})}&\leq C(||h||_{L^{q}(\mathfrak{S})} + || g ||_{L^{q}(\partial\mathfrak{S})}),\quad \forall\,\, i=1,2.
\end{aligned}
\end{align}
Moreover, if $h, g$ have compact supports contained in $|x|\leq 20\sigma/\ve$, then
\begin{align} \label{propertyofTtwo0}
\bigl |\phi(x,z)\bigr|+\bigl|\nabla
\phi(x,z)\bigr|\leq ||\phi||_{L^{\infty}}\ e^{-2\sigma/\ve}\quad \mbox{for}\ |x|>40\sigma/\ve.
\end{align}
\end{proposition}

\proof Note that the problem can be written as
\begin{align*}
\phi_{xx}+\Delta^{\Gamma_{\var}}\phi-\phi+pw^{p-1}\phi=&-p(W^{p-1}-w^{p-1})\phi-\chi B(\phi)
\\
&+\chi h+c(\ve z)\ \chi w_x+d(\ve z)\ \chi Z\quad\,\,\mbox{in}\,\mathfrak{S},
\\
\frac{\partial\phi}{\partial\tau_{\var}}=-\chi g+\chi \widetilde{D}_0(\phi)&+\chi \widetilde{D}(W+\phi)\quad\,\,\mbox{on}\,\partial\mathfrak{S},
\\
\intreal\phi w_x\mathrm{d}x=\Lambda_{1},\ \ \ \ \
\intreal\phi &Z(x)\mathrm{d}x=\Lambda_{2}\quad\,\mbox{in}\,\,\Gamma_{\var}.
\end{align*}

  Let
\begin{align*}
\varphi=T_1\Bigl{(}\chi h-p(W^{p-1}-w^{p-1})\phi-\chi
B(\phi),\,G \Bigr{)},
\end{align*}
where
\begin{align*}
G(\phi)=-\chi g+\chi \widetilde{D}_0(\phi)+\chi \widetilde{D}(W+\phi),
\end{align*}
and $T_1$ is the bounded operator defined by Lemma \ref{lemma0.12}.

The key point is that the operator
$$
\widetilde{B}_4(\phi)=-\chi B(\phi)-p(W^{p-1}-w^{p-1})\phi,
$$
is small in the sense that
$$||\widetilde{B}_4(\phi)||_{q,\varrho}\leq C \sigma ||\phi||_{q,\varrho}.$$
Similar results hold for $G(\phi)$. Hence, the results can be derived by
the invertibility conclusion of Lemma \ref{lemma0.12} if we choose $\sigma$
sufficiently small.

  Since $\chi$ is supported on $|x|<20\sigma/\ve$, then $\phi$ satisfies for $|x|>40\sigma/\ve$
a problem of the form
\begin{align*}\nonumber
\phi_{xx}+\Delta^{\Gamma_{\var}}\phi-(1+o(1))\phi\,=\,0\quad\mbox{for}\quad |x|>40\sigma/\ve,
\qquad
\frac{\partial\phi}{\partial\tau_{\var}}\,=\,0.
\end{align*}
Hence, the validity of formula (\ref{propertyofTtwo0}) can be showed easily.
\qed

\section{Solving the Nonlinear Projection Problem }\label{section7}

In this section, we will solve (\ref{problemforphitwoone})-(\ref{problemforphitwothree})
in $\mathfrak{S}$. A first elementary, but crucial observation is the following:  The term
$$
E_{11}=\ve^3\Delta^{\Gamma}eZ+\ve\lambda_0e Z,
$$
in the decomposition of
$E_1$, has precisely the form $d(\ve z)Z$ and can be absorbed in that
term $ \chi d(\var z) Z$. Then, the equivalent equation of (\ref{problemforphitwoone}) is
$$\mathcal{L}_2(\phi)=\chi E_{12}+\chi \widetilde{N}_3(\phi)+c(\ve z)\ \chi w_x+d(\ve z)\ \chi Z.$$

Let $T_1$ be the bounded operator defined by Proposition \ref{operatorTtwo}. Then the problem
(\ref{problemforphitwoone})-(\ref{problemforphitwothree}) is equivalent to the following fixed point problem
\begin{align}\label{fixedpointproblem}
\phi
\,=\,
T_1\bigl(\, \chi E_{12}+\chi \widetilde{N}_3(\phi),\ \chi g\,\bigr)
\,\equiv\,
\mathcal{A}(\phi).
\end{align}
 We collect some useful facts to find the domain of the operator $\mathcal{A}$ such that
$\mathcal{A}$ becomes a contraction mapping.

  The big difference between $E_{11}$ and $E_{12}$ is
their sizes. From (\ref{boundnessforEoneone}) and (\ref{boundnessforEonetwo})
\begin{align}\label{boundnessofEtwo}
||E_{12}||_{q,\varrho}\leq c_{*}\ \ve^{\frac{1}{2}+4-\frac{3}{q}},
\end{align}
while $E_{11}$ is only of size $O(\ve^{1+1/2-3/q})$.  Similarly, we have
\begin{align}\label{boundnessofEtwon}
||g||_{q,\varrho;\partial\mathfrak{S}} \leq c_{*}\ \ve^{\frac{1}{2}+4-\frac{3}{q}}.
\end{align}

  The operator $T_1$ has a useful property: assume $\hat{h}$ has a support contained in $|x|\leq 20\sigma/\ve$,
then $\phi=T_1(\hat{h})$ satisfies the estimate
\begin{align}\label{propertyofTtwo}
\bigl |\phi\bigr|+\bigl|\nabla
\phi\bigr|\leq ||\phi||_{\infty,\varrho}\ e^{-2\sigma/\ve}\quad \mbox{for}\ |x|>40\sigma/\ve.
\end{align}
  Recall that the operator $\psi(\phi)$ satisfies, as seen directly from its definition
\begin{align*}
||\psi(\phi)||_{\infty,\varrho}
\leq C\ve \Bigl[\ \bigl|\bigl| \ |\phi| +|\nabla\phi|\
\bigl|\bigl|_{L^{\infty}(|x|>20\sigma/\ve)}\,+\,e^{-\sigma/\ve}\ \Bigr],
\end{align*}
and a Lipschitz condition of the form
\begin{align*}
||\psi(\phi)-\psi(\tilde{\phi})||_{\infty,\varrho}
\leq C\ve \Bigl[\ \bigl|\bigl| \ |\phi-\tilde{\phi}| +|\nabla(\phi-\tilde{\phi})|\
\bigl|\bigl|_{L^{\infty}(|x|>20\sigma/\ve)}\ \Bigr].
\end{align*}
  Now, the facts above will allow us to construct a region
where contraction mapping principle applies and then solve the
problem (\ref{problemforphitwoone})-(\ref{problemforphitwothree}). Consider the following closed, bounded subset
\begin{align}
\mathfrak {D} =
\left\{
\phi\in H^2(\mathfrak{S})
\left|
\begin{array}{l}
\ \|\phi\|_{2,q,\varrho} \le \varsigma\ve^{\frac{1}{2}+4-\frac{3}{q}},\medskip
\\
\ \Bigl|\Bigl|\ |\phi| +|\nabla\phi|\  \Bigl|\Bigl|_{L^{\infty}(|x|>40\sigma/\ve)}\ \leq\ ||\phi||_{2,q,\varrho}\ e^{-\sigma/\ve}.
\end{array}
\right.
\right\}
\end{align}

  We claim that if the constant $\varsigma$ is sufficiently large, then the
map $\mathcal{A}$ defined in (\ref{fixedpointproblem}) is a contraction
form $\mathfrak{D}$ into itself. Let us analyze the Lipschitz character of the nonlinear
operator involved in $\mathcal{A}$ for functions in $\mathfrak{D}$
\begin{align*}
\chi \widetilde{N}_3(\phi)&=\chi \widetilde{N}_1(\phi+\psi(\phi))+\chi p W^{p-1}\psi(\phi)
\\
&\equiv\bar{N_3}(\phi)+\chi p W^{p-1}\psi(\phi).
\end{align*}
Note that
$\widetilde{N}_1(\varphi)=p\bigl[(W+t\varphi)^{p-1}-W^{p-1}\bigr]\varphi\ $ for some $t\in(0,1)$.
  From here it follows that
$$|\widetilde{N}_1(\varphi)|\leq C(|\varphi|^p+|\varphi|^2).$$
  Denoting $S_\sigma=\mathfrak{S}\cap\bigl\{|x|<10\sigma/\ve \bigr\}$,
we have that for $\phi\in \mathfrak{D}$
$$ ||\bar{N}_3(\phi)||_{q,\varrho}
\leq C\bigl[\
||\phi||^p_{qp,\varrho}+||\phi||^2_{2q,\varrho}+
||\psi(\phi)||^p_{qp,\varrho;S_\sigma} +||\psi(\phi)||^2_{2q,\varrho;S_\sigma}\ \bigr].$$
Using Sobolev's embedding, we derive
$$||\phi||^p_{qp,\varrho}+||\phi||^2_{2q,\varrho}
\leq C\bigl(\ ||\phi||^p_{2,q,\varrho}+||\phi||^2_{2,q,\varrho}\ \bigr).$$
Using estimates (\ref{boundforpsi}), the facts that $\phi\in \mathfrak{D}$, (\ref{propertyofTtwo}),
that of the area of $S_\sigma$ is of order $O(\sigma/\ve)$ and
Sobolev's embedding, we get
$$||\psi(\phi)||^p_{qp,\varrho;S_\sigma}
+||\psi(\phi)||^2_{2q,\varrho;S_\sigma}
\leq C e^{-\sigma/{4\ve}}\Bigl[\ 1+||\phi||^p_{2,q,\varrho}
+||\phi||^2_{2,q,\varrho}\ \Bigr]. $$
Hence, from the properties of
$W$ and $\psi(\phi)$ we obtain
\begin{align}\label{boundsofNtwo}
||\chi \widetilde{N}_3(\phi)||_{qp,\varrho}\leq C (\ve^{(\frac{1}{2}+4-\frac{3}{q})p}\varsigma^p+\ve^{9-\frac{6}{q}}\varsigma^2 ).
\end{align}

  As for Lipschitz condition, after a direct calculation we find
\begin{multline*}
||\widetilde{N}_1(\varphi_1)-\widetilde{N}_1(\varphi_2)||_{q,\varrho}
\leq C\Bigl[ ||\varphi_1||^{p-1}_{qp,\varrho}+||\varphi_1||_{2q,\varrho}
+||\varphi_2||^{p-1}_{qp,\varrho}+||\varphi_2||_{2q,\varrho}\Bigr]\\
\times\bigl(||\varphi_1-\varphi_2||_{qp,\varrho}+||\varphi_1-\varphi_2||_{2q,\varrho}\bigr).\quad
\end{multline*}
  Hence,
\ba
\begin{split}
||\bar{N}_3(\phi)-\bar{N}_3(\tilde{\phi})||_{q,\varrho}&\leq
||N_1(\phi+\psi(\phi))-N_1(\tilde{\phi}+\psi(\phi))||_{q,\varrho;S_\sigma}\\
&\quad+||N_1(\tilde{\phi}+\psi(\phi))-N_1(\tilde{\phi}+\psi(\tilde{\phi}))||_{q,\varrho;S_\sigma}\\
&\leq \upsilon\Bigl(||\phi-\tilde{\phi}||_{2q,\varrho;S_\sigma}+||\phi-\tilde{\phi}||_{qp,\varrho;S_\sigma}\Bigr)\\
&\quad+ \upsilon\Bigl(||\psi(\phi)-\psi(\tilde{\phi})||_{2q,\varrho;S_\sigma}
       +||\psi(\phi)-\psi(\tilde{\phi})||_{qp,\varrho;S_\sigma}\Bigr),
\end{split}
\ea
where $\upsilon=\upsilon_1+\upsilon_2$ with
$$\upsilon_l=||\phi_l||^{p-1}_{qp,\varrho;S_\sigma}+||\psi(\phi_l)||^{p-1}_{qp,\varrho;S_\sigma}
+||\phi_l||_{2q,\varrho;S_\sigma}+||\psi(\phi_l)||_{2q,\varrho;S_\sigma}.
$$
  Arguing as above and using the Lipschitz dependence of $\psi$ on $\phi$,
it can be derived
\begin{align}\label{contractionpropertyofNtwo}
||\chi \widetilde{N}_3(\phi)-\chi \widetilde{N}_3(\tilde{\phi})||_{q,\varrho}\leq
    C\bigl(\ve^{(\frac{1}{2}+4-\frac{3}{q})(p-1)}\varsigma^{p-1}
    +\ve^{\frac{1}{2}+4-\frac{3}{q}}\varsigma\bigr)||\phi-\tilde{\phi}||_{2,q,\varrho}.
\end{align}

  Now, we can find the solution of (\ref{fixedpointproblem})
in the sequel. Let $\phi\in \mathfrak{D}$ and
$\nu=\mathcal{A}(\phi)$, then from (\ref{boundnessofEtwo})-(\ref{boundnessofEtwon}) and (\ref{boundsofNtwo})
$$
||\nu||_{2,q,\varrho}
  \leq ||T_1||\Bigl[\ c_{*}\ve^{\frac{1}{2}+4-\frac{3}{q}}+C\varsigma^p\ve^{(\frac{1}{2}+4-\frac{3}{q})p}
  +C\varsigma^2\ve^{9-\frac{6}{q}}\ \Bigr].
$$
Choosing any number $\varsigma>C_{*}||T_1||,$ we get that for small $\ve$
$$
||\nu||_{2,q,\varrho}\leq \varsigma \ve^{\frac{1}{2}+4-\frac{3}{q}}.
$$
From (\ref{propertyofTtwo})
$$
\Bigl|\Bigl||\nu|+|\nabla \nu|\Bigr|\Bigr|_{L^{\infty}(|x|>40\sigma/\ve)}
\leq ||\nu||_{\infty}\ e^{-2\sigma/\ve}\leq
||\nu||_{2,q,\varrho}\ e^{-\sigma/\ve}.
$$
  Therefore, $\nu\in \mathfrak{D}.$ $\mathcal{A}$ is
clearly a contraction thanks to (\ref{contractionpropertyofNtwo})
and we can conclude that (\ref{fixedpointproblem}) has a unique
solution in $\mathfrak{D}$.

  The error $E_{12}$ and the operator $T_1$ itself carry
the functions $f_{2}$ and $e$ as parameters. For future reference, we should consider
their Lipschitz dependence on these parameters. (\ref{Lipshitzoferroonfande}) is just the
 formula about the Lipschitz dependence of error $E_{12}$ on these two parameters.
The other task can be realized by careful and direct computations of all terms involved in the
differential operator which will show this dependence is indeed Lipschitz with
respect to the $W^{2,q}$-norm (for all $\ve$).

  Within the operator, consider for instance the following term involving $\Delta^{\Gamma}f_{2}$
$$Q_{f_{2}}(\phi)=\ve^4\Delta^{\Gamma}f_{2}\phi_x.$$
Then we have
$$
||Q_{f_{2}}(\phi)||^q_{L^q(B)}\leq
\ve^{4q-2} \int_{\Gamma}\bigl|\Delta^{\Gamma}f_{2}\bigr|^q\Bigl( \sup_{z}\intreal |\phi_x(x,z)|^q\mathrm{d}x\Bigr).
$$
Let $\mu(z)=\intreal |\phi_x(x,z)|^q\mathrm{d}x$.
Then there holds
\bb
\begin{split}
\sup_{z} \mu(z)&\leq \ve\int_{\mathfrak{S}}|\phi_x|^q\mathrm{d}x+\frac{1}{\ve}\int_{\mathfrak{S}}|\phi_x|^{q-1}|\nabla^{\Gamma}\phi_{x}|\mathrm{d}x\\
&\leq \frac{1}{2}\sup_{z} \mu(z)+\frac{C}{\ve^2}\int_{\mathfrak{S}}|\nabla^{\Gamma}\phi_{x}|^q\mathrm{d}x,
\end{split}
\ee
and we can obtain $$\mu(z)\leq C\ve^{-2}||\phi||^q_{2,q,\varrho}.$$
  Therefore, $$||Q_{f_{2}}(\phi)||_{q,\varrho}\leq C\ve||f_{2}||_a.$$
  Similar estimates can be applied to other terms in the operator involving $\Delta^{\Gamma}f_{2}$.

For the linear operator $T_1$, we have  the following Lipschitz
dependence $$||T_1(f_2)-T_1(\tilde{f}_2)||_{2,q,\varrho}\leq C\ve ||f_2-\tilde{f}_2||_a.$$
   Moreover, the operator $\widetilde{N}_3$
also has Lipschitz dependence on $(f_{2}, e)$. It is easily checked
that for $\phi\in \mathfrak{D}$ we have, with obvious notation
$$
||\chi \widetilde{N}_{3,(f_2, e)}(\phi)-\chi \widetilde{N}_{3,(\tilde{f}_2, \tilde{e})}(\phi)||_{q,\varrho}
\leq C\ve^{\frac{1}{2}+4-\frac{3}{q}}\Bigl[\ ||f_2-\tilde{f}_2||_a+||e-\tilde{e}||_b\ \Bigr].
$$
Hence, from the fixed point characterization we get that
\begin{align}\label{Lipshitzofphionfande}
||\phi(f_2,e)-\phi(\tilde{f}_2,\tilde{e})||_{2,q,\varrho}
\leq C\ve^{\frac{1}{2}+4-\frac{3}{q}}\bigl[\ ||f_2-\tilde{f}_2||_a+||e-\tilde{e}||_b\ \bigr].
\end{align}

As a conclusion, we give the proposition.
\begin{proposition}\label{proposition7.1}
There is a number $\varsigma>0$ such that for all $\ve$ small enough and
all parameters $(f_{2}, e)$ in $\digamma$, problem (\ref{problemforphitwoone})-(\ref{problemforphitwothree})
has a unique solution $\phi=\phi(f_{2}, e)$ which satisfies
$$
 ||\phi||_{2,q,\varrho}\leq \varsigma\ve^{\frac{1}{2}+4-\frac{3}{q}},\quad ||\Lambda_{i}||_{q,\varrho;\Gamma_{\var}}\leq \varsigma\ve^{\frac{1}{2}+4-\frac{3}{q}},\quad \forall\,\, i=1,2,
$$
$$
\Bigl|\Bigl|\ |\phi| +|\nabla\phi|\ \Bigl|\Bigl|_{L^{\infty}(|x|>40\sigma/\ve)}
 \leq ||\phi||_{2,q,\varrho}e^{-\sigma/\ve}.
$$
Moreover,\, $\phi(f_{2},e)$ depends  Lipschitz-continuously on the parameters $f_{2}$ and \,$e$ in
the sense of  the estimate (\ref{Lipshitzofphionfande}).
\end{proposition}
\qed

\section{Estimates of the Projections Against $w_x$ and $Z$ }\label{section8}
\setcounter{equation}{0}

  As we mentioned in Section \ref{section5},
we will set up the system of differential equations for
the unknown parameters $f_{2}$ and $e$ defined on $\Gamma$ which are equivalent to making $\,c,\,d$
zero in the system (\ref{problemforphitwoone})-(\ref{problemforphitwothree}).
On the boundary, we have imposed the conditions (\ref{boundaryf2}) and (\ref{boundarye})
for the parameters $f_{2}$ and $e$, which  are  restated respectively in the following
\begin{align}
\nabla_{\rho}f_{2}+I(\var\theta)f_{2}+\mathfrak{R}(\theta)=0,
\quad
\nabla_{\rho}e+\frac{1}{2}I(\var\theta)e=0.
\end{align}
On the interior of $\Gamma$, these equations are obtained by simply integrating the equations
(\ref{problemforphitwoone}) (only in $x$) against $w_x$ and $Z$ respectively.
It is easy to derive the following equations
\begin{align}\label{6point1}
\begin{aligned}
\int_{\mathbb{R}}\Bigl{[}\chi E_1+\chi \widetilde{N}_3(\phi)+\Delta^{\Gamma_{\var}}\phi+\chi B(\phi)+p\bigl(W^{p-1}-w^{p-1}\bigr)\phi\Bigr{]} w_x {\mathrm{d}}x\,=\,0,
\end{aligned}
\end{align}
\begin{align}\label{6point2}
\begin{aligned}
\int_{\mathbb{R}}\Bigl{[}\chi E_1+\chi\widetilde{N}_3(\phi)
+\Delta^{\Gamma_{\var}}\phi+\lambda_{0}\phi+\chi B(\phi)
+p\bigl(W^{p-1}-w^{p-1}\bigr)\phi\Bigr{]} Z {\mathrm{d}}x\,=\,0,
\end{aligned}
\end{align}
where the error term $E_1$ is defined in (\ref{E1}), the operators $\widetilde{N}_{3}$ and $B$ are defined
in (\ref{problemforphitwoone}) and (\ref{problemafterscaling2}). It is crucial to estimate the terms
$$
\intwx{E_1} \quad \mbox{and} \quad\intz{E_1}.
$$
The same arguments can be applied to other terms  in (\ref{6point1}) and (\ref{6point2}).
Now, we divide the estimates for the components in (\ref{6point1}) and (\ref{6point2})
into several parts.

\subsection{Estimates of the projection against $w_{x}$}
First, multiplying (\ref{E1}) by $w_x$ and integrating over the variable $x$, using the decomposition of $E_1$
in (\ref{E11}) and the fact that  $w_x$ is an odd function in $x$, we obtain
\begin{align*}
\int_{\mathbb{R}}E_{1}w_{x}\mathrm{d}x\,=\,\int_{\mathbb{R}}E_{12}w_{x}\mathrm{d}x.
\end{align*}
More precisely, there holds
\begin{align*}
\int_{\mathbb{R}}E_{12}w_{x}\mathrm{d}x&=\ve^4\int_{\mathbb{R}}\big[\check{\mathbf{B}}_4
+\check{G}_{4}+\check{\mathfrak{D}}_{4}\,+p(p-1)w^{p-2}eZ\phi_{3}\big]w_{x}\mathrm{d}x
\\
&\quad
+\int_{\mathbb{R}} \Big[\sum\limits_{i=4}^{6}\varepsilon^{i}S_{i}+\sum\limits_{i=4}^{6}\varepsilon^{i}T_{i}
+\sum\limits_{i=5}^{10}\varepsilon^{i}G_{i}+\varepsilon^{5}\mathbf{B}_{5}\Big]w_{x}\mathrm{d}x.
\end{align*}
By direct calculation, we derive  that
\begin{align}\label{E1timewx}
\begin{aligned}
\int_{\mathbb{R}}E_{12}w_{x}\mathrm{d}x\,=\, &-\ve^4\sigma_1 \Big[
\Delta^{\Gamma}f_{2} +(k_{1}^{2}+k_{2}^{2})f_{2}\Big]
+ \gamma_1 \ve^{4} e + \gamma_2\ve^{6}\Delta^{\Gamma}e
\\
& +\ve^{5}b_{1\ve}\Delta^{\Gamma}e
+\ve^{5}b_{2\ve}\Delta^{\Gamma}f_{2}
+\sum\limits_{i=5}^{10}\varepsilon^{i}b_{i\ve},
\end{aligned}
\end{align}
where $\gamma_1, \gamma_2$ are two constants and
$$
\sigma_1=\int_{\R}w_{x}^{2}\mathrm{d}x.
$$
Here and below we denote by $b_{l\ve}, l=1, 2, i$, generic, uniformly bounded continuous functions of the form
\[ b_{l\ve}\,=\,b_{l\ve}(f_{2},e,\nabla^{\Gamma} f_{2},\nabla^{\Gamma} e),
\]
where additionally $b_{l\ve}$ is uniformly Lipschitz in its arguments.

\subsection{Projection of terms involving $\phi$}
We will estimate other terms that involve $\phi$ in (\ref{6point1}),
\begin{align}\label{7point6}
\int_{\mathbb{R}}\Bigl{[}\chi N_2(\phi)+\Delta^{\Gamma_{\var}}\phi+\chi B(\phi)+p\bigl(W^{p-1}-w^{p-1}\bigr)\phi\Bigr{]} w_x {\mathrm{d}}x.
\end{align}
Using the condition in (\ref{problemforphitwothree}),
we first estimate $\int_{\mathbb{R}}\Delta^{\Gamma_{\var}}\phi w_x {\mathrm{d}}x$,
the estimate for the term can be done as follows
\begin{align*}
\Upsilon_{1}(z)=\int_{\mathbb{R}}\Delta^{\Gamma_{\var}}\phi w_x {\mathrm{d}}x
=\Delta^{\Gamma_{\var}}\int_{\mathbb{R}}\phi w_x {\mathrm{d}}x
\equiv \Delta^{\Gamma_{\var}}\Lambda_{1}.
\end{align*}
With the help of the Proposition \ref{proposition7.1}, we have the following estimate
\begin{align*}
\|\Upsilon_1\|_{L^{q}(\Gamma_{\var})} \leq C \ve^{\frac{1}{2}+4-\frac{3}{q}}.
\end{align*}

 The last two components in (\ref{7point6}) are
$$
\Upsilon_2\,=\,\int_{\mathbb{R}} B(\phi) w_x\mathrm{d}x
\quad\mbox{and}\quad
\Upsilon_3\,=\,\int_{\mathbb{R}} p\big(W^{p-1} - w^{p-1}\big)\phi w_x\mathrm{d}x.
$$
Here we recalled the definitions of the operator $B$ in (\ref{problemafterscaling2})
and the local approximation $W$ in (\ref{defW}).
We make the following observation: all terms in $B(\phi)$ carry $\ve^{2}$ and involve powers
of $x$ times derivatives of $1, 2$ or two orders of $\phi$. The
conclusion is that since $w_x$ has exponential decay then
\begin{align*}
\int_{\Gamma} |\Upsilon_2(\theta)|^{q} \mathrm{d}\theta \leq
C \ve^{4}
\|\phi\|^{q}_{2,q,\varrho}.
\end{align*}
Hence there holds
\begin{align*}
\|\Upsilon_2\|_{L^{q}(\Gamma)} \leq C \ve^{\frac{1}{2}+4+\frac{1}{q}}.
\end{align*}
In $B(\phi)$ we single out two less regular terms. The one whose
coefficient depends on $\Delta^{\Gamma}f_{2}$ explicitly has the form
\begin{align*}
\Upsilon_{2*}
&\,=\, \ve^k \Delta^{\Gamma} f_{2} \int_{\mathbb{R}} \phi_x Z
\Bigg(1 +\ve\Big(x-\sum\limits_{l=0}^{2}\ve^l f_l\Big)^{-2}\Bigg)
\\
&\,=\,-\ve^k \Delta^{\Gamma} f_{2} \int_{\mathbb{R}} \phi\Bigg\{ Z\Big(1 +\ve\big(x-\sum\limits_{l=0}^{2}\ve^l f_l\big)\Big)^{-2}\Bigg\}_x.
\end{align*}
Since $\phi$ has Lipschitz dependence on $(f_{2}, e)$ in the form
(\ref{Lipshitzofphionfande}), we see that
\begin{align}\label{LipshitzforLemdaone}
\begin{aligned}
\|\Upsilon_{2*}(f_{2}, e) - \Upsilon_{2*}(\tilde{f}_{2},
\tilde{e}) \|_{L^{q}(\Gamma)}
\,\leq\,
C \ve^{\frac{1}{2}+4-\frac{3}{q}}\big(\|f_{2}-\tilde{f}_{2} \|_a + \|e - \tilde{e}\|_b \big).
\end{aligned}
\end{align}
The other arising from second derivative in $y$ for $\phi$ is
\[
\Upsilon_{2**} = \int_{\mathbb{R}} \Delta^{\Gamma} \phi w_{x} \Bigg[1 - \Big(1 + \ve\big(x
-\sum\limits_{l=0}^{2} \ve^lf_l\big)\Big)^{-2}\Bigg]\mathrm{d}x.
\]
We readily see that
\begin{align}
\begin{aligned}\label{LipshitzforLemdatwo}
\|\Upsilon_{2**}(f_{2}, e)-\Upsilon_{2**}(\tilde{f}_{2},\tilde{e})\|_{L^{q}(\Gamma)}
\,\leq\,
C \ve^{\frac{1}{2}+4+\frac{1}{q}}\big(\|f_{2}-\tilde{f}_{2}\|_a +\|e-\tilde{e}\|_b\big).
\end{aligned}
\end{align}
The remainder $\Upsilon_2 - \Upsilon_{2*}-\Upsilon_{2**}$ actually defines for fixed $\ve$ a compact operator
of the pair $(f_{2},e)$ into $L^{q}(\Gamma)$. This is a consequence of the fact that weak convergence
in $W^{2,q}(\mathfrak{S})$ implies local strong convergence in $W^{1,q}(\mathfrak{S})$. If $f_{2, j}$ and $e_j$
are weakly convergent sequences in $W^{2,q}(\mathfrak{S})$ then clearly the functions $\phi(f_{2,j},e_j)$
constitute a bounded sequence in $W^{1,q}(\mathfrak{S})$. In the above remainder one can integrate by parts
if necessary once in $x$. Averaging against $w_x$ which decays exponentially localizes the situation
and the desired result follows.

Let us consider now the term
\[
\Upsilon_3(z) =\int_{\mathbb{R}} p\big(W^{p-1} - w^{p-1}\big)\phi w_x\mathrm{d}x.
\]
Since the term $W=w(x)+\varepsilon eZ+\sum\limits_{i=2}^{3}\varepsilon^{i}\varphi_{i}+\sum\limits_{i=1}^{3}\phi_i$ can be estimated as
\[
\ve |e(\ve z)Z(x)|+\sum_{i=2}^{3}|\ve^i\varphi_{i}|+\sum_{i=1}^{3}|\phi_{i}| \le C\ve(1
+ |x|^2)e^{-|x|},
\]
we easily see that for some $\kappa > 0$ the uniform bound holds
\[
|W^{p-1}-w^{p-1}|\cdot|w_x| \le C \ve e^{-\kappa|x|}.
\]
From here we readily find that
\[
\|\Upsilon_3\|_{L^{q}(\Gamma)} \,\leq\,
C\ve^{\frac{4}{q}}\|\phi\|_{2,q,\varrho} \,\leq\, C\ve^{\frac{1}{2}+4+\frac{1}{q}}.
\]

We observe also that the term in (\ref{7point6}) such as
$$
\Upsilon_4(z)= \int_{\mathbb{R}}\chi \widetilde{N}_{3}(\phi) w_x\mathrm{d}x,
$$
can be estimated similarly. In fact, using the definition of
$\widetilde{N}_{3}(\phi)$ and the exponential decay
of $w_x$ we obtain
\[
\|\Upsilon_4\|_{L^q(\Gamma)} \leq C
\|\phi\|^2_{2,q,\varrho} \leq C\ve^{9-\frac{6}{q}}.
\]
These terms define compact operators similarly as before.

\subsection{Estimates of the projection against $Z$}
We observe that exactly the same estimates can be carried out in the terms
obtained from integration against $Z$. So the remaining thing is to
compute the term $\int_{\mathbb{R}} E_1 Z \mathrm{d}x$.

Multiplying (\ref{E1}) by $Z$ and integrating over the variable $x$
and using the decomposition of $E_1$ in (\ref{E11}), we get
\begin{align*}
\int_{\mathbb{R}}E_1 Z \mathrm{d}x=\int_{\mathbb{R}} E_{11} Z \mathrm{d}x+\int_{\mathbb{R}} E_{12} Z \mathrm{d}x,
\end{align*}
where
\begin{align*}
\int_{\mathbb{R}} E_{11} Z \mathrm{d}x\,=\,\ve(\ve^2\Delta^{\Gamma}e
+\lambda_0 e)\int_{\mathbb{R}} Z^2\mathrm{d}x
\,=\,\ve^3\Delta^{\Gamma}e+\ve\lambda_0 e.
\end{align*}
On the other hand, we have
\begin{align*}
\int_{\mathbb{R}}E_{12} Z \mathrm{d}x&=\ve^4\int_{\mathbb{R}}\big[\hat{\mathbf{B}}_4
+\hat{G}_{4}+\hat{\mathfrak{D}}_{4}\,+p(p-1)w^{p-2}eZ\phi_{3}\big] Z \mathrm{d}x
\\
&\quad
+\int_{\mathbb{R}} \Big[\sum\limits_{i=4}^{6}\varepsilon^{i}S_{i}+\sum\limits_{i=4}^{6}\varepsilon^{i}T_{i}
+\sum\limits_{i=5}^{10}\varepsilon^{i}G_{i}+\varepsilon^{5}\mathbf{B}_{5}\Big] Z \mathrm{d}x.
\end{align*}
The components in $ \hat{\mathbf{B}}_4,\hat{G}_{4}\,\mbox{and}\,\hat{\mathfrak{D}}_{4}$ are even functions
in the variable $x$ and independent of the parameters $f_{2}$.
Here, the function $\phi_{3}$ has the form
\begin{align*}
\phi_{3}(x,z)=\psi_{31}(x,\var z)+\psi_{32}(x,\var z),
\end{align*}
The components in $\psi_{31}$ and $\psi_{32}$ are independent of the parameters $f_{2}$.
Moreover, $\psi_{31}$ is an odd function in the variable $x$ and $\psi_{32}$ is an even function
in the variable $x$. Therefore, adding up all terms together, we conclude that
\begin{align*}
\int_{\mathbb{R}}E_{12} Z \mathrm{d}x\,=\,&\ve^{5}b^1_{1\ve}\Delta^{\Gamma} e
     +\ve^{5}b^2_{1\ve}\Delta^{\Gamma}f_{2}
+\sum\limits_{i=5}^{10}\varepsilon^{i}b_{i\ve}(f_{2},e,\nabla^{\Gamma} f_{2},\nabla^{\Gamma} e).
\end{align*}
Here we denote by $b^{j}_{1\ve}, j=1, 2$ and $b_{i\ve}$, generic, uniformly bounded continuous functions,
moreover, $b_{i\ve}$ is uniformly Lipschitz in its arguments.

\section{Solving the System for $(f_{2}, e)$: Proof of Theorem \ref{main}}\label{section9}
\setcounter{equation}{0}

\subsection{Proof of Theorem \ref{main}}
\noindent Using the estimates in previous section, we find the following
nonlinear, nonlocal system of differential equations for the parameters $(f_{2}, e)$
in the variable $y=\var z\in\Gamma$
\begin{align}\label{nonlocalsystemoffandeone}
\mathcal{L}_1^{*}(f_{2})&\equiv \Delta^{\Gamma}f_{2}+(k_{1}^{2}+k_{2}^{2})f_{2}
     =\gamma_{1}e+\gamma_{2}\epsilon^{2}\Delta^{\Gamma}e+M_{1\varepsilon}\quad \mbox{in }\Gamma,
\\
\label{nonlocalsystemoffandetwo}
\mathcal{L}_{2\ve}^{*}(e)&\equiv -\ve^2 \Delta^{\Gamma}e-\lambda_0 e
    =M_{2\ve}\quad\quad \mbox{in }\Gamma,\quad\ \
\end{align}
with the boundary conditions
 \begin{align}
 \label{nonlocalsystemoffandethree}
     \frac{\partial f_{2}}{\partial\tau}+If_{2}+\mathfrak{R}&=0\quad \quad\mbox{on }\partial\Gamma,
 \\
 \label{nonlocalsystemoffandefour}
     \frac{\partial e}{\partial\tau}+\frac{1}{2}Ie&=0\quad \quad\mbox{on }\partial\Gamma,
\end{align}
where $\gamma_{1}$ and $\gamma_{2}$ are two constants defined in (\ref{E1timewx}),
$I$ and ${\mathfrak R}$ are
smooth functions defined in (\ref{I}) and (\ref{mathfrakR}) and $\tau$ denotes the inward normal of $\partial\Gamma$.
The operators $M_{1\ve}$ and $M_{2\ve}$ can be decomposed in the following form
$$
M_{l\ve}(f_{2},e)=A_{l\ve}(f_{2},e)+K_{l\ve}(f_{2},e),\ \ l=1,2
$$
where $K_{l\ve}$ is uniformly bounded in $L^q(\Gamma)$ for $(f_{2},e)$
in $\digamma$ and is also compact. The operator $A_{l\ve}$
is Lipschitz in this region, see (\ref{LipshitzforLemdaone})-(\ref{LipshitzforLemdatwo}),
\begin{align}\label{LipshitzforA}
||A_{l\ve}(f_{2},e)-A_{l\ve}(\tilde{f}_{2},\tilde{e})||_{L^q(\Gamma)}
\leq C\var^{\frac{1}{2}+\frac{1}{q}} \bigl[\ ||f_{2}-\tilde{f}_2||_a+||e-\tilde{e}||_b \ \bigr].
\end{align}

Before solving
(\ref{nonlocalsystemoffandeone})-(\ref{nonlocalsystemoffandefour}),
some basic facts about the invertibility of corresponding linear
operators are in order. We first consider the following problem
\begin{align}\label{problemforfwithhomogeneousbd}
\begin{aligned}
 \mathcal{L}_1^{*}(f_{2})=\Delta^{\Gamma}f_{2}+(k_{1}^{2}+k_{2}^{2})f_{2}
     \,=\,h\quad \mbox{in }\Gamma,
\qquad
      \frac{\partial f_{2}}{\partial\tau}+If_{2}\,=\,0\quad \mbox{on }\partial\Gamma.
\end{aligned}
\end{align}

\begin{lem}\label{lemma9point1}
Under the non-degeneracy condition of $\Gamma$ in (\ref{Geometriceigenvalueproblem}),
if $h\in L^q(\Gamma)$ then there is a constant $\ve_0$ for each $0<\ve<\ve_0$,
the problem(\ref{problemforfwithhomogeneousbd}) has a unique
solution $f_{2}\in W^{2,q}(\Gamma)$ with the property
\\
$
\qquad\qquad\qquad\qquad||f_{2}||_{L^{\infty}(\Gamma)}+||\bigtriangledown^{\Gamma}f_{2}||_{L^{\infty}(\Gamma)}
+||\bigtriangleup^{\Gamma}f_{2}||_{L^{q}(\Gamma)}\leq C||h||_{L^{q}(\Gamma)}.
$
\end{lem}

\proof
\noindent Under the non-degeneracy condition of $\Gamma$ in
(\ref{Geometriceigenvalueproblem}), the existence and the a
priori estimates can be easily proved.
\qed
\\

We then consider the following problem
\begin{align}\label{problemforewithhomogeneousbd}
\begin{aligned}
\mathcal{L}_{2\ve}^{*}(e)=-\ve^2 \Delta^{\Gamma}e-\lambda_0 e
\,=\,g\quad \mbox{in }\Gamma,
\qquad
\frac{\partial e}{\partial\tau}+\frac{1}{2}Ie\,=\,0\quad\mbox{on }\partial\Gamma.
\end{aligned}
\end{align}

\begin{lem}\label{lem9point2}
If $g\in L^q(\Gamma)$ then there exists
$\ve_0>0$ such that for a sequence $(\varepsilon_{l})_{l}$ with $0<\ve_l<\ve_0$
such that problem (\ref{problemforewithhomogeneousbd})
has a unique solution $e\in W^{2,q}(\Gamma)$ which satisfies
$$
||e||_{b}\leq C\ \ve^{-2}_{l}\ ||g||_{L^q(\Gamma)}.
$$
Moreover, if $g\in
W^{2,q}(\Gamma)$ then
\begin{align}\label{estimatesofe}
\ve^2_{l}||\Delta^{\Gamma}e||_{L^q(\Gamma)}+\ve_{l}||\nabla^{\Gamma}e||_{L^q(\Gamma)}
+||e||_{L^{\infty}(\Gamma)}\leq C||g||_{W^{2,q}(\Gamma)}.
\end{align}
\end{lem}
The proof will be given in subsection \ref{subsection9.2}.\qed
\\

We finally consider the following system
\begin{align}\label{systemoffande}
\begin{aligned}
\mathcal{L}^*(f_{2},e)\equiv\bigl(\ \mathcal{L}_1^{*}(f),\ \mathcal{L}_{2\ve}^{*}(e)\ \bigr)\,=\, \bigr(\, h,\, g\, \bigr)
\quad\mbox{in }\Gamma
\\
\frac{\partial f_{2}}{\partial\tau}+If_{2}=\Xi,
\qquad
 \frac{\partial e}{\partial\tau}+\frac{1}{2}Ie=0
 \quad\mbox{on }\partial\Gamma,
\end{aligned}
\end{align}
where $\Xi$ is a smooth function.

\begin{lem}\label{lem7point4}
Under the nondegeneracy condition (\ref{Geometriceigenvalueproblem}),
if $h,\ g\in L^q(\Gamma)$ then for the sequence of the parameter $\ve$ in Lemma \ref{lem9point2},
there is a unique
solution $(f_{2},e)$ in $W^{2,q}(\Gamma)$
to problem (\ref{systemoffande}) which satisfies
$$
||f_{2}||_{a}+||e||_{b}
\leq
C \bigl[\ ||h||_{L^q(\Gamma)}+\ve^{-1}||g||_{L^q(\Gamma)}+\|\Xi\|_{L^q(\Gamma)}\ \bigr].
$$
\end{lem}

\proof Under the non-degeneracy condition (\ref{Geometriceigenvalueproblem}), there exist $\bar{f_{2}}$ and $e_0$ satisfying
\begin{align*}
\begin{aligned}
\Delta^{\Gamma}\bar{f_{2}}\, =\, 0,
\qquad
\frac{\partial \bar{f_{2}}}{\partial\tau}+I\bar{f_{2}}\,=\,\Xi.
\end{aligned}
\end{align*}
Setting $f_{2}=\tilde{f_{2}}+\bar{f_0}$  
to the system (\ref{systemoffande}),
the final conclusion can be derived from Lemma \ref{lemma9point1} and Lemma \ref{lem9point2}.
\qed
\\

\medskip
\noindent {\bf Proof of Theorem \ref{main}:}
Let $(\tilde{f},\tilde{e})\in \digamma$, where $\digamma$ is defined in (\ref{regionoffande}),  and define
\begin{align*}
\bigl(\ h(f_{2},e),\ g(f_{2},e)\ \bigr)&=\Bigl(\ \ve A_{1\ve}(f_{2},e)+\ve K_{1\ve}(\tilde{f_{2}},\tilde{e}),\
    \ve^2 A_{2\ve}(f_{2},e)+\ve^2 K_{2\ve}(\tilde{f_{2}},\tilde{e})\ \Bigr),
\\
\Xi&=\mathfrak{R}(\tilde{f_{2}},\tilde{e}).
\end{align*}
From (\ref{LipshitzforA}), $A_{1\ve}$ and $A_{2\ve}$ are contraction mappings of its arguments in $\digamma$.
By Banach Contraction Mapping theorem and Lemma \ref{lem7point4}, we can solve the nonlinear problem
$$
\mathcal{L}^*(f,e)\equiv\Bigl(\mathcal{L}^{*}_1(f_{2}),\mathcal{L}^{*}_{2\ve}(e) \Bigr)=(h,g),
$$
with the boundary conditions defined
in (\ref{systemoffande}) on the region $\digamma$.
Hence, we can define a new operator $\mathcal{Z}$ from $\digamma$ into $\digamma$
by $\mathcal{Z}(\tilde{f_{2}},\tilde{e})=(f_{2},e)$. Finding a solution to the problem (\ref{nonlocalsystemoffandeone})-(\ref{nonlocalsystemoffandefour})
is equivalent to locating a fixed point of $\mathcal{Z}$.
Schauder's fixed point theorem applies to finish the proof of its existence.
Hence, by Proposition \ref{proposition7.1} and the arguments followed, we complete
the existence part of Theorem \ref{main}.
Other properties of $u_{\ve}$ in Theorem \ref{main} can be showed easily.
\qed

\subsection{Proof of Lemma \ref{lem9point2}}\label{subsection9.2}
To prove Lemma \ref{lem9point2}, we follow the method introduced in \cite{dkwy},
which relies only on elementary considerations on the variational
characterization of the eigenvalues of the operator
${\mathcal L}^{*}_{2\ve}$ and the Weyl's asymptotic formula in (\ref{weylonmanifold}).
We remark  this approach is a slightly different from \cite{Mal2}
and \cite{MalMon} where Kato's theorems were the main tools.

First, we consider the following eigenvalue problem
\begin{align}\label{eigenvalueproblem}
\mathcal{L}^{*}_{2\ve}(v)= -\ve^2\Delta^{\Gamma}v-\lambda_0 v=\Lambda_{\ve}v\quad\mbox{in }\Gamma,
\qquad
\frac{\partial v}{\partial\tau}\,+\,\frac{1}{2}\,I\,v\,=\,0\quad\mbox{on }\partial\Gamma.
\end{align}
We denote its eigenvalues $\Lambda_{\ve,j}$ in non-decreasing order and counting them with multiplicity.
Here $\lambda_0$ is the unique positive eigenvalue to the eigenvalue problem (\ref{definitionofZ}),
which implies the spectrum of ${\mathcal L}^{*}_{2\ve}$ contains negative or zero eigenvalues.
From the Courant-Fisher
characterization we can write $\Lambda_{\ve,j}$ in two different ways:
\begin{align}
\Lambda_{\ve,j}&\,=\,\sup_{\Xi\,\in\,\Xi_{j-1}}
\Bigg[\inf_{v\perp \Xi,\,v\neq 0}\frac{\int_{\Gamma}v {\mathcal L}^{*}_{2\ve}v}{\int_{\Gamma}v^2}\Bigg],\label{courant2}
\\
\nonumber
\\
\Lambda_{\ve,j}&\,=\,\inf_{\Xi\,\in\,\Xi_j}
\Bigg[\sup_{v\,\in\, \Xi,\, v\neq 0}\frac{\int_{\Gamma}v {\mathcal L}^{*}_{2\ve}v}{\int_{\Gamma}v^2}\Bigg].\label{courant1}
\end{align}
Here $\Xi_j\,(\mbox{resp. }  \Xi_{j-1})$ represents the family of $j$ dimensional (resp. $j-1$ dimensional)
subspaces of $H^2(\Gamma)$ constituting of functions defined on $\Gamma$
with boundary condition in (\ref{eigenvalueproblem}),
and the symbol $\perp$ denotes orthogonality with respect to the $L^2$ scalar product.
There holds the following result for the estimates of gap between two successive eigenvalues.

\begin{lem}\label{lem7point3}
There exits a number $\ve_0>0$ such that for all $0<\ve_1<\ve_2<\ve_0$ and all $j\geq 1$ the following estimate holds.
\begin{align}
\Lambda_{\ve_1,j}\,=\,\frac{\ve_1^2}{\ve_2^2}\Lambda_{\ve_2,j}-\lambda_0\Big(1-\frac{\ve_1^{2}}{\ve_2^{2}}\Big).
\end{align}
In particular, the functions $\ve\in(0,\ve_0)\mapsto\Lambda_{\ve,j}$ are continuous and increasing.
\end{lem}
\proof
Let us consider small numbers $0<\ve_1<\ve_2$. We observe that
for any $v$ with $\int_{\Gamma}v^2=1$, we have
\begin{align*}
\ve_2^{-2}\int_{\Gamma}v {\mathcal L}^{*}_{2\ve_2}v-\ve_1^{-2}\int_{\Gamma}v {\mathcal L}^{*}_{2\ve_1}v
=\lambda_0(\ve_1^{-2}-\ve_2^{-2}).
\end{align*}
Then the result follows.
\qed

\medskip
\noindent{\bf Proof of Lemma \ref{lem9point2}:}
For $\ell\in\mathbb{N}$, choose $\sigma_\ell=2^{-\ell}$. In order to find a sequence of values $\ve_{\ell}\in(\sigma_{\ell+1}, \sigma_{\ell})$
such that the spectrums of the operators ${\mathcal L}^{*}_{2\ve_{\ell}}$, for large $\ell$,  stay away from $0$, we define
\begin{align*}
\mathcal{K}^1_\ell=\{\,\ve\in(\sigma_{\ell+1},\sigma_\ell):\ \mbox{ker}{\mathcal L}^{*}_{2\ve}\neq \varnothing\,\},
\qquad
\mathcal{K}^2_\ell=(\sigma_{\ell+1},\sigma_\ell)\setminus\mathcal{K}^1_\ell.
\end{align*}
It is crucial to estimate the cardinality of $\mathcal{K}^1_\ell$.
If $\ve\in\mathcal{K}_\ell^1$ then for some $j$ we have that $\Lambda_{\ve,j}=0$.
The monotonicity of the function $\ve\in(0,\ve_0)\mapsto\Lambda_{\ve,j}$ implies that   $\Lambda_{\sigma_{\ell+1},j}<0$.
Hence,
\begin{align}\label{cardinicity}
\mbox{card}(\mathcal{K}^1_\ell)\leq N_{\sigma_{\ell+1}},
\end{align}
where ${\mathbf N}_{\ve}$ is the number of negative eigenvalues of the operator $\mathcal{L}^{*}_{2\ve}$.

We now give an asymptotic estimate on the number ${\mathbf N}_{\ve}$ of negative eigenvalues of the differential operator
$\mathcal{L}^{*}_{2\ve}$.
By $(\rho_i)_i$ we will denote the set of eigenvalues of the eigenvalue problem
$$
-\Delta^{\Gamma}\omega=\rho\,\omega\quad\mbox{in}\,\,\,\Gamma,
\qquad
\frac{\partial \omega}{\partial\tau}\,+\,\frac{1}{2}\,I\,\omega\,=\,0\quad\mbox{on }\partial\Gamma.
$$
From the Weyl asymptotic formula as those in (\ref{weylonmanifold}) and
the formula in (\ref{courant1}), one derives
\begin{align*}
{\mathbf N}_{\ve}\,\geq\,C_{\Gamma}\bigl(1+o(1)\bigr)\ve^{-2},
\end{align*}
where $C_{\Gamma}$ is a fixed constant depending on the volume of the manifold $\Gamma$ and its dimension.
To prove a similar upper bound, we choose $i$ to be the first index such
that $\ve^2\rho_i-\lambda_0>0$.
Then from the Weyl formula we find that
$$
i=C_{\Gamma}\bigl(1+o(1)\bigr)\ve^{-2}.
$$
Define $\Xi_{j-1}=\mbox{span}\{\omega_\ell:\,l=1,2,\cdots,j-1 \}$.
For an arbitrary function $v\in H^2{(\Gamma)}$ and $v\perp \Xi_{j-1}$, we can write
$$
v=\sum_{l\geq j}\kappa_\ell\omega_\ell.
$$
Plugging this $v$ into (\ref{courant2}) and using the Weyl formula, we also have
\begin{align*}
{\mathbf N}_{\ve}\,\leq\,C_{\Gamma}\bigl(1+o(1)\bigr)\ve^{-2}.
\end{align*}
Hence we get that
\begin{align*}
{\mathbf N}_{\ve}\sim C_{\Gamma}\ve^{-2} \ \mbox{ as } \ve\rightarrow 0.
\end{align*}

The last inequality and (\ref{cardinicity}) imply that
$\mbox{card}(\mathcal{K}^1_\ell)\leq C\sigma_\ell^{-2}$, and hence there exists an interval $(a_\ell,b_\ell)$ such that
\begin{align}\label{suitableinterval}
(a_\ell,b_\ell)\subset\mathcal{K}^2_\ell,\quad |b_\ell-a_\ell|\geq\frac{\mbox{meas}(\mathcal{K}^2_\ell)}{\mbox{card}(\mathcal{K}^1_\ell)}
\geq 2C_0\sigma_\ell^{3},
\end{align}
for a universal positive constant $C_0$, independent of $\ell$.
By setting $\ve_\ell=(a_\ell+b_\ell)/2$ for all large $\ell\in\mathbb{N}$, we conclude that ${\mathcal L}^{*}_{2\ve_\ell}$ is invertible
and there exists a number $C>0$, independent of $\ell$, such that for all $j\in\mathbb{N}$ there holds
\begin{align}\label{spectral gap}
\big|\Lambda_{\ve_\ell,j}\big|\,\geq\,C\ve_\ell^{2}.
\end{align}
Assume the opposite, namely that for some $j$ we have
\begin{align*}
\big|\Lambda_{\ve_\ell,j}\big|\,<\,\delta\ve_\ell^{2},
\end{align*}
with $\delta$ arbitrarily small. Since $\ve_{\ell}\in\mathcal{K}^2_\ell$,
then $\big|\Lambda_{\ve_\ell,j}\big|>0$.
Let us assume that
\begin{align}\label{con}
0<\Lambda_{\ve_\ell,j}\,<\,\delta\ve_\ell^{2}.
\end{align}
Then from Lemma \ref{lem7point3}, we have
\begin{align}
\Lambda_{a_\ell,j}
=
\Lambda_{\ve_\ell,j}
-\frac{(\ve_\ell^2-a^2_\ell)}{\ve^2_\ell}(\Lambda_{\ve_\ell,j}+\lambda_0).
\end{align}
The inequalities in (\ref{suitableinterval}) and (\ref{con}) imply that
\begin{align*}
\Lambda_{a_\ell,j}
\leq
\delta\ve_\ell^{2}
-C_0\sigma_{\ell}^{2}\frac{\sigma_{\ell}(\ve_\ell+a_\ell)}{\ve^2_\ell}(\Lambda_{\ve_\ell,j}+\lambda_0)<0,
\end{align*}
if $\delta$ is chosen a priori sufficiently small.
It follows from the continuity of the function $\ve\in(0,\ve_0)\mapsto\Lambda_{\ve,j}$
that $\Lambda_{\ve,j}$ must vanish at some $\ve\in(a_{\ell},\ve_{\ell})$,
and we get a contradiction with the choice of the interval $(a_{\ell},b_{\ell})$.\\
The case
$$
-\delta\sigma_\ell^{2}\,<\,\Lambda_{\ve_\ell,j}<0
$$
can be handled similarly.
In fact, we have the inequality
\begin{align*}
\Lambda_{b_\ell,j}
=
\Lambda_{\ve_\ell,j}
+\frac{(b_\ell^2-\ve_\ell^2)}{\ve^2_\ell}(\Lambda_{\ve_\ell,j}+\lambda_0)>0.
\end{align*}
Hence, the proof of (\ref{spectral gap}) for the spectral gap  between critical eigenvalues was complete.

As a consequence, the solution to (\ref{lem9point2}) exists and satisfies
\begin{align}\label{l2est v}
\|e\|_{L^q(\Gamma)}\,\leq\, C \ve_\ell^{-2}\,{\|g\|_{L^q(\Gamma)}}.
\end{align}
From (\ref{l2est v}) by a standard  elliptic argument one can show
\begin{align*}
{\ve_\ell}^{2}\,\|\Delta^{\Gamma}e\|_{L^{q}(\Gamma)}
+{\ve_\ell}\,\|\nabla^{\Gamma}e\|_{L^{q}(\Gamma)}
+\|{e}\|_{L^\mathbf{\infty}(\Gamma)} \,\leq\,
C\ve_\ell^{-2}\,{\|{h_2}\|_{L^{q}(\Gamma)}}.
\end{align*}
The reader can refer \cite{dkwy} for proof of further estimate in (\ref{estimatesofe}).
\qed
\\

\appendix
\section{ A Linear Model Problem I}
\renewcommand{\thesection}{A}
\setcounter{equation}{0}
\renewcommand{\theequation}{A.\arabic{equation}}


\newtheorem{theorem}{Theorem}[section]
\newtheorem{lemma}{Lemma}[theorem]
\newtheorem{corollary0}{Corollary}[theorem]
\setcounter{theorem}{1}

Recall that $w$ is the even function defined in
(\ref{definitionofw})
 and $Z$ is the even eigenfunction defined in the eigenvalue problem (\ref{definitionofZ}).
Recall that, $\mathbf{\Lambda}$ represents the strip in $\mathbb{R}^3$ with the notations
\begin{align}
\begin{aligned}\label{lambda}
\mathbf{\Lambda}\,=\, \{\, (x,\theta,\eta):\, x\in \mathbb{R}, (\theta,\eta)\in\partial\mathbb{D}\times\mathbb{R}_{+}\, \},
\\
\partial \mathbf{\Lambda}\,=\,\{\, (x,\theta,\eta):\, x\in \mathbb{R},\, \theta\in\partial\mathbb{D},\,\eta=0\, \}.
\end{aligned}
\end{align}

We first consider the following linear problem
\begin{align}\label{problemzerodottwo}
\phi^{0}_{xx}+\frac{1}{l^{2}_{1}(\theta)}\phi^{0}_{\theta\theta}
+\frac{1}{l^{2}_{2}(\theta)}\phi^{0}_{\eta\eta}
-K\phi^0+pw^{p-1}\phi^0=0 \qquad \mbox{in}\,\,\mathbf{\Lambda},
\end{align}
\begin{align}\label{problemzerodotthree}
\phi^0_\eta=G(x,\theta) \qquad \mbox{on} \,\,\partial\mathbf{\Lambda},
\end{align}
\begin{align}\label{problemzerodotfour}
\int_\mathbb{R} \phi^0(x,\theta,\eta)w_x(x)\mathrm{d}x=0,
\qquad
\int_\mathbb{R} \phi^0(x,\theta,\eta)Z(x)\mathrm{d}x=0\quad\,\mbox{in}\,\,\Gamma_{\var},
\end{align}
where $K >\lambda_0 +1$ is a large positive constant,
$l_{1}(\theta)$ and $l_{2}(\theta)$ are two smooth positive functions given in (\ref{I}).
Suppose the following orthogonality conditions hold
\begin{align} \label{orthogonality1}
\int_\mathbb{R} G(x,\theta)w_x(x)\mathrm{d}x=0,
\qquad
\int_\mathbb{R} G(x,\theta)Z(x)\mathrm{d}x=0.
\end{align}

\begin{lemma}\label{model0.5}
If $G\in L^2(\partial\mathbf{\Lambda})$, and the orthogonality
conditions (\ref{orthogonality1}) hold, then there is a unique
solution $\phi^0$ to the problem (\ref{problemzerodottwo})-(\ref{problemzerodotfour})
for any large positive constant $K$. Moreover there is a
constant $C>0$, independent of $\ve$, such that the solution to the problem
(\ref{problemzerodottwo})-(\ref{problemzerodotfour})
 satisfies a priori estimate
$$
||\phi^0||_{H^2(\mathbf{\Lambda})}\leq C||G||_{L^2(\partial\mathbf{\Lambda})}.
$$
\end{lemma}

\proof Since $K$ is large,
the proof of  the existence and uniqueness of the solution to (\ref{problemzerodottwo})-(\ref{problemzerodotfour})
and its estimate is standard.
To show the $L^2$-orthogonality (\ref{problemzerodotfour}),
using the equations of $Z(x)$ and $\phi^0$ and also the condition (\ref{orthogonality1}),
for
$$
\varphi(\theta,\eta)=\intreal \phi^0(x,\theta,\eta)Z(x)\mathrm{d}x,
$$
one finds
$$
\frac{1}{l^{2}_{1}(\theta)}\varphi_{\theta\theta}
+\frac{1}{l^{2}_{2}(\theta)}\varphi_{\eta\eta}-(K-1-\lambda_0)\,\varphi=0 \quad\,\mbox{in}\,\,\Gamma_{\var},
\quad \varphi_{\eta}(\theta,0)=0.
$$
Choosing  $K>\lambda_0 +1$,  we deduce that
\begin{align*}
\varphi(\theta,\eta)=\intreal \phi^0(x,\theta,\eta)Z(x)\mathrm{d}x=0\quad\,\mbox{in}\,\,\Gamma_{\var}.
\end{align*}
Similarly we have
\begin{align*}
\intreal \phi^0(x,\theta,\eta)w_x(x)\mathrm{d}x=0\quad\,\mbox{in}\,\,\Gamma_{\var}.
\end{align*}
\qed

A special case of Lemma \ref{model0.5} is the following problem:
finding function $\hat{\phi}\in H^2\bigl(\mathbf{\Lambda}\bigr)$ such that
\begin{align}\label{problemzerodotsix}
\begin{aligned}
\hat{\phi}_{xx}+\frac{1}{l^{2}_{1}(\theta)}\hat{\phi}_{\theta\theta}+\frac{1}{l^{2}_{2}(\theta)}\hat{\phi}_{\eta\eta}&-\tilde{K}\hat{\phi}+pw^{p-1}\hat{\phi}=0 \qquad \mbox{in}\,\,\mathbf{\Lambda},
\\
&\hat{\phi}_\eta=G(x,\theta)  \qquad  \mbox{on} \,\,\partial\mathbf{\Lambda},
\end{aligned}
\end{align}
where $\tilde{K}$ is a large positive constant.
\begin{lemma}\label{lemma0.6}
Suppose the function $\,G(x,\theta)$ is even in the variable $x$,
then there exists a large positive constant $\tilde{K}$  such that the problem
(\ref{problemzerodotsix}) has a unique
solution $\hat{\phi}$, which is an even function in the variable $x$ and satisfies
$$||\hat{\phi}||_{H^2(\mathbf{\Lambda})}\leq C ||G||_{L^2(\partial\mathbf{\Lambda})},$$
Moreover, if $\,G(x,\theta)$ is exponentially decaying in $x$, then
\begin{align}\label{zerodotseven}
\bigl | \hat{\phi}(x,\theta,\eta)\bigr |< C e^{-\alpha|x|},
\end{align}
where $\alpha>0$ and
 the constant $C$ does not depend on $\ve$.
\end{lemma}
\proof
For any ${\tilde K}$ large enough,
the proof of  the existence and uniqueness of the solution to (\ref{problemzerodotsix})
and its estimate is standard.
By uniqueness and evenness of $\,G(x,\theta)$,
$\hat{\phi}$ is an even function in the variable $x$.  By the exponentially decaying of $\,G(x,\theta)$, we also have (\ref{zerodotseven}).
\qed

Next, we consider the following problem
\begin{align}\label{problemzerodoteight}
\begin{aligned}
\mathcal{L}_0(\tilde{\phi})\equiv\tilde{\phi}_{xx}+\frac{1}{l^{2}_{1}(\theta)}\tilde{\phi}_{\theta\theta}
+&\frac{1}{l^{2}_{2}(\theta)}\tilde{\phi}_{\eta\eta}
-\tilde{\phi}+pw^{p-1}\tilde{\phi}=h \qquad  \mbox{in}\,\,\mathbf{\Lambda},
\\
&\tilde{\phi}_\eta=G(x,\theta)      \qquad  \mbox{on} \,\,\partial\mathbf{\Lambda},
\\
\intreal\tilde{\phi}(x,\theta,\eta)w_x(x)\mathrm{d}x=0,
\quad
&\intreal\tilde{\phi}(x,\theta,\eta)Z(x)\mathrm{d}x=0\quad\,\mbox{in}\,\,\Gamma_{\var}.
\end{aligned}
\end{align}
\begin{lemma}\label{lemma0.7}
If $h\in L^2{(\mathbf{\Lambda})},\ G\in L^2(\partial\mathbf{\Lambda})$ and
the orthogonality conditions (\ref{orthogonality1}) hold, then for any solution $\tilde{\phi}$ to
problem (\ref{problemzerodoteight}) we have
$$||\tilde{\phi}||_{H^2(\mathbf{\Lambda})}\leq C\bigl[\ ||h||_{L^2(\mathbf{\Lambda})}
+||G||_{L^2(\partial\mathbf{\Lambda})}\ \bigr]$$
where the constant $C$ does not depend on $h,\ G$ and $\ve$.
Furthermore, if $ |h| + |G| \leq C e^{ -\alpha |x|}$,
then $ |\tilde{\phi}| \leq C e^{-c \alpha |x|}$ for some $C,\, c>0$.
\end{lemma}

\proof Let $\phi^0(x,\theta,\eta)$ be defined in Lemma \ref{model0.5}
 and
$\tilde{\phi}=\phi^0+\phi$.
Then we  have
\begin{align}\label{zerodotnine}
\begin{aligned}
\phi_{xx}+\frac{1}{l^{2}_{1}(\theta)}\phi_{\theta\theta}+\frac{1}{l^{2}_{2}(\theta)}\phi_{\eta\eta}&-\phi
+pw^{p-1}\phi=h+(1-K)\phi^0 \qquad  \mbox{in}\,\,\mathbf{\Lambda},
\\
&\phi_\eta=0   \qquad  \mbox{on} \,\,\partial\mathbf{\Lambda},
\\
\intreal\phi(x,\theta,\eta)w_x(x)\mathrm{d}x&=0,
\qquad
\intreal\phi(x,\theta,\eta)Z(x)\mathrm{d}x=0\quad\,\mbox{in}\,\,\Gamma_{\var}.
\end{aligned}
\end{align}
\noindent
Let $\xi_k,k=1,2,\cdots,$ be the eigenfunctions (corresponding eigenvalues $\lambda_k,k=1,2,\cdots$) of
the following eigenvalue problem
$$
-\frac{1}{l^{2}_{1}(\theta)}\phi_{\theta\theta}-\frac{1}{l^{2}_{2}(\theta)}\phi_{\rho\rho}
=\lambda\xi\quad\,\mbox{in}\,\,\Gamma.
$$
Let us consider Fourier series
decompositions for $h$ and $\phi$ of the form
Let us consider Fourier series decompositions for $h+(1-K)\phi^0$ and $\phi$:
$$\phi(x,\theta,\eta)=\sum^{\infty}_{k=0}\phi_k(x)\xi_k(\var\theta,\var\eta),$$
$$h(x,\theta,\eta)+(1-K)\phi^0(x,\theta,\eta)=\sum^{\infty}_{k=0}h_k(x)\xi_k(\var\theta,\var\eta).$$

From the equation (\ref{zerodotnine}) we arrive at the following equation
\bb\label{zerodotten}
-\ve^2\lambda_k \phi_{k}+\phi_{k,xx}-\phi_k+pw^{p-1}\phi_k=h_k,
\ee
with the orthogonality condition
\begin{align}\label{zerodoteleven}
\intreal\phi_k(x) w_x(x)\mathrm{d}x=0,\ \ \ \ \intreal\phi_k(x) Z(x)\mathrm{d}x=0.
\end{align}
Let us consider the bilinear form in $H^1(\mathbb{R})$
$$
B(\psi, \psi) = \int_{\mathbb{R}} [\ |\psi_x|^2 + |\psi|^2-pw^{p-1}|\psi|^2 \ ]\, \mathrm{d}x \ .
$$
Since (\ref{zerodoteleven}) holds uniformly in $k$ we conclude that
\begin{equation}\nonumber
C[\ \|\phi_{k}\|_{L^2(\mathbb{R})}^2 + \|\phi_{k,x}\|^2_{L^2(\mathbb{R})}\ ] \, \le\,
B(\phi_{k}, \phi_{k}),
\end{equation}
for a constant $C>0$ independent of $k$. Using this fact and  equation (\ref{zerodotten}) we arrive at
\begin{align}\label{zerodottwelve}
(1+\lambda^{2}_{k}\ve^4)\|\phi_{k}\|_{L^2(\mathbb{R})}^2 + \|\phi_{k,x}\|^2_{L^2(\mathbb{R})} \,
\le\, C\|h_{k}\|_{L^2(\mathbb{R})}^2.
\end{align}
Moreover, we see from (\ref{zerodotten}) that $\phi_{k}$ satisfies an
equation of the form
$$
 \phi_{k,xx}  -  \phi_{k} = \tilde h_{k} \qquad  \mbox{on} \,\, \mathbb{R},
$$
where $\|\tilde h_{k}\|_{L^2(\mathbb{R})} \le C \|h_{k}\|_{L^2(\mathbb{R})}$. Hence
it follows that
\begin{equation}\label{zerodotthirteen}
\|\phi_{k,xx}\|^2_{L^2(\mathbb{R})}  \, \le\, C\|h_{k}\|_{L^2(\mathbb{R})}^2.
\end{equation}
Summing up estimates (\ref{zerodottwelve}) and (\ref{zerodotthirteen}) in
$k$, we conclude that
$$
\|D^2 \phi\|^2_{L^2(\mathbf{\Lambda})} +\|D \phi\|^2_{L^2(\mathbf{\Lambda})}
+\|\phi\|^2_{L^2(\mathbf{\Lambda})} \, \le\, C\|h\|_{L^2(\mathbf{\Lambda})}^2.
$$
The final estimate follows from the estimates of $\phi$ and $\phi^0$.
\qed

A corollary of Lemma \ref{lemma0.7} is the following
\begin{corollary0}\label{cor1}
Let $ G\in L^2(\partial\mathbf{\Lambda})$  satisfy
the orthogonality conditions (\ref{orthogonality1}) and $ h=0$.
Then problem (\ref{problemzerodoteight}) has a unique solution $ \tilde{\phi}$ such that
\begin{align*}
\|\tilde{\phi}\|_{H^2({\mathbf{\Lambda}})} \leq C \|G\|_{L^2(\partial{\mathbf{\Lambda}})},
\end{align*}
where the constant $C$ does not depend on $G$ and $\varepsilon$.
Furthermore, if $  |G| \leq C e^{ -\alpha |x|}$,
then $ |\tilde{\phi}| \leq C e^{-c \alpha |x|}$ for some $C, c>0$.
\qed
\end{corollary0}


\section*{Appendix B. A Linear Model Problem II}
\renewcommand{\thesection}{B}
\setcounter{equation}{0}
\renewcommand{\theequation}{B.\arabic{equation}}
\newtheorem{lemma1}{Lemma}[theorem]
\setcounter{theorem}{2}

Recall $\mathfrak{S}$ represents the strip
\begin{align}\label{mathfraks}
 \{\ (x,z): x\in \mathbb{R}, z\in\Gamma_{\var} \}
\end{align}
in $\mathbb{R}^3$. $\partial \mathfrak{S}$ is the component of the boundary of $\mathfrak{S}$,
i.e.
$$\partial \mathfrak{S}=\{\ (x,z): x\in \mathbb{R},\, z\in\partial\Gamma_{\var} \}.$$

We consider the following problem: given $h\in L^q{\bigl(\mathfrak{S}\bigr)}$
and $G\in L^q(\partial\mathfrak{S})$, finding
functions $\ \phi\in W^{2,q}\bigl(\mathfrak{S}\bigr),\ c, d \in L^q(\Gamma)$
and $\Lambda_{1}, \Lambda_{2}$ such that
\begin{align}\label{problem0.27}
\begin{aligned}
\phi_{xx}+\Delta^{\Gamma_{\var}}\phi-\phi&+pw^{p-1}\phi
     =h+c(\ve \eta)\, \chi(\ve |x|)\, w_x+d(\ve \eta)\, \chi(\ve |x|)\, Z \quad  \mbox{in}\,\,\mathfrak{S},
\\
\frac{\partial\phi}{\partial\tau_{\var}}&=G\qquad  \mbox{on} \,\, \partial\mathfrak{S},
\\
\intreal\phi(x,z)&w_x\mathrm{d}x=\Lambda_{1},
\quad
\intreal\phi(x,z)Z(x)\mathrm{d}x=\Lambda_{2}
\quad\,\mbox{in}\,\,\Gamma_{\var},
\end{aligned}
\end{align}
where $\tau_{\var}$ denotes the inward normal of $\partial\mathfrak{S}$, and $\chi(t)$ a  smooth  cut-off function such that $ \chi (t)=1$ for $ |t| \leq 10 \sigma$ and $ \chi (t)=0$
for $t \geq 20 \sigma$, and $\sigma>0$ is a small constant defined in Section \ref{section5}.

\begin{lemma}\label{lemma0.12}
There exist functions $\ c(\ve z),\ d(\ve z)$ with respect to $h$ such that the problem (\ref{problem0.27}) has
a unique solution $\phi=T_1(h,\, G)$.
Moreover,
\begin{align*}
||\phi||_{q,\varrho}&\,\leq\, C(||h||_{L^q(\mathfrak{S})} + || g ||_{L^q (\partial\mathfrak{S})}),
\\
||\Lambda_{i}||_{L^q(\Gamma_{\var})}&\,\leq\, C(||h||_{L^q(\mathfrak{S})} + || g ||_{L^q (\partial\mathfrak{S})}),
\quad \forall\, i=1,2,
\end{align*}
where the constant $C$ does not depend on $h,\, G$ and $\ve$.
\end{lemma}

\proof The proof will be carried out in three steps.
\\
{\bf Step 1:} Let us assume that in problem (\ref{problem0.27}) the terms $G,c(\var z),d(\var z)$ are identically zero. Arguing as in Lemma \ref{lemma0.7}, for any $h\in L^{q}(\mathfrak{S})$ and any solution $\phi\in W^{2,q}(\mathfrak{S})$ of problem (\ref{problem0.27}) we have
\begin{align}\label{phi}
||\phi||_{q,\varrho}\leq C\, ||h||_{L^q(\mathfrak{S})}
\end{align}

{\bf Step 2:} We claim that the a priori estimate obtained in Step 1
is in reality valid for the full problem (\ref{problem0.27}).
We first choose suitable $\Lambda_{1}$ and $\Lambda_{2}$ such that
\begin{align}
\nabla^{\Gamma_\ve}\Lambda_{1}=\intreal g(x,z)w_x\,\mathrm{d}x,
\qquad
\nabla^{\Gamma_\ve}\Lambda_{2}=\intreal g(x,z)Z(x)\,\mathrm{d}x
\quad\,\mbox{in}\,\,\Gamma_{\var}.
\end{align}
Let $\phi^{0}$ be the solution of
\begin{align*}
\phi^{0}_{xx}+\Delta^{\Gamma_{\var}}\phi^{0}-\phi^{0}+pw^{p-1}\phi^{0}=0 \quad  \mbox{in}\,\,\mathfrak{S},
\qquad
\frac{\partial\phi^{0}}{\partial\tau_{\var}}=G\qquad  \mbox{on} \,\, \partial\mathfrak{S},
\end{align*}
Note that we have
\begin{align}\label{phi0}
||\phi^{0}||_{q,\varrho}\leq C\,||G||_{L^q(\partial\mathfrak{S})}
\end{align}
By defining
$$
\bar{\Lambda}_{1}=\intreal\phi^{0}(x,z)w_x\mathrm{d}x,
\qquad
\bar{\Lambda}_{2}=\intreal\phi^{0}(x,z)Z(x)\mathrm{d}x
\quad\,\mbox{in}\,\,\Gamma_{\var},
$$
we have that, to prove the general case it suffices to apply the argument with
$$
\bar{\phi}=\phi-\phi^{0}+\frac{(\bar{\Lambda}_{1}-\Lambda_{1})w_{x}}{\intreal w^{2}_x\mathrm{d}x}
+\frac{(\bar{\Lambda}_{2}-\Lambda_{2})Z(x)}{\intreal Z^{2}(x)\mathrm{d}x}.
$$
The $\bar{\phi}$ satisfies a problem of a similar form with homogeneous Neumann boundary condition and orthogonality condition, as well as $h$ replaced by a function $\bar{h}$ with norm bounded by
$$
||\bar{h}||_{L^q(\mathfrak{S})}\leq C\,( ||h||_{L^q(\mathfrak{S})}
+||G||_{L^q(\partial\mathfrak{S})}).
$$

{\bf Step 3:} We consider the problem
\begin{align*}
\bar{\phi}_{xx}+\Delta^{\Gamma_{\var}}\bar{\phi}&-\bar{\phi}+pw^{p-1}\bar{\phi}=\bar{h}+c\,\chi\,w_x+d\,\chi\,Z\qquad  \mbox{in}\,\,\mathfrak{S},
\\
&\frac{\partial\bar{\phi}}{\partial\tau_{\var}}=0\qquad  \mbox{on} \,\, \partial\mathfrak{S},
\\
 \intreal\bar{\phi}(x,z)&w_x\mathrm{d}x=0,\ \ \ \intreal\bar{\phi}(x,z)Z(x)\mathrm{d}x=0\quad\,\mbox{in}\,\,\Gamma_{\var}.
\end{align*}
The existence of the solutions can be proved similarly as that in Lemma \ref{lemma0.7}.
There also hold the priori estimate
\begin{align*}
||\phi||_{q,\varrho}&\leq C(||h||_{L^q(\mathfrak{S})} + || g ||_{L^q (\partial\mathfrak{S})}),
\\
||\Lambda_{i}||_{L^q(\Gamma_{\var})}&\leq C(||h||_{L^q(\mathfrak{S})} + || g ||_{L^q (\partial\mathfrak{S})}),\quad \forall i=1,2.
\end{align*}
\qed


\end{document}